%% file: Skyrme-final.tex
\documentclass{amsart}
\usepackage{amssymb,amsmath,amsthm,mathtools, geometry, amssymb, verbatim, esint, soul}
\usepackage{a4wide}
\usepackage{float}
\usepackage{graphicx}
\usepackage[utf8]{inputenc} 
\usepackage{a4wide}
\usepackage{tikz}
\usetikzlibrary{calc}
\usepackage{ulem}
\usepackage{stackengine,scalerel,graphicx}
\stackMath
\usepackage{hyperref}
\usepackage{amsfonts} 
\usepackage{latexsym}
\usepackage[font=small,format=hang,labelfont={sf,bf}]{caption}
\usepackage{epsfig}
\usepackage{subfig}
\usepackage{url}
\usepackage{varioref}
\usepackage{bm}
\usepackage{color}
\usepackage{mathrsfs}
\usepackage{latexsym}
\usepackage{bbm}
\usepackage{faktor}
\usepackage{yhmath}
\usepackage{empheq}
\usepackage{mathtools}
\usepackage{color}
\usepackage{marginnote}
\usepackage{comment}
\usepackage{cancel}
\usetikzlibrary{arrows.meta} 
\newcommand\overarc[1]{\ThisStyle{%
  \setbox0=\hbox{$\SavedStyle#1$}%
  \stackon[.5pt]{\SavedStyle#1}{%
  \rotatebox{90}{$\SavedStyle\scaleto{)}{.95\wd0}$}}}}

\allowdisplaybreaks[1]

\newcommand{\diam}{\mathrm{diam}\,}          

\newcommand{\mres}{\mathbin{\vrule height 1.6ex depth 0pt width
0.13ex\vrule height 0.13ex depth 0pt width 1.3ex}}

\newcommand{\calT}{\mathcal{T}}
\newcommand{\calL}{\mathcal{L}}
\newcommand{\calN}{\mathcal{N}}
\newcommand{\calU}{\mathcal{U}}
\newcommand{\calR}{\mathcal{R}}
\newcommand{\calE}{\mathcal{E}}

\newcommand{\calM}{\mathcal{M}}

\newcommand{\calC}{\mathcal{C}}

\newcommand{\calSF}{\mathcal{SF}}

\newcommand{\boldn}{n}

\newcommand{\mathN}{\mathbb{N}}

\newcommand{\mathR}{\mathbb{R}}
\newcommand{\mathZ}{\mathbb{Z}}
\newcommand{\mathS}{\mathbb{S}}

\newcommand{\cone}{\mathrm{cone}}

\newcommand{\eps}{\varepsilon}
\newcommand{\wto}{\overset{\ast}{\rightharpoonup} }

\makeatletter

\@addtoreset{equation}{section}
\makeatother

\newtheorem{theorem}{Theorem}[section]
\newtheorem*{theorem*}{Theorem}
\newtheorem{lemma}[theorem]{Lemma}

\theoremstyle{definition}
\newtheorem{definition}[theorem]{Definition}
\newtheorem{remark}[theorem]{Remark}

\title[concentration in magnetic skyrmions]{Energy concentration in a two-dimensional magnetic skyrmion model: variational analysis of lattice and continuum theories}

\author[L. Briani]{L. Briani}
\address[Luca Briani]{Technische Universit\"at M\"unchen, Boltzmannstrasse 3, 85748 Garching, Germany	}
\email[]{luca.briani@tum.de}

\author[M. Cicalese]{M. Cicalese}
\address[Marco Cicalese]{Technische Universit\"at M\"unchen, Boltzmannstrasse 3, 85748 Garching, Germany	}
\email[]{cicalese@.ma.tum.de}

\author[L. Kreutz]{L. Kreutz}
\address[Leonard Kreutz]{Technische Universit\"at M\"unchen, Boltzmannstrasse 3, 85748 Garching, Germany	}
\email[]{kleo@cit.tum.de}

\date{\today}	

\begin{document}

\begin{abstract}
We investigate the formation of singularities in a baby Skyrme type energy model, which describes magnetic solitons in two-dimensional ferromagnetic systems. In presence of a diverging anisotropy term, which enforces a preferred background state of the magnetization, we establish a weak compactness of its topological charge density, which converges to an atomic measure with quantized weights. We characterize the $\Gamma$-limit of the energies as the total variation of this measure. In the case of lattice type energies, we first need to carefully define a notion of discrete topological charge for $\mathbb{S}^2$-valued maps. We then prove a corresponding compactness and $\Gamma$-convergence result, thereby bridging the discrete and continuum theories.
\end{abstract}
	
\maketitle


\tableofcontents

\section{Introduction}

Magnetic solitons are localized, stable solutions to nonlinear equations governing magnetic systems which play a central role in materials science. These solitons, such as vortices, domain walls, and chiral magnetic skyrmions, are characterized by a nontrivial topology that prevents them from being continuously deformed into a uniform state (a trivial ground state of all such models). Very recent real space observations of these structures on the plane (see for instance \cite{Heinze2011,Muhlbauer2009, Zu2010}) have generated considerable theoretical interest. In two dimensional models for magnetic solitons, one typically considers a magnetization field $u:\Omega\subset\mathbb{R}^2\to\mathbb{S}^2$, where $\Omega\subset\mathbb{R}^2$ is open and bounded and $\mathbb{S}^2$ denotes the unit sphere in 
$\mathbb{R}^3$. An energy functional is then associated to $u$ to drive the system toward equilibrium. The choice of the energy model is often made in order to single out one or more physical phenomenon one aims to focus on. In particular energy concentration on topological singularities can be effectively captured by the so-called baby Skyrme model. This model encapsulates the interplay between energy minimization, localization, and topological properties of solitons, and provides insight into the behavior of more complex systems where additional interactions may be present. Its energy can be written as 
\begin{equation}\label{intro:S}
F^{BS}(u)=\frac{1}{2}\int_{\Omega}|\nabla u|^2\,\mathrm{d}x+\kappa \int_\Omega u\cdot(\nabla\times u)\, \mathrm{d}x+\lambda\int_{\Omega}|u-n|^2\,\mathrm{d}x\,,
\end{equation}
where \(n\in\mathbb{S}^2\) is a fixed unit vector representing the background state and $\kappa$ and $\lambda$ are two positive constants. The Dirichlet energy of $u$ physically models the so called exchange interaction, which favors a uniform (or slowly varying) magnetization. The second term, known as the Dzyaloshinskii-Moriya interaction (DMI), favors twisting of the magnetic texture and stabilizes chiral structures in magnets. The third term, associated with the Zeeman interaction, penalizes deviations of $u$ from the uniform state $n$. In this context, a soliton is characterized by a topological charge defined as the Brouwer degree of the map $u\circ\sigma$, where $\sigma$ is the stereographic projection mapping $n$ to $\infty$, and given by
\begin{equation*}
\deg(u)=\frac{1}{4\pi}\int_{\Omega} u\cdot\bigl(\partial_{x}u\times \partial_{y}u\bigr)\,\mathrm{d}x\,.
\end{equation*}
This integer corresponds to the number of times (with sign determined by the orientation) the domain $\Omega$ is "wrapped" around the target sphere $\mathbb{S}^2$. The quantization of the topological charge translates the physic idea of topological protection: solitons cannot simply "unwind" without encountering a discontinuity, which would cost an infinite amount of energy. \\

Topologically non-trivial energy minimizers of \eqref{intro:S} have been recently associated to local minimizers of the micromagnetic energy \cite{Bogdanov1994,Bogdanov1999,Bogdanov1989,BogdanovJETP1989},  and as such they have been recognized as attractive candidates for information technology applications (see for instance \cite{Kiselev2011,Zhang2020}). From a mathematical perspective, instead, the appeal of models like the baby Skyrme one described above lies in their balance of simplicity and richness: while the energy functionals are relatively straightforward, their minimizers exhibit intricate patterns and nontrivial topological features. 
As a result, over the past decades, much effort has been devoted to understanding both the existence and the qualitative properties of solitons in such models, employing tools from the calculus of variations and geometric measure theory. In this framework the first rigorous mathematical studies of chiral magnetic skyrmions driven by an energy of the type \eqref{intro:S} goes back to \cite{Mel} (see also \cite{Lin2004}) in which the existence of degree-one global energy minimizers in the whole plane is established. However, a comprehensive mathematical characterization of solitons remains incomplete. Recently some relevant results concerning the existence (under confinement or in specific topological classes) or the optimal profile of solitons have been obtained (see for instance \cite{Doring2017, Komineas2020, Komineas2021, Bernand2021, Gustafson2021, Bernand2022, Monteil2023, DiFratta2024, Muratov2024}). 
Continuum models for magnetic solitons, like the baby Skyrme model, not only provide a fertile ground for new mathematical developments and have significant practical applications, but they also serve as crucial benchmarks for numerical simulations. The latter have already been very helpful in  uncovering a vast array of possible local minimizers (see \cite{Foster2019,Kuchkin2020,Kuchkin2023}). However, a mathematically rigorous transition from theoretical models to numerical applications necessitates discretizing the problem and thus passing to a lattice formulation. The discrete framework is also interesting in its own right, as it offers insight into the energetic behavior of classical lattice spin systems and their variational coarse graining that occurs as the lattice spacing vanishes (see for instance \cite{ABCS23} and the reference therein for a recent account of this topic when the energy concentrates at the surface scaling). In both cases, ensuring that the discrete formulation accurately captures the energetic structure and topological invariants of the continuum theory is a challenging task. A particularly delicate issue in the discrete setting is the definition of the topological charge. For $\mathbb{S}^1$-valued maps this has been successfully done in several recent papers concerning the variational analysis of classical spin systems like the ferromagnetic and antiferromagnetic  XY models (\cite{AC09, ADGP14, BCKO22, CanSeg, COR22CPAM, COR22ARMA, BCDP18}), simple models for nematic liquid crystals \cite{BraCicSol} and even models describing the emergence of dislocations in microplasticity (\cite{BCGO24, P07}). In all these cases one can take advantage of the well-understood description of the energy concentration phenomenon in the Ginzburg-Landau setting (see for instance \cite{Bethuel1994,SS}) and relate the computation of the degree of the relevant order parameter of each model to the winding number of a field defined on a lattice. This discrete winding number is then obtained by summing the oriented phase differences along the lattice edges, an approach that leverages the natural cyclic order of $\mathbb{S}^1$. By contrast $\mathbb{S}^2$ lacks such a canonical ordering, making such a task harder. In this paper, our construction aims at assigning to each lattice triangle a signed area and defines the degree as the sum of these contributions. However, such an assignment becomes ambiguous in the case that neighbouring lattice points are mapped into antipodal vectors on the two-dimensional sphere. We point out that such an ambiguity is the result of the finite resolution of the reference lattice and that more work is required to take care of such ambiguities. To make this task even more difficult, our discrete topological charge of a lattice spin field $u$ needs to match the degree of its continuous interpolation leaving both energies (of the discrete and of the continuum field) essentially unchanged. To this end, we make use of a refined interpolation constructed on a finer lattice. The definition of discrete topological charge is discussed in Section~\ref{sec:dis1} and is one of the main contribution of this paper. To the best of our knowledge, such a general construction appears here for the first time (see \cite{berg1981} for a definition when non ambiguites occur).\\

In the rest of the introduction we summarise the main results of the paper. We are interested in the onset of topological singularities and in their description through energy concentration. More precisely, we study baby skyrmion-type energies in both continuum and discrete settings. In the continuum formulation, for a bounded open set \(\Omega\subset\mathbb{R}^2\) and $n\in\mathbb{S}^2$, we consider the space of admissible spins
\[
\calU(\Omega)=\{u\in H^1(\mathbb{R}^2;\mathbb{S}^2) : u(x)=n \text{ for a.e. } x\in\mathbb{R}^2\setminus\Omega\}\,,
\]
and define the energy obtained, for $\eps>0$, from \eqref{intro:S} when $\kappa=0$ and $\lambda=1/\eps^2$, namely
\begin{equation}\label{intro:Feps}
F_\varepsilon(u,\Omega)=\frac{1}{2}\int_{\Omega}|\nabla u|^2\,\mathrm{d}x+\frac{1}{\varepsilon^2}\int_{\Omega}|u-n|^2\,\mathrm{d}x\,.
\end{equation}
As explained in Section \ref{sec:Skyrme}, one can incorporate in the analysis below also the case $\kappa_\eps=\kappa\neq 0$, since the additional DMI term $\int_\Omega u\cdot(\nabla\times u)\,\mathrm{d}x$ is a continuum and vanishing perturbation in the topology chosen to perform the $\Gamma$-limit (for a detailed presentation of this notion we refer to the monographs \cite{Bra,DM}). Our first main result is the compactness Theorem~\ref{th:contcomp} asserting that if a sequence \((u_\varepsilon)\subset\calU(\Omega)\) satisfies
$\sup_{\varepsilon}F_\varepsilon(u_\varepsilon,\Omega)<+\infty$,
then the associated topological charge densities
\[
q(u_\varepsilon)=u_\varepsilon\cdot(\partial_{x}u_\varepsilon\times \partial_{y}u_\varepsilon)
\]
converge (in the weak$^*$ sense of measures) to an atomic measure \(\mu\in\calM(\Omega)\) such that \(\mu/(4\pi)\) has integer coefficients. This result suggests that to track energy concentration one can define an energy which depends on the concentration measure itself, namely $\mathscr{F}_\eps:\calM(\Omega)\to [0,+\infty]$:
\begin{equation*}
\mathscr{F}_\eps(\mu)=\inf\left\{F_\eps(u,\Omega): u\in\calU(\Omega),\, q( u_\eps)\,\mathrm{d}x=\mu\right\}\,.
\end{equation*}
In Theorem ~\ref{th:contmain} we prove that, as $\eps\to 0$, $\mathscr{F}_\eps$ $\Gamma\mbox{-}$converge with respect to the weak$^*$ topology of $\calM(\Omega)$ to $\mathscr{F}:\calM(\Omega)\to [0,+\infty]$ which is finite on those atomic measures $\mu=4\pi\sum_{i=1}^{N}d_i\delta_{x_i},\ d_i\in\mathbb{Z}$, on which it takes the value 
\begin{equation*}
\mathscr{F}(\mu)=|\mu|(\Omega)\,.
\end{equation*}
Theorem ~\ref{th:contmain} can be seen as the extension to low energy configurations of one of the results proved for minimizers in \cite[Theorem 2.5]{Muratov2024}.\\

In the discrete setting, the analysis becomes considerably more intricate. We consider a regular triangular lattice $\calL_\varepsilon\subset\mathbb{R}^2$ of spacing $\eps$ and define the set of admissible spin fields 
\[
\mathcal{SF}_\eps(\Omega)=\{u:\calL_\eps\to \mathS^2:\ u_\eps\equiv n \text{ in }\calL_\eps(\mathR^2\setminus\Omega)\}\,.
\]
For \(u\in\mathcal{SF}_\eps(\Omega)\), we introduce its discrete energy $E_\eps(u,\Omega)$ given by
\[
E_\varepsilon(u,\Omega)=\sum_{T\in\calT_\varepsilon(\Omega)}\Bigl(H_\varepsilon(u,T)+Z_\varepsilon(u,T)\Bigr)\,,
\]
where the sum is extended to all elementary triangular cells $T:=[i,j,k]\subset\Omega$ of the lattice $\calL_\eps$, and where
\[
H_\varepsilon(u,T)=\frac{1}{2}\Bigl(|u(i)-u(j)|^2+|u(j)-u(k)|^2+|u(k)-u(i)|^2\Bigr)\,,
\]
and
\[
Z_\varepsilon(u,T)=\frac{1}{2}\Bigl(|u(i)-n|^2+|u(j)-n|^2+|u(k)-n|^2\Bigr)
\]
are the exchange and anisotropy discrete energies, respectively. If on one side the form of $H_\eps$ and $Z_\eps$ is suggested by the straightforward discretization of the corresponding continuum energies in \eqref{intro:Feps}, on the other side the function space where we want to consider the discrete energies requires a nontrivial definition of discrete topological charge. As already discussed in this introduction, the latter turns out to be one of the major challenges in this framework and it is analyzed in Section~\ref{sec:dis1}. Based on that, one defines a family of admissible discrete topological charges associated with a spin configuration $u$ (see Definition \ref{df:DisTopCha}). With that at hand, in Theorem~\ref{th:discomp} one proves that if a sequence \((u_\varepsilon)\subset\calSF_\varepsilon(\Omega)\) satisfies
$\sup_{\varepsilon}E_\varepsilon(u_\varepsilon,\Omega)<+\infty$, then, up to subsequences, any associated discrete topological charge \(\mu_\varepsilon\in\calM_{\mathrm{adm}}(\Omega;u_\varepsilon)\) converges (in the weak\(^*\) sense) to an atomic measure \(\mu\in\calM(\Omega)\) of the type $\mu=4\pi\sum_{i=1}^{N}d_i\delta_{x_i},\ d_i\in\mathbb{Z}$. Also in this case, the compactness result suggests to define a concentration energy obtained optimizing a certain charge over all possible spin fields associated to the same admissible measure. Hence we set $\mathscr{E}_\eps:\calM(\Omega)\to [0,+\infty]$ to be the family of functionals defined by
\begin{equation*}
\mathscr{E}_\eps(\mu)=\inf\{E_\eps(u)\colon u\in\calSF_\eps(\Omega), \text{ such that } \mu\in\calM_{adm}(\Omega; u) \}\,,
\end{equation*}
and in Theorem~\ref{th:main} we show that the sequence $\mathscr{E}_\eps$ $\Gamma\mbox{-}$converges with respect to the weak$^*$ topology of $\calM(\Omega)$ to the functional $\mathscr{E}:\calM(\Omega)\to[0,+\infty]$ which turns out to be finite on measures $\mu$ of the type $\mu=4\pi\sum_{i=1}^{N}d_i\delta_{x_i},\ d_i\in\mathbb{Z}$ where it takes the form 
\[
\mathscr{E}(\mu)=\sum_{i=1}^{N}\psi(d_i)\,. 
\]
Here $\psi(d)$ is given by the cell formula in \eqref{df:cell} and can be interpreted as the minimal discrete baby Skyrme energy associated with the formation of a singularity of degree $d$. Also in this setting our result is actually more general (see Section \ref{sec:Skyrme}) as it extends to energies including antisymmetric exchange interactions, namely a sum over all the elementary triangles $T=[i,j,k]$ of the lattice of terms of the type 
\[ 
A_\eps(u,T)=\frac12(D_\eps^{ij}\cdot(u(i)\times u(j))+D_\eps^{jk}\cdot(u(j)\times u(k))+D_\eps^{ki}\cdot(u(k)\times u(i)))\,,
\]
provided the (material dependent) interaction coefficients are of the type $D_\eps^{ij}=\eps D(i-j)$. We observe that an appropriate choice of $D$ makes the antisymmetric exchange interaction energy contribution the lattice analogue of the DMI interaction term considered in the continuum setting.\\

As a final remark we find it worth stressing that in both the continuum and the discrete setting our $\Gamma$-convergence result does not give information on the optimal value of the degree of minimizers. This appears to be a delicate question which requires a better understanding of the interactions between the degree one minimizers, namely the magnetic skyrmions.  \\

The rest of the paper is organized as follows. In Section~\ref{sec:bas} we introduce the basic notation and recall essential background material on Sobolev maps, Radon measures, and the topological degree. In Section~\ref{section:main-res} we state our main results precisely, both in the continuum and in the discrete setting. In particular, Theorem~\ref{th:contcomp} provides the compactness and concentration properties of the continuum energy \(F_\varepsilon\), while Theorem~\ref{th:contmain} establishes its \(\Gamma\)-convergence. These two theorems are then proved in Section~\ref{sec:cont} using blow-up and covering arguments. In Section~\ref{sec:dis1} we develop the discrete framework by introducing admissible interpolation regions and surfaces, which are crucial for defining a discrete notion of topological charge (see Definition~\ref{df:DisTopCha}). Section~\ref{sec:dis2} is then devoted to the proof of the compactness and \(\Gamma\)-convergence results in the discrete setting (Theorems~\ref{th:discomp} and \ref{th:main}). Finally, Appendix~\ref{sec:app1} and Appendix~\ref{sec:app2}  collect several technical lemmas that underpin the analysis throughout the paper.

\section{Basic Notation}\label{sec:bas}

\subsubsection*{Basic Notation}
We denote with $\mathS^2=\{y\in\mathR^3\colon\ |y|=1\}$ the set of unit vectors of $\mathR^3$ and we set $n=(0,0,1)$. Given $p,q\in\mathS^2$, we say that $p, q$ are antipodal if $p=-q$. If $p$ and $q$ are not antipodal we denote by $\gamma_{p,q}\colon[0,1]\to\mathS^2$ the unique lenght minimizing geodesic curve satisfying $\gamma(0)=p$ and $\gamma(1)=q$. With a slight abuse of notation,  we denote by $\gamma_{p,q}$ also the image of the curve in $\mathbb{S}^2$.

We denote by $\{e_1, e_2\}$ the standard basis of $\mathR^2$ and for every $x\in\mathR^2$ we set $x^{\perp}=(-x\cdot e_2, x\cdot e_1)$.
For $0<r<s$, we denote
\[
\begin{split}
& B_r=\{z\in\mathR^2\colon\ |x-z|<r\}\,,\quad A_{R,r}=\{z\in\mathR^2\colon\ r<|x-z|<R\}\,,\\
& Q_r=\left\{z\in\mathR^2\colon\ |z\cdot e_i|\le \frac r 2, \text{ for }i=1,2\right\}\,,
\end{split}
\]
and for $x\in\mathR^2$ we set $B_r(x)=B_r+x$, $A_{r,R}(x)=A_{r,R}+x$ and $Q_r(x)=Q_r+x$.

Given $\Omega$ an  open set of $\mathR^2$, we denote by $\calM(\Omega)$ the set of real valued Radon measures on $\Omega$ endowed with the usual notion of weak$^*$ convergence, and we call $\mu\in\calM(\Omega)$ an atomic  measure if 
\[\mu=\sum_{x\in E} c_x \delta_x\,,\]
for some finite set $E\subset\mathR^2$ and some family of coefficents $c_x\in\mathR$. We also say that $\mu$ has integer coefficents if $c_x\in\mathZ$ for every $x\in E$.
We observe that the set
\[\left\{\sum_{x\in E} c_x \delta_x\colon\, E\hbox{ finite},\, c_x\in\{\pm 1\}\right\}
\] 
is weakly$^*$ dense in the set of Dirac measures with integer coefficents.
Given $\mu\in\calM(\Omega)$ and $U\subseteq\Omega$ an open set, we define its flat norm $\|\mu\|_{\text{flat}, U}$ as follows: 
\[
\|\mu\|_{\text{flat}, U}=\sup \left\{ \int_U \phi\, {\rm d} \mu\colon\, \phi\in C^{\infty}_c(U),\, \|\phi\|_{L^\infty(U)}\le 1,\, \|\nabla \phi\|_{L^\infty(U)}\le 1 \right\}\,. 
\]
Whenever $\|\mu \|_{\text{flat}, U}$ is finite, we can equivalently compute the right hand side of the previous identity by taking the supremum among Lipschitz functions compactley suported in $U$. Note that if $(\mu_n)\subset\calM(\Omega)$ is such that $\sup_n|\mu_n|<\infty$ then $\|\mu_n\|_{\text{flat}, \Omega}\to 0$ if and only if $\mu_n\wto 0$.

\subsubsection*{Topological charge}
Throughout the paper, as it is customary, we set
\[
H^1(\mathR^2;\mathS^2)=\{u\in H^1(\mathR^2;\mathR^3),\,|u|=1 \text{ a.e. in }\mathR^2\}\,.
\]
Given $\Omega\subseteq\mathR^2$ a bounded open set  we then define $\calU(\Omega)\subset H^1(\mathR^2;\mathS^2)$ as follows
\begin{equation}\label{df:calU}
\calU(\Omega)=\{u\in H^1(\mathR^2;\mathS^2),\, u(x)=n \text{ for a.e. } x\in\mathR^2\setminus\Omega\}\,.
\end{equation}
For $u\in\mathcal{U}(\Omega)$  we define the  the topological charge density as
\[
q(u)=u\cdot \partial_x u\times \partial_y u\,.
\]
The topological charge or Skyrmion number of $u$ is then  given by 
\[
\mathrm{deg}(u)=\frac 1 {4\pi}\int_{\mathR^2}q(u)\,\mathrm{d}x=\frac 1 {4\pi}\int_{\Omega}q(u)\, \mathrm{d}x\,.
\]
For $u\in \calU(\Omega)\cap C^1(\bar\Omega;\mathS^2)$, the quantity $\mathrm{deg}(u)$ corresponds to the Brouwer degree of the self map of $\mathS^2$, $u\circ \sigma:\mathS^2\to\mathS^2$, where $\sigma$ is the stereographic projection mapping $n$ to $\infty$. In particular, $\mathrm{deg}(u)\in\mathZ$. This remains true for a general $u\in\calU(\Omega)$ by approximation, see \cite{ScUh}  (see also \cite[Lemma A.1]{BrCo}).
Moreover notice that for $u\in\calU(\Omega)$ it holds, by Schwarz inequality, $q(u)\in L^1(\Omega)$ so that $q(u)\mathrm{d}x\in\calM(\Omega)$.
\begin{remark}
For a Lipschitz map $u:\mathR^2\to\mathR^3$, with $|u|\equiv 1$, a simple algebraic computation shows that the Jacobian determinant, 
\[
J(u)=\sqrt{\det(\nabla u^*\nabla u)}
\] 
can be expressed as
\[
J(u)=|\partial_xu\times \partial_y u|=|q(u)|\,.
\]
In particular, by the area formula it holds
\begin{equation}\label{eq:area}
\int_{u(R)}\#\{x\in R:\, u(x)=y\}\,{\rm d} \mathcal{H}^2(y) =\int_{R}|q(u)|\,\mathrm{d}x\,,\quad \hbox{ for every measurable set }R\subset\mathR^2\,.
\end{equation}
\end{remark}

\subsubsection*{Simple regions}
In the sequel we will use the follwing notion of \textit{simple region},  see for instance \cite[Sec. 4.5]{Do}. Given an orientable surface $M$ (in this paper $M=\mathS^2$), a relatively open and connected subset $S\subset M$  is said to be a simple region on $M$ if it is homeomorphic to the disk and its boundary $\partial S$ is the trace of a simple closed piecewise regular  curve $\alpha: I\to \partial S$, being $I\subset\mathR$ some closed interval. 
In this case, it is customary to say that  $\alpha$ is positively oriented with respect to $S$, if the vector  $h(t)=n_M(\alpha(t))\times \alpha'(t)$, where $n_M(\alpha(t))$ is the normal vector field of $M$ at $\alpha(t)$,
points toward $S$, that is: for every $\eps>0$ and every differentiable curve  $\beta:(0,\eps)\to S$ such that $\beta(0)=\alpha(t)$ and $\beta'(0)\neq\alpha'(t)$ it holds
\[
h(t)\cdot \beta'(0)\ge 0\,.
\] 
The curve $\alpha$ is said to be negatively oriented with respect to $S$ if the previous inequality holds with the opposite sign. 
When $M=\mathR^2$, the curve $\alpha$ is positively oriented with respect to $S$, if it parametrizes $\partial S$ counterclockwise. 

We notice that the following are simple regions of $\mathS^2$:
\begin{itemize}
\item[-] Given three linear independent unitary vectors $p_1,p_2,p_3\in\mathS^2$, the geodesic triangle with vertices $p_1,p_2,p_3$, given by  
\[
S=S(p_1,p_2,p_3)=\left\{\frac{\sum_{i=1}^3 \lambda_i p_i}{\left|\sum_{i=1}^{3} \lambda_i p_i\right|}: \lambda_{i}> 0\right\}\subset\mathS^2\,.
\] 
In this case $\partial S=\bigcup_{i\neq j}\gamma_{p_i,p_j}$.
\item[-] Given four unitary vectors $p_1,p_2,p_3,p_4\in\mathS^2$, such that respectively $p_1,p_3,p_4$ and $p_2,p_3,p_4$ are linearly independent, the union of the two adjacent geodesic triangle 
\[
S=S(p_1,p_3,p_4)\cup S(p_2,p_3,p_4)\,.
\]
Here $\partial S=\gamma_{p_1,p_3}\cup \gamma_{p_3,p_2}\cup \gamma_{p_2,p_4}\cup \gamma_{p_1,p_2}$.
\item[-]For a unitary vector $h\in\mathS^2$, any hemisphere centered in $h$, that is 
\[
S=\{q\in\mathS^2\colon\, q\cdot h> 0\}\,.
\]
In this latter case $\partial S=\{q\in\mathS^2:\, q\cdot h= 0\}$, is commonly called a great-circle of $\mathS^2$.
\end{itemize}

\subsubsection*{Triangular lattice}
We consider the triangular lattice $\calL\subset\mathR^2$, given by
\begin{equation}\label{df:lattice}
\calL=
\left\{a
 e_1
+b\hat e_2
\colon a,b\in\mathZ\right\}\,,
\quad \hbox{where }\, \hat e_2=\frac{1}{2} (e_1+\sqrt{3}e_2)\,.
\end{equation}
For $i,j,k\in\mathR^2$ we denote by $[i,j,k]$ and $[i,j]$  the closed convex hull of $i,j,k$ and of $i,j$ respectively. We set $\calT$ to be  the family of triangles subordinated to the lattice $\calL$
\[
\calT=\{[i, j, k]\colon i,j,k\in \calL,\ |i-j|=|i-k|=|j-k|=1\}\,,
\]
and  $\calE$ to be the family of edges subordinated to the lattice $\calL$
\[
\calE=\{[i,j]\colon\ i,j\in\calL,\, |i-j|=1\}\,.
\]
For $\eps>0$ we consider the rescaled versions of $\calL$, $\calT$ and $\calE$, namely $\calL_\eps=\eps\calL$, $\calT_\eps=\eps\calT$ and $\calE_\eps=\eps\calE$. 
For $A\subseteq \mathR^2$ we use the following notation:
\[
\begin{split}
&\calL_\eps(A)=\{i\in\calL_\eps(\mathR^2)\colon i\in A\}\,,\\
&\calT_\eps(A)=\{T\in\calT_\eps(\mathR^2)\colon T\subseteq A\}\,,\\
&\calE_\eps(A)=\{\Gamma\in\calE_\eps(\mathR^2)\colon \Gamma\subseteq A\}\,.
\end{split} 
\]

In what follows, we say that a closed subset $R\subset\mathR^2$ is a $\eps$-triangular region, or simply a triangular region, if it is obtained as finite union of triangles in $\calT_\eps$. 

We denote by $\mathcal{SF}_\eps$ the set of \textit{spin fields}:
\[
\mathcal{SF}_\eps=\{u\colon\calL_\eps\to \mathS^2\}\,.
\]
Given  $\Omega\subset\mathR^2$ a bounded open set, we define, in analogy with \eqref{df:calU}, the following set of spin fields satisfying a discrete Dirichlet boundary condition
\[
\mathcal{SF}_\eps(\Omega)=\{u\colon\calL_\eps\to \mathS^2\colon\ u_\eps\equiv n \text{ in }\calL_\eps(\mathR^2\setminus\Omega)\}\,.
\]

\section{Statment of the main results}\label{section:main-res}

In this paper, for a given bounded  open set $\Omega\subset\mathR^2$  and for $\eps>0$,  we consider two different types of energies: one defined for maps in $H^1(\mathR^2;\mathS^2)$ and one defined for spin fields in $\calSF_\eps$.

More precisely, for $u\in H^1(\mathR^2;\mathS^2)$ we define 
\begin{equation}\label{df:F}
F_\eps(u,\Omega)=\frac 1 2 \int_\Omega |\nabla u|^2\,\mathrm{d}x+\frac{1}{\eps^2}\int_\Omega |u-n|^2\, \mathrm{d}x\,.
\end{equation}
Instead, for $u\in\calSF_\eps$, we define
\begin{equation}\label{df:E}
E_\eps(u,\Omega)=\sum_{T\in\calT_\eps(\Omega)}H_\eps(u, T)+\sum_{T\in\calT_\eps(\Omega)}Z_\eps(u, T)\,,
\end{equation}
where for $T=[i,j,k]\in\calT_\eps(\mathR^2)$, $H_\eps(u,T)$ and $Z_\eps(u,T)$ are respectively given by
\[
H_\eps(u,T)=\frac{1}{2}(|u(i)-u(j)|^2+|u(j)-u(k)|^2+|u(k)-u(i)|^2)
\]
and
\[
Z_\eps(u,T)=\frac{1}{2}(|u(i)-n|^2+|u(j)-n|^2+|u(k)-n|^2)\,.
\]
For notation convenience, when in the  right hand side of \eqref{df:F}, respectively \eqref{df:E}, $\Omega$ is replaced by some measurable set $A\subset\mathR^2$ we write $F_\eps(u, A)$, respectively $E_{\eps}(u,A)$.

Our first result is as follows.   
\begin{theorem}\label{th:contcomp}
Let $\Omega\subset\mathR^2$ be a bounded open set and let $(u_\eps)\subset H^1(\mathR^2;\mathS^2)$ satisfy
\[
\sup_{\varepsilon>0} F_\eps(u_\eps, \Omega)<+\infty\,.
\] Then, there exists an atomic measure $\mu\in\calM(\Omega)$  such that (up to subsequences) $q(u_\eps)\mathrm{d}x$ weak$^*$ converge to $\mu$. Moreover $\frac \mu {4\pi}$ has integer coefficents.
\end{theorem}
Building up from Theorem \ref{th:contcomp} we deduce the following $\Gamma$-convergence results  for the family of functionals $\mathscr{F}_\eps:\calM(\Omega)\to [0,+\infty]$:
\begin{equation}\label{df:scrFeps}
\mathscr{F}_\eps(\mu)=\inf\left\{F_\eps(u)\colon u\in\calU(\Omega),\quad q( u_\eps)\, \mathrm{d} x=\mu\right\}\,.
\footnote{Hereinafter the convention $\inf\emptyset=\infty$ is adopted .} 
\end{equation}
\begin{theorem}\label{th:contmain}
Let $\Omega\subset\mathR^2$ be a bounded open set. Let $\mathscr{F}:\calM(\Omega)\to [0,\infty]$ be the functional defined as
\begin{equation*}
\mathscr{F}(\mu)=\begin{cases}
|\mu|(\Omega) &\text{if }\frac{\mu}{4\pi} \text{ is an atomic measure with integer coefficents,}\\
+\infty & \text{ otherwise}.
\end{cases}
\end{equation*}
Then, the sequence of functionals $\mathscr{F}_\eps$ $\Gamma\mbox{-}$converge with respect to the weak$^*$ topology of $\calM(\Omega)$ to $\mathscr{F}$ as $\eps\to 0$.
\end{theorem}

Theorem  \ref{th:contmain} provides the protototypical result that we aim to achieve in the discrete setting. 
However, in this setting a major difficulty appears already in order to give a rigorous notion of topological charge density for a spin field $u\in\calSF_\eps$.  We do this by partioning the plane into triangular regions whose boundaries do not contain pairs of nearest neighborhoods on which the spin field $u$ assumes antipodal  values. On each of such region it is then possible to define a proper continuos inteprolation and accordingly a notion of topological charge. Nevertheless, this procedure cannot be free from a certain level of arbitrariness as explained below. As a consequence we end up defining, for each $u\in\calSF_\eps$, a set of admissible topological charges $\calM_{\text{adm}}(\Omega; u_\eps)\subset\calM(\Omega)$  (see Definition~\ref{df:DisTopCha}).
We then prove the following compactness result.

\begin{theorem}\label{th:discomp}
Let $\Omega\subset\mathR^2$ a bounded open set and let $(u_\eps)\subset\calSF_\eps(\Omega)$ satisfy
\[
\sup_{\varepsilon >0} E_\eps(u_\eps,\Omega)<+\infty\,.
\]
Let $(\mu_\eps)\subset\calM(\Omega)$ be a family of atomic measures such that for every $\eps>0$ it holds $\mu_\eps\in\calM_{adm}(\Omega; u_\eps)$. 
Then, there exists an atomic measure $\mu\in\calM(\Omega)$  such that (up to subsequences) $\mu_\eps$ weak$^*$ converge to $\mu$. Moreover $\frac \mu {4\pi}$ has integer coefficents.  
\end{theorem}
At last, letting $\mathscr{E}_\eps:\calM(\Omega)\to [0,+\infty]$ to be the family of functionals defined by
\begin{equation}\label{df.scrEeps}
\mathscr{E}_\eps(\mu)=\inf\{E_\eps(u)\colon\ u\in\calSF_\eps(\Omega),\, \text{ such that } \mu\in\calM_{adm}(\Omega; u) \}\,,
\end{equation}
we show that, for every $d\in\mathZ$,  the following quantity is well defined
\begin{equation}\label{df:cell}
\psi(d)=\lim_{\eps\to 0}\left(\inf\left\{E_{\eps}(w, B_1)\colon\, 
\begin{aligned}
& w\in\calSF_\eps(B_{1-3\eps}), \text{ s.t. there exists }\zeta\in\calM_{adm}(\Omega; w)\\
& \text{ with } \zeta(B_1)=4\pi d
\end{aligned}
\right\}\right)\,,
\end{equation}
and we prove the following convergence result.

\begin{theorem}\label{th:main}
Let $\Omega\subset\mathR^2$ be a bounded open set. Let $\mathscr{E}\colon\calM(\Omega)\to [0,+\infty]$ be the functional defined by
\[
\mathscr{E}(\mu)=
\begin{cases}
\sum_{i=1}^N \psi(d_i) & \hbox{if }\frac{\mu}{4\pi }=\sum_i d_i \delta_{x_i}, d_i\in\mathZ\,, \\
+\infty &\text{otherwise.}
\end{cases}
\]
Then, the sequence of functionals $\mathscr{E}_\eps$ $\Gamma\mbox{-}$converge with respect to the weak$^*$ topology of $\calM(\Omega)$ to $\mathscr{E}$ as $\eps\to 0$. 
\end{theorem}

\subsection{The baby Skyrme model with diverging Zeeman term} \label{sec:Skyrme}
In this paragraph we briefly show how to extend the results stated above to cover the full baby Skyrme model both in the continuum and in the discrete setting.

For $u\in H^1(\mathR^2;\mathS^2)$ the DMI energy $D(u,\Omega)$ and the full baby Skyrme energy $F^{\rm BS}_\eps(u,\Omega)$ are given by
\begin{equation*}
D(u,\Omega)=\int_\Omega u\cdot(\nabla \times u)\,\mathrm{d}x\,.
\end{equation*}
and 
\begin{equation*}
 F^{\rm BS}_\eps(u,\Omega)=F_\eps(u,\Omega)+D(u,\Omega)\,,
\end{equation*}
respectively. We observe that by Young's inequality 
\[
|D(u,\Omega)|\leq  \frac 1 2 F_\eps(u_\eps,\Omega)+4|\Omega|\,,
\] 
hence 
\[
F^{\rm BS}_\eps(u_\eps,\Omega)\geq  \frac 1 2 F_\eps(u_\eps,\Omega)-4|\Omega|\,.
\]
As a result Theorem~\ref{th:contcomp} applies with $F^{BS}_\eps$ in place of $F_\eps$. We further observe that for sequences $(u_\eps)$ with $\sup_{\varepsilon >0} F_\eps(u_\eps,\Omega)<+\infty$, by the weak $L^2$ convergence of $\nabla u_\eps$ and the strong $L^2$ convergence of $u_\eps$ to $n$, $D(u_\eps,\Omega)$ vanishes. As a result Theorem \ref{th:contmain} holds with $F^{\rm BS}_\eps$ in place of $F_\eps$.\\

\noindent As for the discrete setting, for $u\in\calSF_\eps$, we define the DMI interaction, or antisymmetric exchange term,  as
\begin{equation*}
A_\eps(u,\Omega)=\sum_{T\in\calT_\eps(\Omega)}A_\eps(u, T)\,,
\end{equation*}
where for $T=[i,j,k]\in\calT_\eps$ we have set 
\[ 
A_\eps(u,T)=\frac12(D_\eps^{ij}\cdot(u(i)\times u(j))+D_\eps^{jk}\cdot(u(j)\times u(k))+D_\eps^{ki}\cdot(u(k)\times u(i)))\,,
\]
being $D_\eps^{ij}$ a family of vectors in $\mathR^3$. Furthermore, we introduce the discrete baby Skyrme functional 
\[ 
E^{\rm BS}_\eps(u,\Omega)=E_\eps(u,\Omega)+A_\eps(u,\Omega)\,.
\]
Similarly to what was noted above we have
\[
\left|A_\eps(u,T)\right|\le \frac 1 2 H_\eps(u, T)+ \frac{2\sqrt{3}} {\eps^2} \|D_\eps\|^2_\infty|T|\,,
\]
which implies
\[ 
\left|A_\eps(u,\Omega)\right|\le  \frac 1 2 H_\eps(u, \Omega)+\frac{2\sqrt{3}} {\eps^2}\, \|D_\eps\|^2_\infty|\Omega|\,.
\]
Assuming 
\[
 D_\eps^{ij}=\eps D(i-j),\quad \text{for some }D:\{\pm e_1, \pm \hat e_2, \pm (e_1 - \hat e_2)\}\to \mathbb{R}^3\,,
\] 
one can argue as above and conclude that both Theorem \ref{th:discomp} and Theorem\ref{th:main} apply with $E^{\rm BS}_\eps$ in place of $E_\eps$.

\section{Concentration of topological charge}\label{sec:cont}

We prove Theorem \ref{th:contcomp} by combining a blow-up argument, proved in Lemma~\ref{lm:blowup} and a suitable covering argument. Throughout the section we assume $\Omega$ to be a bounded  open set and $(u_\eps)\subset  H^1(\Omega,\mathS^2)$ to be a sequence satisfying 
\begin{equation}\label{eq:eqbndF}
\sup_{\varepsilon>0} F_\eps(u_\eps, \Omega)<+\infty\,.
\end{equation}
Within the proofs of the results of this section we use the symbol $C$ to denote some universal constant $C>0$ which may depend on $\Omega$ and on $\sup_{\eps>0}F_{\eps}(u_\eps, \Omega)$ and whose value may change from line to line. 

\begin{lemma}\label{lm:blowup}
Let  $(u_\eps)\subset H^1(\mathR^2;\mathS^2)$ satisfying \eqref{eq:eqbndF}, be such that $q(u_\eps)\,\mathrm{d}x\wto \mu$ for some $\mu\in\calM(\Omega)$. Then, $\mu(\{x\})\in 4\pi \mathZ$ 
for every $x\in\Omega$.
\end{lemma}

\begin{proof}
Let $x\in\Omega$ and suppose without loss of generality that $x=0$. 
We set $\nu_\eps=|\nabla u_\eps|^2\,\mathrm{d}x\in\calM(\Omega)$, and observe that \eqref{eq:eqbndF} implies $\sup|\nu_\eps|(\Omega)<+\infty$. Up to a not relabelled subsequence, we have that $\nu_\eps\wto\nu$, for some $\nu\in\calM(\Omega)$. 
We select a sequence of balls $(B_{r_m})$ with vanishing radii  such that $B_{r_m}\subset\Omega$ for every $m\in\mathN$. We also require that
\begin{equation}\label{eq:perilim1}
\mu(\partial B_{r_m})=0\,,
\end{equation}
and
\begin{equation}\label{eq:perilim2}
\nu(\partial B_{r_m})=0\,,
\end{equation}
for every $m\in\mathN$. Thanks to \eqref{eq:perilim1}, we have
\[
\begin{split}
\mu(\{0\})=
\lim_{m\to +\infty}\lim_{\eps\to 0}\int_{B_{r_m}}q(u_\eps)\,\mathrm{d}x=\lim_{m\to +\infty}\lim_{\eps\to 0}\int_{B_1}q(u^m_\eps)\,\mathrm{d}x\,.
\end{split}
\]
where $u_\eps^m:B_1\to \mathS^2$ is defined by $u_\eps^m(x)=u_\eps(r_mx)$. Thus in order to conclude it is enough to show that right hand side of the previous identity belongs to $4\pi\mathZ$. 

For the sake of brevity we set $\nu_\eps^m=|\nabla u^m_\eps|^2dx$, and denote $\nu^m$  the measure on $B_1$, defined on every measurable set $A\subset B_1$ as
\[
\nu^m(A)=\nu(r_m A)\,.
\]
It is straightforward to check that $\nu^m_\eps\wto\nu^m$ when $\eps \to 0$. Furthermore \eqref{eq:perilim2} implies
\begin{equation}\label{eq:perillim3}
\nu^m(\partial B_1)=0\,,
\end{equation}
for every $m\in\mathN$.

Let us introduce a second familiy of parameters $(\rho_m)\subset(0,1/4)$, chosen in such a way that $\rho_m\searrow 0$  and
\begin{equation}\label{eq:perillim4}
\nu^m(\partial B_{1-2\rho_m})=0\,,\quad \text{and}\quad \nu^m(\bar A_{1,1-2\rho_m})<2^{-m}\,,
\end{equation}
for every $m\in\mathN$. 
We employ for the sake of brevity the following notation
\[
f_{m,\eps}:(1-2\rho_m,1-\rho_m)\to \mathR\,,\quad f_{m,\eps}(s):=\int_{\partial B_s}\left(|\nabla u^m_\eps|^2 +\frac{r_m^2}{\eps^2} |u^m_\eps-n|^2\right)\,\mathrm{d} \mathcal{H}^1\,.
\]
Notice that the coarea formula together with a simple change of variables and \eqref{eq:eqbndF} give
\begin{equation}\label{eq:coareasemplice}
\begin{split}
C\geq F_{\eps}(u_\eps, r_m A_{1-\rho_m,1-2\rho_m})\ge\int_{1-2\rho_m}^{1-\rho_m}f_{m,\eps}(s)\, \mathrm{d}s\,.
\end{split}
\end{equation}
For every $t>0$, we let
\[
I_t=\left\{s\in (1-2\rho_m, 1-\rho_m)\colon f_{m,\eps}(s)\le t \right\}\,.
\]
Markov-Chebyshev inequality gives:
\[
\mathcal{H}^1(I_t)\ge \rho_m-\frac{1}{t}\int_{1-2\rho_m}^{1-\rho_m}f_{m,\eps}(s)\,\mathrm{d}s\,.
\]
By letting $t_{m,\eps}=\frac{2}{\rho_m}\int_{1-2\rho_m}^{1-\rho_m}f_{m,\eps}(s)\, \mathrm{d}s$, the previous inequality and \eqref{eq:coareasemplice} respectively give
\begin{equation}\label{eq:Markov}
\mathcal{H}^1(I_{t_{m,\eps}})\ge \frac{\rho_m}{2}\,,\quad t_{m,\eps}\le \frac{2C}{\rho_m}\,.
\end{equation}
By applying again the coarea formula we select $\overline\rho_{m,\eps}\in (1-\rho_m,1-2\rho_m)\cap I_{t_{m,\eps}}$ such that
\[
\begin{split}
&\int_{\partial B_{\overline\rho_{m,\eps}}}|\nabla u^m_\eps|^2\,\mathrm{d}\mathcal{H}^1\le  \frac{1}{\mathcal{H}^1(I_{t_{m,\eps}})}\int_{I_{t_{m,\eps}}}\int_{\partial B_s}|\nabla u^m_\eps|^2\,\mathrm{d}\mathcal{H}^1\,\mathrm{d}s\,.
\end{split}
\]
In particular, thanks to \eqref{eq:eqbndF}, \eqref{eq:perillim3}, \eqref{eq:perillim4} and \eqref{eq:Markov},  for $\varepsilon>0$ small enough, we get 
\begin{equation}\label{eq:civuole}
\int_{\partial B_{\overline\rho_{m,\eps}}}|\nabla u^m_\eps|^2\,\mathrm{d}\mathcal{H}^1\le \frac{C}{\rho_m}2^{-m}\,.
\end{equation}
Moreover, since  $\overline\rho_{m,\eps}\in I_{t_{m,\eps}}$, by \eqref{eq:Markov} we also have
\[
\begin{split}
&\int_{\partial B_{\overline\rho_{m,\eps}}}|\nabla u^m_\eps|^2\,\mathrm{d}\mathcal{H}^1+\frac{r_m^2}{\eps^2}\int_{\partial B_{\bar\rho_{m,\eps}}}|u^m_\eps-\boldn|^2\,\mathrm{d}\mathcal{H}^1=f_{m,\eps}(\overline\rho_{m,\eps})\le t_{\eps,m}\le \frac{2C}{\rho_m}\,.
\end{split}
\]
This,  implies
\begin{equation}\label{eq:civuole2}
\frac{r_m^2}{\eps^2}\int_{\partial B_{\bar\rho_{m,\eps}}}|u^m_\eps-\boldn|^2\,\mathrm{d}\mathcal{H}^1\le \frac{C}{\rho_m}\,.
\end{equation}
Let $\tilde{u}^m_\eps\colon\partial B_1\to\mathS^2$ be the rescaled maps defined by $\tilde u^m_\eps(x)=u^m_\eps(\bar\rho_{m,\eps}x)$. From \eqref{eq:civuole2} we deduce that, as $\eps\to 0$, the sequence
$(\tilde{u}^m_\eps-\boldn)$ converges to $0$ strongly in $L^2(\partial B_1)$. Moreover, thanks to \eqref{eq:civuole}, we  have that the sequence $(\tilde{u}^m_\eps-\boldn)$ is  bounded in $H^1(\partial B_1)$, uniformly in $\eps$.   In view of the embedding $H^1(\partial B_1)\hookrightarrow C^{0,1/2}(\partial B_1)$, we thus have that $(\tilde{u}^m_\eps-\boldn)$ converges up to some not relabelled subsequence uniformly in $\partial B_1$, as $\eps\to 0$. This in turn implies that  $\|\tilde u_\eps^m-\boldn\|_{L^{\infty}(\partial B_1)}\to 0$ as $\eps\to 0$ and therefore that
\begin{equation}\label{eq:civuole3}
\|u_\eps^m-\boldn\|_{L^{\infty}(\partial B_{\overline{\rho}_{m,\eps}})}\to 0
\end{equation}
as $\eps\to 0$.
Let $\varphi_{\rho_m} \in C_c^\infty(B_1)$ be such that $\varphi_{\rho_m} \equiv 1$ on $B_{1-\rho_m}$ and $\|\nabla \varphi_{\rho_m}\|_{L^\infty(B_1)} \leq C \rho_m^{-1}$. We then define $v_\varepsilon \colon A_{1, \overline{\rho}_{m,\varepsilon}} \to \mathbb{R}^3$ to be
\begin{align*}
v_\varepsilon(x) = \varphi_{\rho_m}(x) u_\varepsilon\left(\frac{{\rho}_{m,\varepsilon}x}{|x|}\right) + (1-\varphi_{\rho_m}(x)) \boldn\,.
\end{align*}
From \eqref{eq:civuole3} we have $|v^m_\eps|\ge \frac{1}{2}$ for $\eps$ small enough. In particular, for $\eps$ small enough, the map $V^m_\eps\in H^1(B_1; \mathR^3)$ given by
\[
V^m_\eps(x)=
\begin{cases}
\frac{v^m_\eps(x)}{|v^m_\eps(x)|}&\text{if } x\in A_{\bar\rho_{m,\eps},1}\,,\\
\boldn &\text{if } x\in \partial B_1\,,\\
u^m_\eps(x) &\text{if } x \in B_{\bar\rho_{m,\eps}}\,.
\end{cases}
\]
is well defined and we have $V^m_\eps\in\calU(B_1)$.  As noticed in Section \ref{sec:bas}, this provides
\[
\int_{B_1}q(V^m_\eps)\,dx\in 4\pi\mathZ \,.
\]
Note that
\[
\int_{B_1} q(u^m_\eps)\,\mathrm{d}x=\int_{B_1} q( V^m_\eps)\,\mathrm{d}x-\int_{A_{1,\bar\rho_{m,\eps}}} q(V^m_\eps)\,\mathrm{d}x+\int_{A_{1,\bar\rho_{m,\eps}}} q(u^m_\eps)\,\mathrm{d}x\,,
\]
Hence to conclude it is enough to prove that the last two terms of the previous identity vanish in the limit as first $\eps\to 0$ and then $m\to\infty$. With this aim, we first observe that, Schwarz inequality, \eqref{eq:perillim3} and \eqref{eq:perillim4} give
\[
\lim_{m\to\infty}\lim_{\eps\to 0}\left|\int_{A_{1,\bar\rho_{m,\eps}}} q(u^m_\eps)\,\mathrm{d}x\right|\le \lim_{m\to\infty}\lim_{\eps\to 0} \int_{A_{1,1-2\rho_m}}|\nabla u^m_\eps|^2\,\mathrm{d}x=\lim_{m\to\infty}\nu^m(A_{1,1-2\rho_m})=0\,.
\]
Additionally, note that for every $x\in A_{1,\overline{\rho}_{m,\eps}}$ we have
\begin{align*}
|\nabla v_\varepsilon(x)|^2 &\leq 2|\nabla \varphi_{\rho_m}(x)|^2 \left|u_\varepsilon\left(\overline{\rho}_{m,\varepsilon}\frac{x}{|x|}\right)-\boldn\right|^2 + \frac{4\overline\rho_{m,\eps}}{|x|}\left|\nabla u_\varepsilon\left(\overline{\rho}_{m,\varepsilon}\frac{x}{|x|}\right)\right|^2 \\&\leq \frac{C}{\rho_m^2}\left|u_\varepsilon\left(\overline{\rho}_{m,\varepsilon}\frac{x}{|x|}\right)-\boldn\right|^2 + C\left|\nabla u_\varepsilon\left(\overline{\rho}_{m,\varepsilon}\frac{x}{|x|}\right)\right|^2\,.
\end{align*} 
Integrating the previous inequality on $A_{1,\bar\rho_{m,\eps}}$ and using the coarea formula and \eqref{eq:civuole} we get
\begin{align*}
\int_{A_{1,\bar\rho_{m,\eps}}} |\nabla v^m_\eps|^2\,\mathrm{d}x&\leq \frac{C}{\rho_m} \|u_\eps^m-\boldn\|_{L^{\infty}(\partial B_{\overline{\rho}_{m,\eps}})}^2 + C\rho_m\int_{\partial B_{\overline{\rho}_{m,\eps}}}|\nabla u_\varepsilon|^2\,\mathrm{d}\mathcal{H}^1 \\& \leq \frac{C}{\rho_m} \|u_\eps^m-\boldn\|_{L^{\infty}(\partial B_{\overline{\rho}_{m,\eps}})}^2 +C2^{-m}\,.
\end{align*}
 This implies
\[
\lim_{m\to\infty}\lim_{\eps\to 0}\int_{A_{1,\bar\rho_{m,\eps}}} |\nabla v^m_\eps|^2\,\mathrm{d}x=0\,.
\] 
This in turn provides 
\[
\begin{split}\lim_{m\to\infty}\lim_{\eps\to 0}\left|\int_{A_{1,\bar\rho_{m,\eps}}} q(V^m_\eps)\,\mathrm{d}x\right|&\le \lim_{m\to\infty}\lim_{\eps\to 0}\int_{A_{1,\bar\rho_{m,\eps}}} |\nabla V^m_\eps|^2\,\mathrm{d}x \\&\le 4\lim_{m\to\infty}\lim_{\eps\to 0}\int_{A_{1,\bar\rho_{m,\eps}}} |\nabla v^m_\eps|^2\,\mathrm{d}x=0
\end{split}
\]
and concludes the proof. 
\end{proof}

\begin{proof}[Proof of Theorem \ref{th:contcomp}]
The existence of $\mu$ follows by weak$^*$ compactness. Thanks to Lemma \ref{lm:blowup}, to show the desired result, it is enough to prove that there exist  $N\in\mathN$ and $\{x_1,\dots, x_{N}\}\subseteq\Omega$ such that 
\begin{equation}\label{eq:LASTIMA2}
|\mu|(\Omega\setminus \{x_1,\dots, x_{N}\})=0\,.
\end{equation}  
Furthermore, we can assume that $(u_\eps)\subset C^{\infty}(\mathR^2;\mathS^2)$. This follows at once, using a standard diagonal argument, being $C^{\infty}(\mathR^2;\mathS^2)$ dense into $H^1(\mathR^2;\mathS^2)$ with respect to the strong $H^1$ topology, see \cite{ScUh} (see also \cite[Lemma A.1 ]{BrCo}). For the sake of readability, we divide the proof into several steps. \\

\noindent\textit{STEP 1.} We show that there exists $\eps_0>0$, which depends on $\Omega$ and $\sup_{\eps>0} F_\eps(u_\eps, \Omega)$, such that for every $\eps<\eps_0$, we can construct a family of disjoint balls $\mathbf{B}^{\eps}$ such that
\begin{equation}\label{eq:ballconstr1}
\{x\in\Omega\colon|u_\eps(x)-\boldn|^2>\sqrt{\eps}\}\subseteq \bigcup_{B\in\mathbf{B}^{\eps}}B
\end{equation}
and
\begin{equation}\label{eq:ballconstr2}
\sum_{B\in\mathbf{B}^{\eps}}r_B\le  C\sqrt{\eps} F_\eps(u_\eps;\Omega)\,,
\end{equation}
where $C>0$ is a constant depending only on $\Omega$  and for $B\in\mathbf{B}^{\eps}$, $r_B$ denotes the radius of $B$.
We proceed as follows. We apply Young inequality to get
\[
\begin{split}
&2\eps  F_\eps(u_\eps;\Omega)\ge \eps \int_\Omega |\nabla u_\eps|^2\,\mathrm{d}x+\frac 1 \eps \int_\Omega |u_\eps-\boldn|^2\,\mathrm{d}x\ge \int_\Omega |(u_\eps-\boldn)\nabla u_\eps|\,\mathrm{d}x= \frac 1 2 \int_\Omega |\nabla |u_\eps-\boldn|^2|\,\mathrm{d}x\,.
\end{split}
\]
By coarea formula, the latter provides
\begin{equation*}
\int_{0}^{\infty}P(\{|u_\eps-\boldn|^2>t\};\Omega)\,\mathrm{d}t\le 4\eps  F_\eps(u_\eps;\Omega)\,.
\end{equation*}
Thus, there exists $\tau_\eps\in (\sqrt{\eps}/2,\sqrt{\eps})$ such that
\begin{equation}\label{eq:perbnd}
P(\{|u_\eps-\boldn|^2>\tau_\eps\};\Omega)\le 8\sqrt{\eps} F_\eps(u_\eps;\Omega)\,.
\end{equation}
Furthermore, by Sard's theorem we can assume $\tau_\eps$ to be such that the open set $\{|u_\eps-\boldn|^2>\tau_\eps\}$ has smooth boundary. 
By the very definition of $F_\eps$ we also have
\begin{equation}\label{eq:volbnd}
|\{|u_\eps-\boldn|^2>\tau_\eps\}|\le 2\eps\sqrt{\eps} F_\eps(u_\eps;\Omega)\,.
\end{equation}
Let $\eps_0$ be such that  
\[2\eps_0\sqrt{\eps_0} \sup_\eps F_{\eps}(u_\eps; \Omega)< |\Omega|/2\,.
\]
Let $M_\eps$ be the number of connected components of $\{|u_\eps-\boldn|^2>\tau_\eps\}$ and denote by $\{|u_\eps-\boldn|^2>\tau_\eps\}_n$ denote the $n$-th. 
Then, for every $\eps<\eps_0$, \eqref{eq:volbnd} and the relative isoperimetric inequality  give
\begin{equation}\label{eq:relaiso}
\sqrt{|\{|u_\eps-\boldn|^2>\tau_\eps\}_n|}\le C P(\{|u_\eps-\boldn|^2>\tau_\eps\}_n;\Omega)\,.
\end{equation}
Since, for every connected open set of finite perimeter $E\subset\Omega$, it holds
\begin{equation}\label{eq:stimaperdia}
\mathrm{diam}(E)\le C (\sqrt{|E|}+P(E;\Omega))\,,
\end{equation}
by combining \eqref{eq:relaiso} and \eqref{eq:stimaperdia}, we deduce that  there exists $C>0$ depending only on $\Omega$ such that 
\begin{equation}\label{eq:stimaperdia2}
\mathrm{diam}(\{|u_\eps-\boldn|^2>\tau_\eps\}_n)\le C P(\{|u_\eps-\boldn|^2>\tau_\eps\}_n;\Omega)\,.
\end{equation}
Notice that achiving the inequality \eqref{eq:stimaperdia} is standard when $\Omega=\mathR^2$, while we refer to  Lemma \ref{lm:relper} in Appendix \ref{sec:app1} to treat the slighty more delicate case when $\Omega$ is a regular  bounded open set.  

 We follow a standard merging procedure. 
We denote by $r_n$ the  left hand side of \eqref{eq:stimaperdia2} and for every $n=1,\dots, M_\eps$ we select $x_n\in \{|u_\eps-\boldn|^2>\tau\}_n$ and let 
\[
\mathbf{B}^\eps_1=\{B_{r_n}(x_n):\ n=1,\dots, M_\eps\}\,.
\]
Note that by construction we have $\{|u_\eps-\boldn|^2>\tau\}_n\subset B_{r_n}(x_n)$.
If $\mathbf{B}^\eps_1$ is made of disjoint balls we set $\mathbf{B}^\eps=\mathbf{B}^\eps_1$. If this is not the case, we take any pair of indices $n_1,n_2\in\{1,\dots, M_\eps\}$, $n_1\neq n_2$, for which $B_{r_{n_1}}(x_{n_1})\cap B_{r_{n_2}}(x_{n_2})\neq\emptyset$ and we consider the ball $B_{r_{n_1+n_2}}(x_{n_1,n_2})$, where 
\[
x_{n_1,n_2}=\frac{r_{n_1}}{r_{n_1}+r_{n_2}}x_{n_1}+\frac{r_{n_2}}{r_{n_1}+r_{n_2}}x_{n_2}\,.
\]
It is easy to verify that 
\[B_{r_{n_1}}(x_{n_1})\cup B_{r_{n_2}}(x_{n_2})\subset B_{r_{n_1}+r_{n_2}}(x_{n_1,n_2})\,.
\]
We then define 
\[
\mathbf{B}^\eps_2=\{B_{r_{n_k}}(x_{n_k})\colon\, k\neq 1,2\}\cup\{B_{r_{n_1}+r_{n_2}}(x_{n_1,n_2})\}\,.
\]
If $\mathbf{B}^\eps_2$ is made of disjoint balls we set $\mathbf{B}^\eps=\mathbf{B}^\eps_2$. Otherwise we repeat the merging procedure. We iterate this scheme $m$ times until we find a family of disjoint balls. We then set $\mathbf{B}^\eps=\mathbf{B}^\eps_m$. By construction we have
\[
\sum_{B\in\mathbf{B}^\eps} r_B= \sum_{n=1}^{M_\eps} \mathrm{diam}(\{|u_\eps-\boldn|^2>\tau\}_n;\Omega)\,.
\] 
By combining the latter inequality with \eqref{eq:perbnd} and \eqref{eq:stimaperdia2} we deduce
 \eqref{eq:ballconstr2}.
Moreover, by recalling that $\tau<\sqrt{\eps}$, we have
\[
\{|u_\eps-\boldn|^2>\sqrt{\eps}\}\subset\{|u_\eps-\boldn|^2>\tau\}\subset\bigcup_{B\in \mathbf{B^\eps}} B\,.
\]
So \eqref{eq:ballconstr1} holds as well. \\

\noindent\textit{STEP 2.}
For every $\eta>0$ we define a sub-family of $\mathbf{B}^\eps$ as follows:
\begin{equation}\label{eq:perbnd2}
\mathbf{B}^{ \eps}_{\eta}=\{B\in \mathbf{B}^{\eps}\colon F_\eps(u_\eps; B)> \eta \}\,.
\end{equation}
Notice that, by a simple energetic consideration, the number of balls in $\mathbf{B}^{\eps}_{\eta}$ is finite.
Taking into account \eqref{eq:ballconstr2} we conclude that the Hausdorff limit, as $\eps\to 0$, of the set $ \bigcup_{\mathbf{B}^\eps_\eta} \bar B$ consists of a finite family of points. We call them  $\{x_1,\dots, x_{N_\eta}\}$. 
\\

\noindent\textit{STEP 3.}
 We fix $\phi\in C^{\infty}_c(\Omega\setminus \{x_1,\dots, x_{N_\eta}\})$.
Let $\sigma>0$ and
\[
\Omega_{\sigma}= \Omega\setminus \bigcup_{i}^{N_\eta} \bar B_\sigma(x_i)\,.
\]
 be such that the $\mathrm{spt}\phi\Subset\Omega_\sigma$.
By Hausdorff convergence, there exists  $0<\eps_1<\eps_0$ such that
\[
B\subset \bigcup_{i=1}^{N_\eta} B_\sigma(x_i)\,,\hbox{ for every }\eps<\eps_1\,, \quad B\in\mathbf{B}^{\eps,}_{\eta}\,. 
\]
Being $\mathrm{spt}\phi\Subset\Omega_{\sigma}$, we also select $0<\eps_1<\eps_0$ to be such that
for every $B\in\mathbf{B}^{\eps}$, with $B\cap\mathrm{spt}\phi\neq\emptyset$ it holds $B\Subset\Omega_\sigma$ and we relabel $\eps_0=\min\{\eps_1,\eps_2\}$.  Note that at this stage $\eps_0$ may depend on $\eta$ and $\phi$. We now deduce a preliminary estimate, see \eqref{eq:LASTIMA} below. 

In the rest of the proof  we will write $\bigcup B$ and $\sum$ without subscript to denote respectively the union and the summation taken among balls in  $\mathbf{B}^{\eps}\setminus\mathbf{B}^{\eps}_{\eta}$ that intersect $\mathrm{spt}\phi$, and we will indicate by $u_\eps^1$ and $u_\eps^2$ the first two components of the vector field $u_\eps$.
By using the identites
\begin{equation}\label{eq:identitavettoriali}
\begin{split}
&q(u_\eps)=q(u_\eps-\boldn)+\boldn(\partial_x u_\eps\times \partial_y u_\eps)\\
&\boldn(\partial_x u_\eps\times \partial_y u_\eps)=\nabla \times j(u_\eps^1,u_\eps^2)\quad \hbox{with } j(u^1_\eps,u_\eps^2)=\frac{ u_\eps^1\nabla u_\eps^2-u_\eps^2 \nabla u_\eps^1}{2}\,,
\end{split}
\end{equation}
and Stoke's theorem we get
\[
\begin{split}
\int_{\bigcup B} q(u_\eps)\phi\,\mathrm{d}x=\int_{\bigcup B}q(u_\eps-\boldn)\phi\,\mathrm{d}x-\int_{\bigcup B} j(u^1_\eps,u^2_\eps)\nabla^{\perp}\phi\,\mathrm{d}x+\sum\int_{\partial B}\phi j(u^1_\eps,u^2_\eps) \tau_{\partial B}\,\mathrm{d}\mathcal{H}^1\,.
\end{split}
\]
Similarly,
\[
\begin{split}
&\int_{\Omega_\sigma\setminus \bigcup B} q(u_\eps)\phi\,\mathrm{d}x=\\
&\int_{\Omega_\sigma\setminus \bigcup B} q(u_\eps-\boldn)\phi\,\mathrm{d}x-\int_{\Omega_\sigma\setminus \bigcup B} j(u^1_\eps,u^2_\eps)\nabla ^{\perp}\phi\,\mathrm{d}x-\sum\int_{\partial B}\phi j(u^1_\eps,u^2_\eps) \tau_{\partial B}\,\mathrm{d}\mathcal{H}^1\,.
\end{split}
\]
Notice that in the previous equalities, the sign in the last terms kept track of the boundaries orientations.  
Adding up the two previous identities we obtain:
\begin{equation}\label{eq:miviencosi}
\begin{split}
&\int_\Omega q(u_\eps)\phi\,\mathrm{d}x=\int_{\Omega_\sigma} j(u^1_\eps,u^2_\eps)\nabla ^{\perp}\phi\,\mathrm{d}x+\int_{\Omega_\sigma\setminus \bigcup B} q(u_\eps-\boldn)\phi\,\mathrm{d}x+\sum \int_{ B} q(u_\eps-\boldn)\phi\,\mathrm{d}x\,.
\end{split}
\end{equation}
Denoting by $\phi_B$ the average of $\phi$ in the ball $B$, for every $x\in B$  we have 
\[
|\phi(x)-\phi_B|\le 2r_B \|\nabla \phi\|_{L^{\infty}(\Omega)}\,,\quad |\phi_B|\le \|\phi\|_{L^\infty(\Omega)}\,.
\]
In particular, on adding and subtracting $\phi_B$ and using \eqref{eq:ballconstr2} we get
\begin{equation}\label{eq:launo}
\left|\sum \int_{B} q(u_\eps-\boldn)\phi\,\mathrm{d}x\right|\le C \|\nabla \phi\|_{L^{\infty}(\Omega)} F_\eps(u_\eps,\Omega)\sqrt{\eps} +\|\phi\|_{L^\infty(\Omega)}\sum \left|\int_{B} q(u_\eps-\boldn)\,\mathrm{d}x\right|\,.
\end{equation}
On the other hand the second term in the right hand side of \eqref{eq:miviencosi} is easily estimate thanks to the definition of $\mathbf{B^{\eps}}$:
\begin{equation*}
\left|\int_{\Omega_\sigma\setminus \bigcup B} (u_\eps-\boldn)(\partial_x u_\eps\times \partial_x u_\eps)\phi\,\mathrm{d}x\right|\le \eps^{1/4}\|\phi\|_{L^{\infty}(\Omega)}F_{\eps}(u_\eps;\Omega)\,.
\end{equation*}
Thanks to the energetic bound and up to subsequences (not relabelled), $u_\eps$ satisfies
\[
\int_\Omega|u_\eps^i|^2\,\mathrm{d}x\le \int_\Omega|u_\eps-\boldn|^2\,\mathrm{d}x\le \eps^2 F_\eps(u_\eps,\Omega)\,,\quad \hbox{for }i=1,2\,,
\]
so that, H\"older inequality gives
\begin{equation}\label{eq:latre}
\left|\int_{\Omega_\sigma} j(u^1_\eps,u^2_\eps)\nabla ^{\perp}\phi\,\mathrm{d}x\right|\le C\|\nabla \phi\|_{L^\infty(\Omega)} F_\eps(u_\eps,\Omega)\eps \,.
\end{equation}
By combining \eqref{eq:launo}-\eqref{eq:latre} with \eqref{eq:miviencosi}, we can conclude
\begin{equation}\label{eq:LASTIMA}
\begin{split}
\left|\int_\Omega q(u_\eps)\phi\,\mathrm{d}x\right|&\le C F_\eps(u_\eps,\Omega)(\|\nabla \phi\|_{L^{\infty}(\Omega)}(\eps+ \sqrt{\eps}) +\|\phi\|_{L^\infty(\Omega)}\sqrt[4]{\eps})\\&\quad+\|\phi\|_{L^\infty(\Omega)}\sum\left|\int_{B} q(u_\eps-\boldn)\,\mathrm{d}x\right|\,.
\end{split}
\end{equation}
\\

\noindent\textit{STEP 4.} We  refine \eqref{eq:LASTIMA}, by estimating the last term in the right hand side, assuming $\eta$ small enough. 
More precisely, we introduce for every $B\in \mathbf{B}^\eps\setminus \mathbf{B}^\eps_\eta$, the rescaled map $u^B_\eps \colon B_1\to \mathS^2$  defined as
\[
u^B_\eps(x)=u_\eps(r_B (x+x_B))\,.
\]
where $x_B$ denotes the center of $B$.
Notice that, conformal invariance and \eqref{eq:perbnd2} give 
\[
\int_{B_1} q(u^B_\eps-\boldn)\,\mathrm{d}x=\int_{B} q(u_\eps-\boldn)\,\mathrm{d}x\quad\hbox{ and } \int_{B_1} |\nabla u^B_\eps|^2\,\mathrm{d}x=\int_{B}|\nabla u_\eps|^2\,\mathrm{d}x\le \eta\,.
\]
Furthermore, by the very definition of $\mathbf{B}^\eps$, we have
\[
|u^B_\eps(x)-n|\le \eps^{1/4}\quad \hbox{ for every }x\in \partial B_1\,.
\]
We extend $(u^B_\eps(x)-n)$ to a map $U^B_\eps\colon B_2\to \mathR^3$ in such a way that
\[
\begin{split}
& U^B_\eps\equiv u^B_\eps-\boldn \hbox{ on }B_1\,,\\
& U^B_\eps \hbox{ constant on }\partial B_2\,,\\
& \int_{B_2} |\nabla  U^B_\eps|^2\,\mathrm{d}x\le C \int_\Omega |\nabla u^B_\eps|^2\,\mathrm{d}x\le C F_\eps(u^B_\eps, B)\,,\quad \|U^B_\eps\|_{L^{\infty}(B_2)}\le C\|u^B_\eps-\boldn\|_{L^{\infty}(\partial B_1)}\le C\eps^{1/4}\,.
\end{split}
\]
This can be done for instance by using a suitable reflection and cut-off procedure on the harmonic extension of $(u^B_\eps(x)- n)$ on $B_1$, see 
 Lemma \ref{lm:cutoff} for the details. 
By possibly taking a smaller $\eps_0$ we can define $V^B_\eps\colon B_2\to \mathS^2$ as 
\[
V^B_\eps(x)=
\begin{cases}
u_\eps^B(x) & \hbox{if }x\in B_1\,,\\
\frac{U^B_\eps(x)+\boldn}{|U^B_\eps(x)+\boldn|}  & \hbox{if }x\in A_{2,1}\,.
\end{cases}
\]
 Since $V^B_\eps$ assume a constant value on $\partial B_2$ we have that
\[
\int_{B_2} q(V^B_\eps)\,\mathrm{d}x\in4\pi\mathZ\,.
\]
On the other hand we have $\|\nabla V^B_\eps\|_{L^2(B_2)}<C F_{\eps}(u^B_\eps, B)\le C\eta$. Therefore, up to chosing $C \eta <1$ it holds
\[
\int_{B_2} q(V^B_\eps)\,\mathrm{d}x=0\,.
\]

Exploiting \eqref{eq:identitavettoriali} (with $V_\eps^B$ instead of $u_\eps$) we deduce
\[
\begin{split}
\left| \int_{B} q(u_\eps-\boldn)\,\mathrm{d}x\right|&=\left|\int_{B_1} q(u^B_\eps-\boldn)\,\mathrm{d}x\right|\le \left|\int_{B_2}q(V^B_\eps-\boldn)\,\mathrm{d}x\right|+\left|\int_{A_{2,1}}q(V^B_\eps-\boldn)\,\mathrm{d}x\right|\\
&=\left|\int_{B_2}\nabla \times j((V^B_\eps)^1, (V^B_\eps)^2)\,\mathrm{d}x\right|+\left|\int_{A_{2,1}}q(V^B_\eps-\boldn)\,\mathrm{d}x\right|\\
&=\left|\int_{\partial B_2}j((V^B_\eps)^1,(V^B_\eps)^2) \tau_{\partial B}\,\mathrm{d}\mathcal{H}^1\right|+\left|\int_{A_{2,1}}q(V^B_\eps-\boldn)\,\mathrm{d}x\right|
\\
&=\left|\int_{A_{2,1}}q(V^B_\eps-\boldn)\,\mathrm{d}x\right|\,,
\end{split}
\]
where in the last equality we used the fact that $V^B_\eps$ is constant on $\partial B_2$. Observe that easy manipulation leads to
\begin{equation*}
|V^B_\eps(x)-\boldn|\le C \eps^{1/4}\hbox{ for } x\in A_{2,1}\,,
\end{equation*}
which, plugged into the previous chain of inequalities, provides 
\begin{equation}\label{eq:LASTIMA3}
\begin{split}
\left| \int_{B} q(u_\eps-\boldn)\,\mathrm{d}x\right|&\le \left|\int_{A_{2,1}}q(V^B_\eps-\boldn)\,\mathrm{d}x \right|\\
&=\left|\int_{A_{2,1}}(V^B_\eps-\boldn)\cdot\partial_x V_\eps^B\times \partial_y V_\eps^B  \,\mathrm{d}x\right| \le C\eps^{1/4} F_{\eps}(u_\eps, B)\,.
\end{split}
\end{equation}\\

\noindent\textit{STEP 5.} By combining \eqref{eq:LASTIMA3} and \eqref{eq:LASTIMA}, and by taking the limit for $\eps\to 0$, we deduce that
\[
\left|\int \phi\, \mathrm{d}\mu\right|=\lim_{\eps\to 0}\left|\int \phi\, q(u_\eps)\, \mathrm{d}x\right|=0\,.
\]
By the arbitrariness of  $\phi\in C^{\infty}_c(\Omega\setminus \{x_1,\dots,x_{N_\eta}\})$ we conclude 
\[
|\mu|(\Omega\setminus \{x_1,\dots, x_{N_\eta}\})=0\,.
\]  
This proves \eqref{eq:LASTIMA2}.
\end{proof}

\begin{proof}[Proof of Theorem \ref{th:contmain}]
Let $(\mu_\eps)\subset\calM(\Omega)$ be such that $\mu_\eps\wto\mu$ and $\liminf_\eps\mathscr{F}_\eps(\mu_\eps)<C$, for some $C>0$. Then by Theorem \ref{th:contcomp}, $\mu$ is an atomic measure with coefficent in $4\pi\mathZ$. The $\liminf$ inequality then follows from the weak lower semicontinuity of the total variation of $\mu_\eps$ in $\Omega$ once observed that, for every $u_\eps\in\calU(\Omega)$ for which $\mu_\eps=q(u_\eps)\,dx$ one also has
\[
|\mu_\eps|(\Omega)=\int_\Omega |q(u_\eps)|\,\mathrm{d}x\le F_\eps(u_\eps; \Omega)\,,
\]
so that, by the very definition \eqref{df:scrFeps}, it holds
\[
|\mu_\eps|(\Omega)\le \mathscr{F}_\eps(\mu_\eps)\,.
\]
To prove the $\limsup$ inequality  it is  enough to adapt some computations in \cite[Lemma 3.1]{Mel}.

We first assume $\mu=4\pi\delta_{0}$. 
Let $\theta\colon\mathR^2\to\mathS^2$ be the inverse stereographic projection mapping $+\infty$ to $\boldn$, that is
\[
\theta(x)=\frac{1}{1+|x|^2}\left(2x\cdot e_1, x\cdot e_2, |x|^2-1\right)\,.
\]
A straightforward computation shows that
\begin{equation}\label{eq:stereo}
\frac 1 2 |\nabla \theta(x)|^2=\frac{4}{(1+|x|^2)^2}\,,\quad |\theta(x)-\boldn|^2=\frac{4}{1+|x|^2}\,.
\end{equation}
We let $\eps>0$ and select an auxiliary parameter $\lambda_\eps>0$ such that $\lim_{\eps\to 0}\eps\lambda_\eps=+\infty$. We set
\[
\theta_{\eps}(x):=\theta(\lambda_\eps x)\,,\quad \phi_{\eps}(x)=\boldn+\left(\frac{2\eps-|x|}{\eps}\right)\left(\theta_\eps\left(\frac{\eps x}{|x|}\right)-\boldn\right)\,,
\]
and 
\[
\Phi_{\eps}(x)=\frac{\phi_{\eps}(x)}{|\phi_{\eps}(x)|} \hbox{ for }x\in  A_{\eps,2\eps}\,.
\]
Notice that, by \eqref{eq:stereo}
\begin{equation}\label{eq:stereo2}
\left|\theta_\eps\left(\frac{\eps x}{|x|}\right)-\boldn\right|^2=\frac{4}{1+\lambda_\eps^2\eps^2}\,.
\end{equation}
Hence, for $\eps$ small enough we have
\begin{equation}\label{eq:boundsulmodulo}
|\phi_{\eps}(x)|\ge \frac 1 2 \hbox{ for }x\in  A_{\eps,2\eps}\,,
\end{equation}
which implies that $\Phi_{\eps}$ is well defined and takes values on $\mathS^2$. 
We define $u_{\eps}\in \calU(\Omega)$ as
\[
u_{\eps}(x)=
\begin{cases}
\theta_{\lambda_\eps}(x)& \hbox{ if }x\in B_\eps\,,\\
\Phi_{\eps}(x)& \hbox{ if }x\in A_{\eps,2\eps}\,,\\
\boldn &\hbox{ if }x\in \Omega\setminus B_{2\eps}\,.
\end{cases}
\]
By \eqref{eq:stereo}, \eqref{eq:stereo2} and \eqref{eq:boundsulmodulo} we have
\begin{align}\label{eq:parteBeps}
&\int_{B_\eps}|\nabla \theta_\eps|^2\,\mathrm{d}x=4\pi-\frac{4\pi}{1+\lambda_\eps^2\eps^2}\,,\quad \frac{1}{\eps^2}\int_{B_\eps}|\theta_\eps-\boldn|^2\,\mathrm{d}x=\frac{4\pi\log(1+\lambda_\eps^2\eps^2)}{\lambda_\eps^2\eps^2}\,,\\
\label{eq:parteAeps1}
&\int_{A_{\eps,2\eps}}|\nabla \Phi_{\eps}|^2\,\mathrm{d}x\le C\frac{4}{1+\lambda_\eps^2\eps^2}+C\frac{\lambda_\eps^2\eps^2}{(1+\lambda_\eps^2\eps^2)^2}\,,\\
\label{eq:parteAeps2}
& \frac{1}{\eps^2}\int_{A_{\eps,2\eps}}\left|\phi_{\eps}-\boldn\right|^2\, \mathrm{d}x\le \frac{C}{1+\lambda_\eps^2\eps^2}\,.
\end{align}
Thanks to \eqref{eq:parteBeps}-\eqref{eq:parteAeps2}  we conclude
\[
\lim_{\eps\to 0}F_\eps(u_{\eps}, \Omega)= 4\pi=|\mu|(\Omega)\,.
\]
We first check that $q(u_\eps)\,dx\wto d\delta_0$ for some $d\in\mathR$. This immediately follows from the fact that, being $\mathrm{spt} q(u_\eps)\subset B_{2\eps}$,  for every $\phi\in C^{\infty}_c(\Omega)$ such that $0\notin\mathrm{spt}\phi$, it holds that
\[
\lim_{\eps\to 0}\left(\int_{\Omega}q(u_{\eps})\phi \,\mathrm{d}x\right)=0\,.
\] 
To identify $d$ we choose $\phi\in C_c(\Omega)$ to be constant in some ball $B$ centered at the origin and $\eps$ so small that $B_{2\eps}\Subset B$. We have
\[
\int_{\Omega}q(u_{\eps})\phi\,\mathrm{d}x=\phi(0)\int_{B}q(u_{\eps})\,\mathrm{d}x=\phi(u)\mathrm{deg}(u_{\eps})=4\pi \phi(0)\,,
\] 
since the Brower degree of $u_{\eps}$ is $1$. Taking the limit for $\eps\to 0$ we deduce $d =4\pi $ as expected. 

Replacing $\theta$ with $\bar\theta\colon\mathR^2\to\mathS^2$,
$\bar \theta(x)=\frac{1}{1+|x|^2}\left(2x^{\perp}, |x|^2-1\right)$
permits to build a recovery sequence for $\mu =-4\pi\delta_0$. The full $\Gamma\mbox{-}\limsup$ result follows by a standard density argument.  
\end{proof}

\section{A discrete notion of Topological charge}\label{sec:dis1}

The goal of this section is to define a proper notion of topological charge for a spin field, this is done in Definition \ref{df:DisTopCha}. We fix throughout this subsection $\eps>0$ and $u\in\calSF_\eps(\Omega)$. 

\subsection{Admissible interpolation region}

The spin field $u$ introduces a natural selection  of  edges on the lattice $\calL_\eps$:
\[
\calN_\eps(u)=\{[i,j]\colon\ i,j\in\calL_\eps,\ |i-j|=\eps,\ u(i)=-u(j)\}\subseteq\calE_\eps\,,\quad 
\calC_\eps(u)=\calE_\eps\setminus \calN_\eps(u)\,.
\]

Notice that since it is not possible to have three mutually antipodal points on the sphere, we have that every triangle in $\calT_\eps$ contains at least one edge in  $\calC_\eps(u)$. 
We introduce some terminology.   We say that two triangles $T,T'\in\calT_\eps$ are neighbors if their intersection belongs to $\calN_\eps(u)$, and that are connected, if there exist  triangles  $T_1,\dots,T_N\in\calT_\eps$  with 
$T_1=T$ and $T_N=T'$ and such that $ T_n$ and $ T_{n+1}$ are neighbors. Similarly, a triangular region is said to be connected if it is union of pairwise connected triangles.
We call a triangular region  $R \subset \mathbb{R}^2$ an admissible interpolation region for the spin field $u$, if it is connected and maximal with respect to the set inclusion 
and we set
\[
\calR_\eps(u)=\{R\subseteq \mathR^2\colon \text{admissible interpolation region for } u\}\,.
\] 
Being $u\in\calSF_\eps(\Omega)$ we have that every $T\in \calT_\eps$ not compactly contained in $\Omega$ is an admissible interpolation region for $u$, being $u\equiv\boldn$ on the vertices of $T$. In particular, being $\Omega$ bounded, every admissible interpolation contains finitely many triangles.
Also, it is clear from the definition that every $T\in\calT_\eps$ is contained in a unique $R\in\calR_\eps(u)$.  
Given $R\in\calR_\eps(u)$ we say that $[i,j]\in \calC_\eps(u)$ is  a boundary edge of $R$ if it is the common edge of two triangles $T,T'\in\calT_\eps$, such that $T\subseteq R$ and $T'\subseteq \mathR^2\setminus \mathrm{int}(R)$.

Before proceeding further we need to show few elementary properties for the triangular regions defined above.

\begin{figure}[htp]
\center
\begin{tikzpicture}[scale=.8]
\begin{scope}[shift={(-4,0)}]
\clip(0:2.5)--++(120:2.5)--++(180:1.5)--++(240:2.5)--++(300:1.5)--++(0:2.5)--++(60:1.5);
\foreach \j in {0,...,6}{
\draw[ultra thin,gray](0:\j)++(0:-2)++(60:-3)--++(60:6);
\draw[ultra thin](0:\j)++(0:-2)++(120:-3)--++(120:6);
\draw[ultra thin](60:\j)++(60:-2)++(0:-3)--++(0:6);
}
\draw[fill=gray!20!white](0:0)--++(0:1)--++(120:1)--(0:0);
\foreach \j in {-6,...,6}{
\foreach \i in {-6,...,6}{
\draw[fill=black](0:\j)++(60:\i) circle(.02);
}
}
\draw(60:1)++(330:.5) node{$R$};
\end{scope}
\begin{scope}[shift={(2,0)}]
\clip(0:2.5)--++(120:2.5)--++(180:2.5)--++(240:2.5)--++(300:2.5)--++(0:2.5)--++(60:2.5);
\foreach \j in {0,...,6}{
\draw[ultra thin,gray](0:\j)++(0:-2)++(60:-3)--++(60:6);
\draw[ultra thin,gray](0:\j)++(0:-2)++(120:-3)--++(120:6);
\draw[ultra thin,gray](60:\j)++(60:-2)++(0:-3)--++(0:6);
}
\draw[fill=gray!20!white](0:0)++(0:1)--++(120:1)--++(180:1)--++(240:1)--++(300:1)--++(0:1)--++(60:1);
\foreach \j in {0,...,5}{
\draw[dashed](0:0)--++(60*\j:1);
}
\foreach \j in {-6,...,6}{
\foreach \i in {-6,...,6}{
\draw[fill=black](0:\j)++(60:\i) circle(.02);
}
}
\draw(0:1)++(120:1)++(330:.5) node{$R$};
\end{scope}
\begin{scope}[shift={(8,1)}]
\clip(0:-2)++(300:1)++(240:1.5)--++(0:5.5)--++(60:1.5)--++(120:2.5)--++(180:4.5)--++(240:2.5)--++(300:1.5);
\foreach \j in {-4,...,8}{
\draw[ultra thin,gray](0:\j)++(0:-2)++(60:-3)--++(60:8);
\draw[ultra thin,gray](0:\j)++(0:-2)++(120:-3)--++(120:8);
\draw[ultra thin,gray](60:\j)++(60:-2)++(0:-3)--++(0:8);
\draw[ultra thin,gray](120:\j)++(120:-2)++(0:-3)--++(0:8);
}
\draw[fill=gray!20!white](0:-2)++(300:1)--++(0:4)--++(120:1)--++(180:3)--++(240:1);
\draw[dashed](0:-1)--++(300:1)--++(60:1)--++(300:1)--++(60:1)--++(300:1)--++(60:1);
\foreach \j in {-4,...,4}{
\foreach \i in {-4,...,4}{
\draw[fill=black](0:\j)++(60:\i) circle(.02);
}
}
\draw(0:2)++(330:.5) node{$R$};
\end{scope}
\end{tikzpicture}
\caption{Examples of the different possible interpolation regions $R$ according to Lemma~\ref{lm:mistero1}. The dashed and the continuous lines denote edges in $\calN_\eps(u)$ and $\calC_\eps(u)$, respectively.}
\label{fig:regions}
\end{figure}
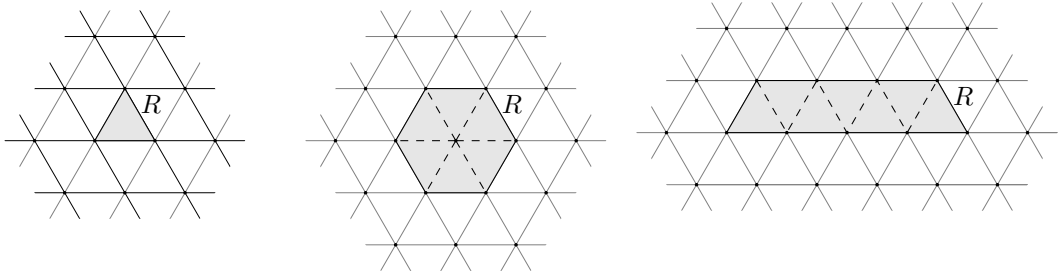

\begin{lemma}\label{lm:mistero1}
Let $R\in \calR_\eps(u)$. Then,  one of the following mutually exclusive conditions hold:
\begin{itemize}
\item [$(i)$]$R\in\calT_\eps$;
\item[$(ii)$] every triangle $T\in\calT_\eps$, $T\subset R$ has two edges in $\calN_\eps(u)$;
\item[$(iii)$] there are exactly two triangles $T\in\calT_\eps$, $T\subset R$ with two edges in $\calC_\eps(u)$, all the others triangles have two edges in $\calN_\eps(u)$.
\end{itemize} 
\end{lemma}

\begin{proof}
Suppose that $R$ is not a triangle, so that $(i)$ does not hold. Replacing every triangle $T$ with a vertex and every edge $\calN_\eps(u)$ with an arc, we obtain a connected planar graph (see Figure \ref{fig:auxiliary_graph}) with finitely many arcs  in which every vertex has at most degree two (being the degree of a vertex  the number of edges that are incident to that vertex). For such a graph it is easy to show, by induction, that the set of vertices with degree one is either empty or contains two elements. If the first case occurs then $(ii)$ holds, otherwise $(iii)$ holds.

\begin{figure}[htp]
\begin{tikzpicture}
\begin{scope}
\clip(0:-2)++(300:1)++(240:1.5)--++(0:4.5)--++(60:3.5)--++(120:2.5)--++(180:2.5)--++(240:2)--++(180:1)--++(240:2.5)--++(300:1.5);
\foreach \j in {-4,...,8}{
\draw[ultra thin,gray](0:\j)++(0:-2)++(60:-3)--++(60:8);
\draw[ultra thin,gray](0:\j)++(0:-2)++(120:-3)--++(120:8);
\draw[ultra thin,gray](60:\j)++(60:-2)++(0:-3)--++(0:8);
\draw[ultra thin,gray](120:\j)++(120:-2)++(0:-3)--++(0:8);
}
\draw[fill=gray!20!white](0:-2)++(300:1)--++(0:3)--++(60:1)--++(120:1)--++(0:1)--++(120:1)--++(180:1)--++(240:1)--++(300:1)--++(180:2)--++(240:1);
\draw[dashed](0:-1)--++(300:1)--++(60:1)--++(300:1)--++(60:1)--++(300:1)--++(60:1)--++(180:1)--++(60:1)--++(180:1)--++(60:1)--++(300:1)--++(60:1);
\foreach \j in {-4,...,4}{
\foreach \i in {-4,...,4}{
\draw[fill=black](0:\j)++(60:\i) circle(.02);
}
}
\draw(0:2)++(330:.5) node{$R$};
\draw[fill=black](240:1)++(180:1)++(0:1/3)++(60:1/3)circle(.05)++(0:1)circle(.05)++(0:1)circle(.05)++(60:1)circle(.05)++(120:1)circle(.05)++(0:1)circle(.05);
\draw[fill=black](240:1)++(180:1)++(0:1)++(60:1/3)++(120:1/3)circle(.05)++(0:1)circle(.05)++(0:1)circle(.05)++(120:1)circle(.05)++(60:1)circle(.05);
\draw[thick](240:1)++(180:1)++(0:1/3)++(60:1/3)--($(240:1)+(180:1)+(0:1)+(60:1/3)+(120:1/3)$)--($(0:1)+(240:1)+(180:1)+(0:1/3)+(60:1/3)$)--($(0:1)+(240:1)+(180:1)+(0:1)+(60:1/3)+(120:1/3)$)--($(0:2)+(240:1)+(180:1)+(0:1/3)+(60:1/3)$)--($(0:2)+(240:1)+(180:1)+(0:1)+(60:1/3)+(120:1/3)$)--($(0:2)+(240:1)+(180:1)+(0:1/3)+(60:1/3)+(60:1)$)--($(0:2)+(240:1)+(180:1)+(0:1)+(60:1/3)+(120:1/3)+(120:1)$)--($(0:2)+(240:1)+(180:1)+(0:1/3)+(60:1/3)+(60:1)+(120:1)$)--($(0:2)+(240:1)+(180:1)+(0:1)+(60:1/3)+(120:1/3)+(120:1)+(60:1)$)--($(0:3)+(240:1)+(180:1)+(0:1/3)+(60:1/3)+(60:1)+(120:1)$);
\end{scope}
\end{tikzpicture}
\caption{Construction of the planar graph in the proof of Lemma \ref{lm:mistero1} (same notation as in Figure \ref{fig:regions}).}
\label{fig:auxiliary_graph}
\end{figure}
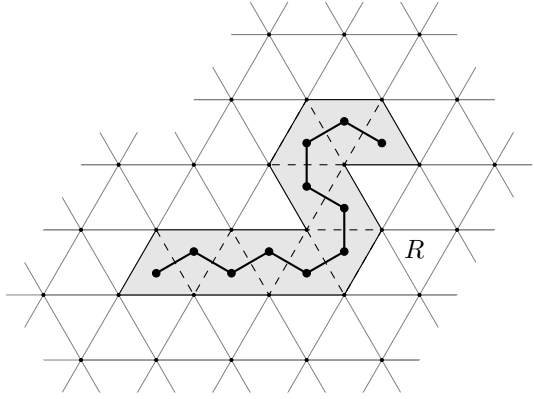
\end{proof}

\begin{lemma}\label{lm:mistero2}
Let $R\in \calR_\eps(u)$.

Then, the set
\[
\{u(i)\colon\ i\in \calL_\eps,\, i\in\Gamma \text{ for some }\Gamma\in\calN_\eps(u),\, \Gamma\subset R \}\,.
\]
is either empty or contains exactly two antipodal points of $\mathS^2$.

\end{lemma}

\begin{proof}
Suppose that $R$ is not a triangle, since otherwise the statement is trivial. We consider the same planar graph as defined in the proof of Lemma \ref{lm:mistero1} above. For every pair of adjancent arcs of such a graph the corresponding edges in $\calN_\eps$ must intersect in some point of $\calL_\eps(R)$. By exploiting this property, toghether with the connectedness of the graph and the definition of $\calN_\eps(u)$, it is easy to conclude.
\end{proof}

\begin{lemma}\label{lm:mistero3}
Let $R\in \calR_\eps(u)$. Then $u(\calL_\eps(R))$ contains at most four points. 
\end{lemma}

\begin{proof}
It follows immediately by combining Lemma \ref{lm:mistero1} and Lemma \ref{lm:mistero2}.
\end{proof}

\begin{lemma}\label{lm:mistero4}
Let $R\in\calR_\eps(u)$ and $T\in\calT_\eps$, $T\subset R$. The the following facts hold:
\begin{itemize}
\item [(a)] if $T$ has a single edge in $\calN_\eps(u)$, then it contains at least a boundary edge of $R$;
\item [(b)]  if $T$ has a single edge in $\calN_\eps(u)$ and $u(\calL_\eps(R))$ contains four points, then $T$ contains two boundary edges of $R$.
\end{itemize} 
\end{lemma}

\begin{proof}
We prove (a). To fix the notation we suppose $T=[i,j,k]$ and $\Gamma_1=[i,j]$, $ \Gamma_2=[i,k]$,  $u(j)=-u(k)$ and $u(i)\notin\{u(j),u(k)\}$. We need to show that at least one between $\Gamma_1$ or $\Gamma_2$ is a boundary edge of $R$. 
Assume by contradiction that this is not the case. This means that given $T_1, T_2\in\calT_\eps$ such that
\[
T_1\cap T=\Gamma_1\,,\quad T_2\cap T=\Gamma_2\,,
\]
we have $T_1,T_2\in\calT_\eps(R)$. 
By Lemma \ref{lm:mistero1} either $T_1$ contains exactly two edges in $\calN_\eps(u)$ or $T_2$ does. In any case $i$ is the endpoint of some edge in $\calN_\eps(u)$, contained in $R$. By Lemma \ref{lm:mistero2} this implies $u(i)\in\{u(j),u(k)\}$ which is a contradiction. 

Let us prove (b). With the same notation as above and by (a), we can assume that, for instance, $T_1\subseteq \mathR^2\setminus \mathrm{int} R$. Then, suppose by contradiction $T_2\subseteq R$. Arguing as above we deduce that $T_2$ cannot contains exactly two edges of $\calN_\eps(u)$. On the other hand, suppose that it contains exactly one edge of $\calN_\eps(u)$. Then by Lemma \ref{lm:mistero1}, $T$ and $T_2$ are the only two triangles having only one edge in $\calN_\eps(u)$. By Lemma \ref{lm:mistero2}, this implies that $u(\calL_\eps(R))$ contains three elements, which is not the case. Therefore, $T_2$ contains no edges of $\calN_\eps(u)$. But this contradicts $T\subset R$. Hence, the lemma is proved. 
\end{proof}

\begin{lemma}\label{lm:misteroaggiunto}
Let $R\in\calR_\eps(u)$ and assume  that $u(\calL_\eps(R))$ contains four points.  Suppose $i,j\in\calL_\eps$ are such that $[i,j]\subset R$ and $u(i)\notin\{u(j),-u(j)\}$, then $[i,j]$ is a boundary edge for $R$.
\end{lemma}

\begin{proof}
Let $k\in\calL_\eps$ be such that $[i,j,k]\subset R$. Notice that it cannot be $[i,k]\in\calN_\eps(u)$ and $[j,k]\in\calN_\eps(u)$ since otherwise we would have $u(i)\in\{u(j),-u(j)\}$. The conclusion then follows by point (b) in Lemma \ref{lm:mistero4}.
\end{proof}

\subsection{Admissible interpolations surfaces} 

For every $R\in\calR_\eps(u)$, we consider the cone generated by the set $u(\calL_\eps(R))$, that is:
\[
\mathrm{cone}_R(u)=\left\{\sum_{w\in u(\calL_\eps(R))}\lambda_w w\colon \lambda_w\geq 0\right\}\,.
\]
We call $S_R(u)\subseteq \mathS^2$ an admissible interpolation surface for $u$ in $R$ if one of the next conditions is satisfied. If $\mathrm{cone}_R(u)$  is not a linear subspace of $\mathR^3$, we require
\begin{equation}\label{df:SR1}
S_R(u)=\left\{\frac{\sum_{w\in u(\calL_\eps(R))}\lambda_w w}{\left|\sum_{w\in u(\calL_\eps(R))|}\lambda_w w\right|}\colon\ \lambda_w\ge 0,\, \sum\limits_{w\in u(\calL_\eps(R))} \lambda_w>0\right\}\,.
\end{equation}
If $\mathrm{cone}_R(u)$  is a linear subspace of $\mathR^3$, we require 
\begin{equation}\label{df:SR2}
S_R(u)=\left\{\frac{\sum_{w\in u(\calL_\eps(R))}\lambda_w w+\tau h}{\left|\sum_{w\in u(\calL_\eps(R))|}\lambda_w w+\tau h\right|}\colon\ \lambda_w,\tau\ge 0,\, \sum\limits_{w\in u(\calL_\eps(R))} \lambda_w+\tau>0\right\}\,.
\end{equation}
for some $h\in\mathS^2$ orthogonal to $\mathrm{cone}_R(u)$.\\

The next lemma deals with the elementary characterization of the admissible interpolation surfaces.

\begin{lemma}\label{rem:geometrico}
Let $R\in\calR_\eps(u)$ and $S_R(u)\subseteq \mathS^2$ be an admissible interpolation surface for $u$ in $R$. Then $S_R(u)$ is either a semi great-circle or an arc of geodesic or a geodesic triangle or the union of two adjacent geodesic triangle or an hemisphere.
\end{lemma}

\begin{proof}
Let $N$ be the number of points in $u(\calL_\eps(R))$.  
Assume $N=1$, then $S_R(u)$ is a single point. 
Assume $N=2$. Here we distinguish two cases. Either $u(\calL_\eps(R))=\{p,-p\}$ for some $p\in\mathS^2$.  Then $S_R(u)$ is determined as in \eqref{df:SR2}, for some $h\in\mathS^2$, such that $h\cdot p=0$, and it is straightforward to verify that
\[
\begin{split}
&S_R(u)=\left\{\frac{\lambda p+\tau h}{|\lambda p+\tau h|}\colon  \tau,\lambda\ge 0,\, \tau+\lambda>0\right\}\cup\left\{\frac{-\lambda p+\tau h}{|-\lambda p+\tau h|}\colon  \tau,\lambda\ge 0,\,\tau+\lambda>0\right\}\,.
\end{split}
\]
This provides
\[
S_R(u)=\gamma_{p,h}\cup\gamma_{-p,h}\,,
\]
that is $S_R(u)$ is a semi great-circle.
Otherwise, $u(\calL_\eps(R))=\{p,q\}$ for some $p,q\in \mathbb{S}^2$ with $p\neq -q$. In this case, $S_R(u)$ is an arc of geoedesic, precisely:
\begin{align*}
S_R(u) = \gamma_{p,q}\,.
\end{align*} 

Assume  $N=3$, then $u(\calL_\eps(R))=\{p_1,p_2,p_3\}$ for some $p_n\in\mathS^2$. If the $p_n$'s are linearly independent then \eqref{df:SR1} holds, and $S_R(u)$ is uniquely determined as the geodesic triangle with vertices $p_1,p_2,p_3$, that is:
\[
S_R(u)=S(p_1,p_2,p_3)\,.
\] 
If the $p_n$'s are not linearly independent, then we distinguish again two cases: either $\mathrm{cone}_R(u)$ is an hyperplane or it is not. In the first case $S_R(u)$ is an hemisphere centered at $h\in\mathS^2$, where $h$ is orthogonal to $\mathrm{cone}_R(u)$, otherwise we have that $S_R(u)$ is a semi great-circle or a geodesic arc, specifically: 
\begin{align*}
S_R(u) =\begin{cases} \gamma_{p,q} \cup \gamma_{-p,q} &\text{if } u(\calL_\eps(R))=\{-p,p,q\}\,,\\
\gamma_{p_i,p_j} &\text{otherwise, where } p_i\cdot p_j \leq p_k\cdot p_l \text{ for all } l\neq k\,.
\end{cases}
\end{align*} 
Recalling Lemma \ref{lm:mistero3}, it remains only to discuss the case where $u(\calL_\eps(R))=\{p_1,p_2,p_3, p_4\}$ for some $p_n\in\mathS^2$. By definition of $R$ we can assume that $p_1=-p_2$. Then, if $\mathrm{cone}_R(u)$ is an hyperplane, $S_R(u)$  is an hemisphere centered at $h\in\mathS^2$, where $h$ is orthogonal to $\mathrm{cone}_R(u)$. If $\mathrm{cone}_R(u)$ is not an hyperplane, but it lies on an hyperplane then $S_R(u)$ is the semi great circle connecting $p_1$ to $-p_1$ and passing through $p_3$ and $p_4$. At last, if $\mathrm{cone}_R(u)$ does not lie on a hyperplane, then $p_1,p_3,p_4$ are linearly independent as well as $p_2,p_3,p_4$. In this case \eqref{df:SR1} holds, and being $p_1=-p_2$, it holds 
\[
\begin{split}
S_R(u) &=\left\{\frac{\lambda_1 p_1+\lambda_3 p_3+\lambda_4 p_4}{| \lambda_1p_1+\lambda_3 p_3+\lambda_4 p_4|}\colon  
\, \lambda_i\ge 0,\quad \lambda_1+\lambda_3+\lambda_4>0
\right\} \\&\quad
\cup \left\{\frac{\lambda_2 p_2+\lambda_3 p_3+\lambda_4 p_4}{| \lambda_2p_2+\lambda_3 p_3+\lambda_4 p_4|}\colon  
\, \lambda_i\ge 0,\quad \lambda_2+\lambda_3+\lambda_4>0
\right\}\,,
\end{split}
\] 
which  gives
\[
S_R(u)=S(p_1,p_3,p_4)\cup S(-p_1,p_3,p_4)\,.
\]
that is $S_R(u)$ is the union of two adjacent geodesic triangles.
\end{proof}

\begin{figure}[htp]
 \def\svgwidth{14.8cm}
   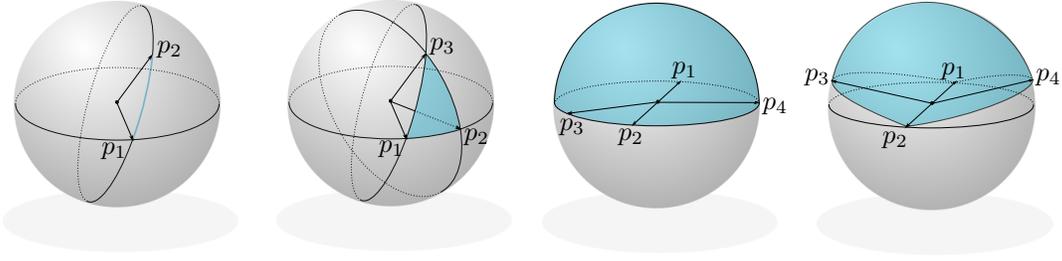
  \caption{The region $S_R(u)$ in the cases $N=2,3,4$.}
\end{figure}

In view of Lemma \ref{rem:geometrico} above and  recalling  the discussion in Section \ref{sec:bas}, we have that the relative interior of every admissible surface $S_R(u)$ for $u$ in $R$ with $\mathcal{H}^2(S_R(u))>0$ is a simple region on $\mathS^2$. In fact, the spin field $u$ induce a natural orientation on $\partial S_R(u)$ as we show below.

\begin{figure}[htp]
\center
\begin{tikzpicture}
\tikzset{>={Latex[width=1mm,length=1mm]}};
\draw[ultra thin,gray](0:0)--++(300:1)--++(60:1)--++(300:1)--++(60:1)--++(180:1)--++(60:1)--++(300:1);
\draw[ultra thin,gray](0:0)--++(300:1)--++(0:1);
\draw[fill=gray!20!white,gray!20!white](0:0)--++(0:1)--++(60:1)--++(180:1)--++(240:1);
 \draw(0:0)--++(0:1)--++(60:1);
  \draw[dashed](0:0)++(0:1)--++(120:1);
   \draw[fill=black](60:1)circle(.02);
     \draw[fill=black](300:1)circle(.02)++(0:1)circle(.02)++(60:1)circle(.02); 
  \draw[fill=black](0:0) circle(.02)++(0:1)circle(.02)++(60:1)circle(.02);
  \draw[->](0:1)++(60:.5) arc (60:180:.5);
  \draw(0:0) node[anchor=north]{$i$};
    \draw(0:1) node[anchor=north]{$j$};
        \draw(0:1)++(60:1) node[anchor=west]{$k$};
\draw(60:.75)++(0:.125)++(30:.5)++(270:.125) node{$R$};
\end{tikzpicture}
\caption{The edge $(i,j)_R$ precedes the edge $(j,k)_R$.}
\label{fig:predecessor}
\end{figure}
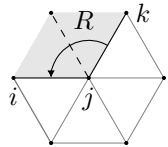

We find it convenient convenient to introduce some more terminology. Let $R\in\calR_\eps(u)$ and $\Gamma=[i,j]\in\calC_\eps(u)$ a boundary edge of $R$. We say that $[i,j]$  is oriented from $i$ towards $j$, and we write $[i,j]=(i,j)_R$, if the lattice point $k\in\calL_\eps$  obtained by rotating counterclockwise the edge $[i,j]$ around $i$ is such that the triangle $[i,j,k]$ lies in $R$ (see Figure \ref{fig:predecessor}). Otherwise we say that it is oriented from $j$ to $i$ and we write $[i,j]=(j,i)_R$.  
We define the geodesic curve
\begin{equation}\label{df:latticegeodesic}
\gamma_{\Gamma, R, \, u}:[0,1]\to \mathS^2,\quad \gamma_{\Gamma, R, \, u}(t)=
\begin{cases}
\gamma_{u(i),u(j)}(t) &\hbox{ if  } \Gamma=(i,j)_R\,,\\
\gamma_{u(j),u(i)}(t) & \hbox{ if } \Gamma=(j,i)_R\,.
\end{cases}
\end{equation}

 Given two oriented boundary edges of $R$, $(i,j)_R$ and $(j,k)_R$ we say that $(i,j)_R$ precedes $(j,k)_R$, and we write $(i,j)_R\prec_R(j,k)_R$, if $[i,j]$ can be obtained by rotating in the counterclockwise direction $[j,k]$ around $j$ without intersecting $\mathR^2\setminus R$. It is then easy to verify that for every $\Gamma$ boundary edge  of $R$ there exist unique $\Gamma_1$ $\Gamma_2$ boundary edges of $R$ such that $\Gamma_1\prec_R\Gamma\prec_R\Gamma_2$, and that it is not possible to have $\Gamma\prec_R \Gamma$.

\begin{lemma}\label{lm:misteriosissimo}
Let $R\in\calR_\eps(u)$ and $S_R(u)\subseteq \mathS^2$ be an admissible interpolation surface for $u$ in $R$. Assume $\mathcal{H}^2(S_R(u))>0$. Then, the relative interior of $S_R(u)$ is a simple region on $\mathS^2$. Moreover there exist $\Gamma_1,\dots,\Gamma_N\in\calC_\eps(u)$ boundary edges of $R$ with $N=3,4$, and
\[
\Gamma_1\prec\dots\prec\Gamma_N\prec\Gamma_1\,,
\] 
such that, the curve 
\[
\begin{split}
&\gamma_{\, R, \, u}\colon[0,N]\to \mathS^2\,,\\ 
&\gamma_{\, R, \, u}(t)=\gamma_{\Gamma_{m+1}\,,\,  R,\, u}(t-m)\,, \quad \hbox{ \text{if }}t\in [m,m+1]\,, \quad 0\le m\le N-1\,.
\end{split}
\]
is a parametrization of the relative boundary of $S_R(u)$. 
\end{lemma}

\begin{proof}
The first part of the lemma can be deduced by Lemma \ref{rem:geometrico} as already noticed. 
In particular either $S_R(u)$ is a geodesic triangle or the union of two adjacent geodesic triangles or an hemisphere. Let $N$ be the number of points in $u(\calL_\eps(R))$. We can assume that either $N=3$ or $N=4$, since otherwise the condition $\mathcal{H}^2(S_R(u))>0$ is violated (see the proof of Lemma \ref{rem:geometrico}).

We discuss the first case, in particular either $S_R(u)$ is a geodesic triangle or an hemisphere and 
\begin{equation*}
u(\calL_\eps(R))=\{w_1,w_2,w_3\}\,,\text{ with }w_n\neq -w_m \text{ for } n,m \in \{1,2,3\}\,.
\end{equation*}
The latter implies that $R\in\calT_\eps$ (since it does not contains edges in $\calN_\eps(u)$). In particular $R=[i,j,k]$, for some lattices points $i,j,k\in\calL_\eps$ and $[i,j]$, $[j,k]$, $[i,k]$ are all boundary edges of $R$. Up to relabelling indices we can assume  $[i,j]=(i,j)_R$, $[j,k]=(j,k)_R$, $[i,k]=(k,i)_R$. By setting $\Gamma_1=[i,j]$, $\Gamma_2=[j,k]$ and $\Gamma_3=[i,k]$, we  conclude.

Assume instead that $N=4$. By referring again to the proof of Lemma \ref{rem:geometrico} we deduce that
\begin{equation*}
u(\calL_\eps(R))=\{w_1,w_2,w, -w\}\,.
\end{equation*}
with $w_n \neq \pm w$ for $n\in \{1,2\}$ and that  $S_R(u)$ is either the union of the geodesic triangle with veritces $S(w,w_1,w_2)$ and $S(-w,w_1,w_2)$ when $\cone_R(u)$ is not a linear subspace, or  an hemisphere when $\cone_R(u)$ is an hyperplane. 

By definition of $R$,  there exist $i,j,k\in\calL_\eps$ such that $[i,j,k]\subset R$ and $[i,k]\in\calN_\eps(u)$, with $u(j)=w_1$, $u(i)=w$ and $u(k)=-w$.  Exploiting Lemma \ref{lm:mistero4} we deduce that $[i,j]$ and $[j,k]$ are boundary edges of $R$. Similarly, there exist $i',j',k'\in\calL_\eps$ such that $[i',j',k']\subset R$, $u(j')=w_2$,  $u(i')=w$ and $u(k')=-w$, and we can infer that $[i',j']$ and $[j',k']$ are boundary edges of $R$ as well.
Moreover, Lemma \ref{lm:mistero2} gives $\{u(i), u(k)\}=\{u(i'), u(k')\}$, while combining Lemma \ref{lm:mistero1} and Lemma \ref{lm:mistero2} we get $u(\mathcal{L}_\varepsilon(R)) \setminus \{u(j), u(j')\}=\{-w,w\}$. 
Up to relabelling indices, we can assume that $[i,j]=(i,j)_R$, $[j,k]=(j,k)_R$ and $[i',j']=(i',j')_R$, $[j',k']=(j',k')_R$.
\par\noindent Now, we consider the sequence $(\Gamma_k)$ of boundary edges, inductively defined as
\[
\Gamma_1=[j,k]\,,\quad  \Gamma_n\prec_R \Gamma_{n+1}\,.
\]
Since $R$ contains a finite union of boundary edges, it has to be $\Gamma_{n_0}\prec_{R} \Gamma_1$, for some $n_0>0$, that is $\Gamma_{n_0}=(i,j)_R$. On the other hand $u$ does not take antipodal values on the endpoints of each $\Gamma_n$. Hence, there exists $n_1\le n_0-1$ such that  $\Gamma_{n_1}=(m,l)_R$,  $u(\{m,l\})=\{u(i), q\}$  and $q\notin\{ u(k), u(i)\}$. This forces $q=u(j')$ and so $\Gamma_{n_1}=(j',k')_R$. As a consequence $\Gamma_{n_1-1}=(i',j')_R$. 
In conclusion the following order holds true:
\[
\begin{split}
\Gamma_1&=(j,k)_R\prec_R\dots\prec_R\Gamma_{n_1-1}=(i',j')_R\prec_R\Gamma_{n_1}\\&=(j',k')_R\prec_R\dots\prec_R\Gamma_{n_0}=(i,j)_R\prec_R(j,k)_R=\Gamma_1\,.
\end{split}
\]
Notice that for every $1\le n\le n_0$, $n\notin \{1,n_1-1,n_1,n_0\}$ the set  $u(\calL_\eps(\Gamma_n))$ is a singleton. Therefore, by setting
\[
\Gamma_{n_1-1}=\Gamma_2\,,\quad \Gamma_{n_1}=\Gamma_3\,,\quad \Gamma_{n_0}=\Gamma_4\,,
\]
we get
\[
\gamma_{R,u}([0,4])=\gamma_{w_1w}\cup\gamma_{ww_2}\cup\gamma_{-ww_2}\cup\gamma_{-ww_1},
\]
which proves the lemma.
\end{proof}

\subsection{Discrete topological charge}  This subsection is devoted to the definition of discrete topological charge and to prove some properties of the spherical extension of a spin field having such charge. 
\begin{definition}[Discrete topological charge]\label{df:DisTopCha}
Let $\eps>0$ and $u\in\calSF_\eps(\Omega)$. An atomic measure $\mu\in\calM(\Omega)$ is called an admissible topological charge for the spin field $u$, briefly $\mu\in\calM_{\text{adm}}(\Omega; u),$ if
\[
\mu=\sum_{R\in\calR_\eps(u)}c_R\,\delta_{x_R}\,,
\]
and the following conditions holds:
\begin{itemize}
\item[(i)]  $x_R\in\mathrm{int}(R)$;
\item[(ii)] $c_R=\sigma(\gamma_{R,u}, S_R(u))\mathcal{H}^2(S_R(u))$, where $S_R(u)$ is an admissible interpolation surface for $u$ in $R$,  $\sigma(\gamma_R, S_R(u))\in \{1,-1\}$ and whenever $\mathcal{H}^2(S_R(u))>0$ it holds $\sigma(\gamma_R, S_R(u))=1$ if the curve $\gamma_{R,u}$ given in Lemma \ref{lm:misteriosissimo} is positively oriented with respect to $S_R(u)$, while $\sigma(\gamma_{R, u}, S_R(u))=-1$ otherwise.
\end{itemize}
\end{definition}

In what follows we show that if $\mu\in \calM_{\text{adm}}(\Omega, u)$ we can construct an interpolation  vector field $\overarc u\in\calU(\Omega)$, such that $\overarc u=u$ on $\calL_\eps$ and for which  $q(\overarc u)dx$ is close in flat norm to $\mu$.

\begin{theorem}\label{theo:interpolation}
Let $\Omega\subset\mathR^2$ be a bounded open set, $u\in\calSF_\eps(\Omega)$ and $\mu\in\calM_{\text{adm}}(\Omega; u)$, be given, according to Definition \ref{df:DisTopCha}, by
\[
\,\quad \mu=\sum_{R\in\calR_\eps(u)}\sigma(\gamma_{R, u}, S_R(u)) \mathcal{H}^2(S_R(u))\delta_{x_R}\,.
\]
Then, there exists $\overarc u\in\calU(\Omega)$, such that $\overarc u=u$ on $\calL_\eps$ and
\begin{equation}\label{eq:intFE}
\frac{1}{2}\int_\Omega |\nabla \overarc u|^2\,\mathrm{d}x\leq C_1 H_\varepsilon(u,\Omega)\,, \quad \text{and} \quad \frac{1}{\eps^2}\int_\Omega |\overarc u-\boldn|^2\,\mathrm{d}x \le C_1 \, Z_\varepsilon(u,\Omega)\,,
\end{equation}
where $C_1>0$ is a universal constant. 
Moreover, for every $R\in\calR_\eps(u)$ we have
\begin{equation}\label{eq:intcarica}
\int_R q(\overarc u)\,\mathrm{d}x=\sigma(\gamma_R, S_R(u))\mathcal{H}^2(S_R(u))\,.
\end{equation}
Finally, $\mu(\mathbb{R}^2) \in 4\pi\mathbb{Z}$ and
\begin{equation}\label{eq:intflat}
\|q(\overarc u)dx-\mu\|_{\mathrm{flat}, \Omega}\le C_2\, \eps  (H_\eps(u,\Omega)+1)H_\eps(u,\Omega)\,,
\end{equation}
where $C_2>0$ is a universal constant. 
\end{theorem}

\begin{proof}
 For the sake of readability we split the proof into different steps; the symbol $C$
denotes a universal constant, which may depend
on $\Omega$  and whose value may change from line to line. \\

\noindent\textit{STEP 1.}
First, we define a linear interpolation  $v:\mathR^2\to \mathR^3$ of the spin field $u$. For every edge $\Gamma=[i,j]\in\calC_\eps(u)$ we define $u_\Gamma\in\mathS^2$ as
\[
u_\Gamma=\frac{u(i)+u(j)}{|u(i)+u(j)|}\,.
\]
Note that $u_\Gamma$ is well defined being $\Gamma\in\calC_\eps(u)$ and that
it holds
\begin{equation}\label{eq:intLm}
u_\Gamma\cdot u(i)\ge 0\hbox{ and  }u_\Gamma\cdot u(j)\ge 0\,.
\end{equation}
Then, for every $R\in\calR_\eps(u)$, we select $u_R\in S_R(u)$ in such a way that
\begin{equation}\label{eq:intLc}
u_R\cdot z\ge 0,\hbox{ for every }z\in S_R(u)\,.
\end{equation}
This can be done by exploiting Lemma \ref{lm:costruzioneu_R}.
We consider an auxiliary lattice of points $\calL'_\eps\subset\mathR^2$ (see Figure \ref{figura})
\[
\begin{split}
\calL_\eps'=\calL_\eps\cup\calL^c_{\eps}\cup \calL^m_{\eps}\,,
\end{split}
\]
where
\[
\calL^c_{\eps}=\{x_T: T\in\calT_\eps\}\,,\quad \calL^m_{\eps}=\{x_\Gamma:\ \Gamma\in\calE_\eps\}\,.
\]
and above $x_T$ and $x_\Gamma$ denote respectively the center of the triangle $T$ and the midpoint of the edge $\Gamma$.
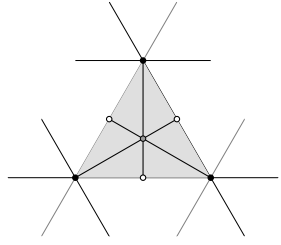
\begin{figure}[htp]
\begin{center}
\begin{tikzpicture}[scale=1.8]
\begin{scope}[shift={(-4,0)}]
\clip(0:1.5)--++(120:1.5)--++(180:.5)--++(240:1.5)--++(300:.5)--++(0:1.5)--++(60:.5);
\foreach \j in {0,...,6}{
\draw[ultra thin,gray](0:\j)++(0:-2)++(60:-3)--++(60:6);
\draw[ultra thin](0:\j)++(0:-2)++(120:-3)--++(120:6);
\draw[ultra thin](60:\j)++(60:-2)++(0:-3)--++(0:6);
}
\fill[gray!25]
($(60:1/3)+(0:1/3)$)--($(60:1/2)+(0:1/2)$)--(60:1)-- cycle;
\fill[gray!25]
($(60:1/3)+(0:1/3)$)--(60:1/2)--(60:1)-- cycle;
\fill[gray!25]
(0:0)--($(60:1/3)+(0:1/3)$)--(60:1/2)-- cycle;
\fill[gray!25]
(0:0)--($(60:1/3)+(0:1/3)$)--(0:1/2)-- cycle;
\fill[gray!25]
(0:1/2)--(0:1)--($(60:1/3)+(0:1/3)$)--cycle;
\fill[gray!25]
(0:1)--($(60:1/3)+(0:1/3)$)--($(60:1/2)+(0:1/2)$)--cycle;
\foreach \j in {0,...,2}{
\draw[thin](0:0)++($(60:1/3)+(0:1/3)$)--++($(120+\j*120:1/3)+(60+\j*120:1/3)$);
\draw[thin](0:0)++($(60:1/3)+(0:1/3)$)--++($(60+\j*120:1/6)+(\j*120:1/6)$);
\draw[thin, fill=white](0:0)++($(60:1/3)+(0:1/3)$)++($(60+\j*120:1/6)+(\j*120:1/6)$) circle(.02);
}
\draw[ultra thin, fill=gray!80](0:0)++($(60:1/3)+(0:1/3)$) circle(.02);
\foreach \j in {-6,...,6}{
\foreach \i in {-6,...,6}{
\draw[fill=black](0:\j)++(60:\i) circle(.02);
}
}
\end{scope}
\end{tikzpicture}
\end{center}
\caption{Subdivision of an elementary cell of ${\mathcal L}_\eps$ into cells of ${\mathcal L'_\eps}$. The points of ${\mathcal L^m_\eps}$ are indicated in white, while those in ${\mathcal L^c_\eps}$ are indicated in gray.} \label{figura}
\end{figure}
Correspondingly, we define the family of triangles $\calT'_\eps$ subordinated to $\calL'_\eps$ as
\[
\calT'_\eps=\left\{[i,j,k]:\, i\in\calL_\eps,\, j\in \calL^c_{\eps},\, k\in \calL^m_{\eps},\,|i-k|=\frac{\eps}{2}\,, |i-j|=\frac{\eps}{\sqrt{3}}\,, |j-k|=\frac{\eps}{2\sqrt{3}}\right\}\,.
\]
We extend, without relabelling it, the spin field $u$ to a map $u:\calL_\eps'\to\mathS^2$  defined by setting  for every $i\in\calL_\eps'$
\begin{equation}\label{eq:intext}
u(i)=
\begin{cases}
u(i) & \text{ if }i\in\calL_\eps\,,\\
u_{R} & \text{ if }i\in\calL^c_{\eps},\, i\in R,\, \text{ with }R\in\calR_\eps(u)\,,\\
u_{\Gamma}& \text{ if }i\in\calL^m_{\eps},\,i\in \Gamma,\,\text{ and }  \Gamma\in\calC_\eps\,,\\
u_{R} & \text{ if }i\in\calL^m_{\eps},\, i\in\Gamma,\,\text{ and } \Gamma\in\calN_\eps\,.
\end{cases}
\end{equation} 
Note that taking advantage of \eqref{eq:intLm} and \eqref{eq:intLc}, for every $[i_1,i_2,i_3]\in\calT'_\eps$, with $i_1\in\calL_\eps,\, i_2\in \calL^c_{\eps},\, i_3\in \calL^m_{\eps}$,  we have
\begin{equation}\label{eq:leonard}
u(i_n)\cdot u(i_m)\ge 0
\end{equation}
for every $n,m=1,2,3$. 
We define $v:\mathR^2\to \mathR^3$ to be the continuous and piecewise affine map defined in every $T\in\calT'_\eps$ as the unique affine map that satisfies $v\equiv u$ on the vertices of $T$. 
For every $x\in \mathR^2$ we have
\begin{equation}\label{eq:boundbasso}
|v(x)|\ge \frac{1}{\sqrt{3}}\,.
\end{equation} 
Indeed, for every $T=[i_1,i_2,i_3]\in\calT'_\eps$, with $i_1\in\calL_\eps,\, i_2\in \calL^c_{\eps},\, i_3\in \calL^m_{\eps}$ and for every $x\in T$ it holds
\[
v(x)=\sum_{n=1}^3\lambda_{n} u(i_n)\,,\]
where $\lambda_1,\lambda_2,\lambda_3\in (0,1)$ and $\lambda_1+\lambda_2+\lambda_3=1$. In particular,  Jensen's inequality and \eqref{eq:leonard} provide
\[
\begin{split}
&|v(x)|^2=\sum_n\lambda^2_{n}+\sum_{n\neq m}\lambda_n\lambda_mu(i_n)\cdot u(i_m)\ge\frac{1}{3}+\sum_{n\neq m}\lambda_n\lambda_mu(i_n)\cdot u(i_m)\ge\frac 1 3 \,,
\end{split}
\]
which is \eqref{eq:boundbasso}. Moreover, we note that
\begin{equation}\label{eq:boundalto}
|\nabla v(x)|\le C\eps^{-1}\text{ for }\, x\in \mathR^2\,,
\end{equation}
which follows since $v$ is a piecewise affine interpolation on the lattice $\calL'_\eps$ of values in $\mathS^2$. 
\par\noindent Now, let us show that the following hold:
\begin{equation}\label{eq:boundperv}
\frac{1}{2}\int_\Omega |\nabla v|^2\,\mathrm{d} x\leq C_1 H_\varepsilon(u,\Omega)\,,\quad \text{and} \quad\frac{1}{\eps^2}\int_\Omega |v-\boldn|^2\,\mathrm{d} x \le C_1 \, Z_\eps(u,\Omega)\,.
\end{equation}
With this aim let $R\in\calR_\eps(u)$. We distinguish two cases. 
Suppose that for some $T=[i,j,k]\in\calT_\eps$, $T\subset R$ there holds $u(i)\cdot u(j)\le 0$. Then, an explicit computation shows
$1 \le H_\varepsilon(u,T)$ and $\frac{1}{2}\leq Z_\varepsilon(u,T)$. Summing over all triangles $T \subset R$, using \eqref{eq:boundalto} we get
\begin{equation}\label{eq:intuglyR}
\frac{1}{2}\int_R |\nabla v|^2\,\mathrm{d}x \leq C\, H_\eps(u, R)\,,  \quad \frac{1}{\varepsilon^2}\int_R |v-\boldn|^2\,\mathrm{d}x \leq C\, Z_\varepsilon(u,R)\,.
\end{equation} 
Suppose instead that $R\in\calR_\eps$ and that for every $T=[i,j,k]\in\calT_\eps$, $T\subset R$ it holds $u(i)\cdot u(j)> 0$. Then $R=[i_1,i_2,i_3]$
for some $i_1,i_2,i_3\in\calL_\eps$ and $\cone_R(u)$ is not a linear subspace of $\mathR^3$.
Accordingly, for $T'\in\calT_\eps'$, $T'\subset T$, given as
\[
T'=\left[i_1, j_1, k_1\right]\,, \text{ where } j_1=\frac{i_1}{2}+\frac{i_2}{2}\,, \,k_1=\frac{i_1}{3}+\frac{i_2}{3}+\frac{i_3}{3}\,,
\]
taking into account of \eqref{eq:intext}, we have
\[
u(j_1)=\frac{u(i_1)+u(i_2)}{|u(i_1)+u(i_2)|}\,,\quad u(k_1)=\frac{\sum_{n=1}^3 \mu_n u(i_n)}{\left|\sum_{n=1}^3 \mu_n u(i_n)\right|}\,, \text{ where }  0<\mu_1,\mu_2,\mu_3<1 \text{ and }\sum_{n=1}^3\mu_n=1\,.
\]
An elementary computation shows that
\[
|u(i_1)-u(j_1)|\le C\left|u(i_1)-u(i_2)\right|\,,\quad 
|u(i_1)-u(k_1)|\le C(|u(i_1)-u(i_2)| +|u(i_1)-u(i_3)|)\,.
\]
The latter implies
\[
|\nabla v(x)|^2\le \frac{C}{\eps^2}\left(|u(i_1)-u(j_1)|^2+|u(i_1)-u(k_1)|^2+|u(j_1)-u(k_1)|^2\right)\le \frac{C}{\eps^2}\, H_\varepsilon(u, T)\,.
\]
for every $x\in T$. On the other hand 
\[
|v(x)-\boldn|^2\le C(|u(i_1)-\boldn|^2+|u(j_1)-\boldn|^2+|u(k_1)-\boldn|^2)\le C\, Z_\varepsilon(u, T)\,.
\]
These inequalities provide
\begin{equation}\label{eq:intgoodR}
\frac{1}{2}\int_R |\nabla v|^2\,\mathrm{d}x \leq C\, H_\varepsilon(u,R)\,,\quad \text{and} \quad \frac{1}{\eps^2}\int_R |v-\boldn|^2\mathrm{d}x\le C\, Z_\eps(u, R)\,.
\end{equation}
By combining \eqref{eq:intuglyR} and \eqref{eq:intgoodR} we get \eqref{eq:boundperv}.\\

\noindent\textit{STEP 2.}
We define
$\overarc u:\mathR^2\to\mathS^2$ as  $\overarc u=\frac{v}{|v|}$.
Taking into account \eqref{eq:boundbasso}, \eqref{eq:boundalto} and \eqref{eq:boundperv} we have that $\overarc u$ is a  well defined Lipschitz continuous map, $\overarc u\in\calU(\Omega)$ and \eqref{eq:intFE} holds. Furthermore, we have that $\overarc u$ restricted to every $R\in\calR_\eps(u)$ satisfies
\begin{equation}\label{eq:intSR}
\overarc u(R)=S_R(u)\,,\text{ and }\, \#\{x \in R\colon u(x) = y\}=1\, \text{ for } \mathcal{H}^2\text{-a.e. } y \in S_R(u)\,.
\end{equation}
as it can be deduced by the definition of $S_R(u)$ (we refer the reader to Lemma \ref{lm:uHat} for a detailed proof of these latter elementary identities).

Next we prove \eqref{eq:intcarica}. Let $R\in\calR_\eps(u)$. When $\mathcal{H}^2(S_R(u))=0$, formula \eqref{eq:intcarica} follows at once by using \eqref{eq:intSR} and the area formula \eqref{eq:area}.
Hence, we can assume $\mathcal{H}^2(S_R(u))>0$.
Notice that for every $A\in\mathR^{3\times 2}$ and $b\in\mathR^3$, one has 
\begin{equation*}
q(Ax+b)=b\cdot(A^1\times A^2)\,\text{ and }\, q\left(\frac{Ax+b}{|Ax+b|}\right)=\frac{b\cdot(A^1\times A^2)}{|Ax+b|^3}\,\text{ if }\, |Ax+b|>0\,.
\end{equation*}
Thus, for every $T\in\calT'_\eps$, $q(\overarc u)$ is either identically null or has constant sign. 
In fact, we claim the following
\begin{equation}\label{eq:claimgeom}
\text{for every }T'\in \calT_\eps',\, T'\subset R \text{ either } q(\overarc u)\equiv 0 \text{ in }T' \text{ or } \mathrm{sign}(q(\overarc u))\equiv\sigma(\gamma_R, R) \text{ in } T'\,. 
\end{equation} 
To prove the claim, assume $\sigma(\gamma_R, S_R(u))$ positive (the other case being similar) and $T'\in\calT'_\eps(R)$. By construction, there is a unique $\Gamma\in\calE_\eps$ such that $\Gamma=[i,j]$ contains one edge of $T'$. We thus have three possibilities. 
If $\Gamma\in\calN_\eps(u)$, taking into account \eqref{eq:intext}, we note that $v$ interpolates in $T'$ between only two values. This implies $q(v)\equiv 0$ on $T'$, hence $q(\overarc u)\equiv 0$ on $T'$ as well. Similarly, if $\Gamma\in\calC_\eps(u)$ and $u(i)=u(j)$, we deduce $q(\overarc u)\equiv 0$ on $T'$. Suppose instead that $\Gamma\in\calC_\eps(u)$ and $u(i)\neq u(j)$. Notice that being $\mathcal{H}^2(S_R(u))>0$ either $R$ is a triangle or $u(\calL_\eps(R))$ contains four points (see Lemma \ref{rem:geometrico}). Therefore, by Lemma 
\ref{lm:misteroaggiunto}, $\Gamma$ is a boundary edge of $R$. Using the notation from the previous subsection we suppose for instance that $\Gamma=(i,j)_R$, and we let
\[
\tau=\frac{j-i}{|j-i|}\,.
\]
Let $x$ be a point in the relative interior of $\Gamma$ with $x\in T'$. For $\delta>0$ small enough and $i=1,2$ we consider the curves $\alpha_i:(0,\delta)\to R$ and $\beta_i:(0,\delta)\to \mathbb{S}^2$ defined as
\[
\begin{split}
\alpha_1(s)=x+s \tau\,,\, \alpha_2(s)=x+s \tau^\perp\,,\\
\beta_i(s)=\overarc u\circ \alpha_i(s)\,\qquad\text{for }i=1,2\,.
\end{split}
\]
Note that the curve $\beta_1$ is a geodesic and in particular there exists $t\in (0,1)$ and $\lambda>0$ such that 
\[
\gamma_{\Gamma, R, u}(t)=\overarc u(x)\,,\qquad \gamma_{\Gamma, R, u}'(t)=\lambda \nabla \overarc u(x)\tau\,,
\]
being $\gamma_{\Gamma, R, u}$ defined as in \eqref{df:latticegeodesic}.
Note moreover that the curve $\beta_2$ takes values in $S_R(u)$ and satisfies
\[
\beta(0)=\overarc u(x)\,,\quad \beta'(0)=\nabla \overarc u(x)\tau^{\perp}\,. 
\]
Recalling our assumption on the coefficent $\sigma(\gamma_R, S_R)$ we have
\[
0\le (\gamma_{\Gamma, R, u}(t_x)\times \gamma_{\Gamma, R, u}'(t_x))\cdot \beta(0)=(\overarc u(x)\times\lambda\nabla \overarc u(x)\tau)\cdot \nabla \overarc u(x)\tau^\perp=\lambda|\tau|^2q(\overarc u)(x)\,.
\]
By the arbitrariness of $x$,  $q(\overarc u)$ is positive on $\Gamma\cap T'$ and from what noticed above, it is positive on $T'$. Thus, \eqref{eq:claimgeom} is proved.  
\par\noindent The identity \eqref{eq:intcarica} is now a simple consequence of \eqref{eq:intSR} and \eqref{eq:area}: 
\[
\begin{split}
\mathcal{H}^2(S_R(u))=\int_R|q(\overarc u)|\,\mathrm{d}x&=\sum_{T\in\calT'_\eps,\, T\subseteq R} \int_T \mathrm{sign}(q(\overarc u))q(\overarc u)\,\mathrm{d}x\\
&=\sum_{T\in\calT'_\eps,\, T\subseteq R} \sigma(\gamma_R, R)\int_T q(\overarc u)\,\mathrm{d}x=\sigma(\gamma_R, R) \int_R q(\overarc u)\,\mathrm{d}x\,.
\end{split}
\]
This implies also that $\mu(\mathbb{R}^2)\in 4\pi\mathbb{Z}$. \\

\noindent\textit{STEP 3.} At last we show \eqref{eq:intflat}. 
Let $\phi\in\mathrm{Lip}(\Omega)$ with $\|\nabla \phi\|_{L^\infty(\Omega)}\le 1$. From the previous step of the proof we have
\[
\begin{split}
\left|\int_\Omega \phi\, q(\overarc u)\,\mathrm{d}x-\int_\Omega \phi\, \mathrm{d}\mu\right|&\le \sum_{R\in\calR_\eps(u)}\left|\int_R \phi\, q(\overarc u)\,\mathrm{d}x-\sigma(\gamma_R, S_R)\mathcal{H}^2(S_R)\phi(x_R)\right|\\
&=\sum_{R\in\calR_\eps(u)}\left|\int_R q(\overarc u)(x)(\phi(x)-\phi(x_R))\,\mathrm{d}x\right|\\
& \le \sum_{R\in\calR_\eps(u)}   \diam(R)\|\nabla \phi\|_{L^{\infty}(\Omega)}\int_R|\nabla \overarc u|^2\, \mathrm{d}x\,.
\end{split}
\] 
Since  $\mathrm{diam}(R)\le C \eps (\# \{\Gamma\in \calN_\eps:\, \Gamma\subset R\}+1)$ and $\#\{\Gamma\in \calN_\eps:\, \Gamma\subset R\}\le H_\eps(u, R)$
we get \eqref{eq:intflat}.
\end{proof}

\section{Concentration of discrete topological charge}\label{sec:dis2}

The aim of this section is to prove Theorem \ref{th:discomp} and Theorem \ref{th:main}. Therefore, for a given $\Omega \subset \mathbb{R}^2$ bounded and regular open set, we consider sequences of spin fields $(u_\eps)\in\calSF_\eps(\Omega)$ such that
\begin{equation*}
\sup_{\eps>0} E_\eps(u_\eps,\Omega)<\infty\,.
\end{equation*}
Within the proofs of the results of this section we use the symbol $C$ to denote some universal constant $C>0$ which may depend on $\Omega$ and on $\sup_{\eps>0}E_{\eps}(u_\eps, \Omega)$ and whose value may change from line to line. 
We start with the proof of the compactness result, which is an easy consequence of Theorem \ref{theo:interpolation} and Theorem \ref{th:contcomp}.

\begin{proof}[Proof of Theorem \ref{th:discomp}]
Let $\overarc u_\eps\in\calU(\Omega)$ be defined from $u_\varepsilon$ as in Theorem \ref{theo:interpolation}. By \eqref{eq:intFE}, we have
\[
\sup_{\eps>0} F_\eps(\overarc u_\eps, \Omega)<\infty\,.
\] 
By Theorem \ref{th:contcomp}, up to subsequences, $q(\overarc u_\eps)dx\wto \mu$ as $\eps\to 0$, where   $\mu\in\calM(\Omega)$ is an atomic measure with coefficents in $4\pi\mathZ$.
By taking the limit as $\eps\to 0$ in \eqref{eq:intflat} written for $\overarc u=\overarc u_\eps$ and $\mu=\mu_\eps$, we deduce that, up to subsequence,  $\mu_\eps\wto \mu$ as well. This concludes the proof.  
\end{proof}
We define $\phi_\eps :(3\eps,+\infty)\times \mathZ\times \mathR^2\to [0,+\infty]$ as
\begin{equation}\label{df:cell2}
\phi_\eps(s,d,x)=\inf\left\{ E_{\eps}(w, B_s(x))\colon\, w\in\calSF_\eps,\,\text{s.t.} \quad \begin{aligned}
& w\equiv n, \text{ in }\calL_\eps(\mathR^2\setminus B_{s-3\eps}(x))\,,\\
& \zeta(B_s(x))=d\,, \text{ for some }\zeta\in\calM_{adm}(\Omega; w)
\end{aligned}
\right\}\,.
\end{equation}

In order to prove Theorem \ref{th:main} we begin with the following Lemma.

\begin{lemma}\label{lm:cell}
$\phi_\eps(s,d,x)$ does not depend on $x$. 
Moreover the limit $\lim_{\eps\to 0}\phi_\eps(s,d,x)$ exists, and it does not depend on $s$.
\end{lemma}

\begin{proof}
First, we notice  that $\phi_\eps(T,d,x)\le \phi_\varepsilon(t,d,x)$, for every $T>t>3\eps$. On the other hand, for every $x\in\mathR^2$ let us set $a^x_1,a^x_2\in\mathR$ to be such that $x=a^x_1e_1+a^x_2\hat e_2$, where $\hat e_2$ is given as in \eqref{df:lattice} and set $x_\eps\in\calL_\eps$ to be
\[
x_\eps=\eps\left\lfloor\frac{a^x_1}{\eps}\right\rfloor e_1+\eps\left\lfloor\frac{a^x_2}{\eps}\right\rfloor \hat e_2\,.
\]
Notice that $|x-x_\eps|<2\eps$. In particular we get that, every spin field which is admissible for $\phi_\eps(s,d,x)$ is admissible also for $\phi_\varepsilon(s+2\eps,d,x_\eps)$, and viceversa every admissible test spin field for $\phi_\eps(s,d,x_\eps)$ is admissible for $\phi_\varepsilon(s+2\eps,d,x)$. This and the aforementioned monotonicity imply:
\[
\phi_\eps(s,d,x)\le \phi_\varepsilon(s+2\eps,d,x_\eps)\le \phi_\eps(s,d,x_\eps)\,,
\]
and 
\[
\phi_\eps(s,d,x_\eps)\le \phi_\varepsilon(s+2\eps,d,x)\le \phi_\eps(s,d,x)\,,
\]
A simple translation argument also shows that 
\[
\phi_\eps(s,d,x_\eps)=\phi_\eps(s,d,0)\,.
\]
Hence, for every $s>0$ and $x\in\mathR^2$, it holds
\begin{equation*}
\phi_\eps(s,d,x)= \phi_\eps(s,d,0)\,,
\end{equation*}
proving that $\phi_\eps(s,d,x)$ does not depend on $x$. 

From now on, we assume $x=0$. Now we select $(\eps_n^1)$ and $(\eps_n^2)$ vanishing sequences such that
\[
\lim_{n\to\infty}\phi_{\eps_n^1}(s,d,0)=\liminf_{\eps\to 0}\phi_{\eps}(s,d,0)\,,\quad \lim_{n\to\infty}\phi_{\eps_n^2}(s,d,0)=\limsup_{\eps\to 0}\phi_{\eps}(s,d,0)\,, 
\]
and 
\[
\eps_n^1 \gg\eps_n^2\,.
\]
Let $w_{n,1}\in\calSF_{\eps_n^1}$ be admissible for $\phi_{\eps_n^1}(s,d,0)$. We consider the rescaled spin field $w_{n,2}\in\calSF_{\eps_n^2}$, defined by
\[
w_{n,2}(k)=w_n\left(\frac{\eps_n^1}{\eps_n^2} k\right).
\]
 For every $k\in \calL_\eps(B_{s-3\eps_n^2})$ we have $w_{n,2}(k)=\boldn$, being
\[
\left|\frac{\eps_n^1}{\eps_n^2} k\right|\ge \frac{\eps_n^1}{\eps_n^2}s-3\eps_n^1\ge s-3\eps_n^2\,.
\]
Moreover
\[
E_{\eps_n^1}(w_{n,1}, B_{s})=E_{\eps_n^2}\left(w_{n,2}, B_{\frac{\eps_n^1}{\eps_n^2}s}\right)\,,
\]
and, given $\zeta_{n,1}\in\calM_{\text{adm}}(B_{s},w_{n,1})$, with $\zeta_{n,1}(B_s)=d$,  the measure $\zeta_{n,2}\in\calM(\Omega)$, defined by setting 
\[
\zeta_{n,2}(E)=\zeta_{n,1}\left(\frac{\eps_n^1}{\eps_n^2}E\right)
\]
for every $E\subset \mathR^2$, is such that 
\[\zeta_{n,2}\in\calM_{\text{adm}}\left(B_{\frac{\eps_n^1}{\eps_n^2}s},w_{n,2}\right)\,,\text{ and } \zeta_{n,2}\left(B_\frac{\eps_n^1}{\eps_n^2}\right)=d\,.\]
Hence, $w_{n,2}$ is admissible for the minimization problem $\phi_{\eps_n^2}\left(\frac{\eps_n^1}{\eps_n^2}s, d, 0\right)$,
Therefore
\[
\phi_{\eps_n^1}(s,d,0)\ge \phi_{\eps_n^2}\left(\frac{\eps_n^1}{\eps_n^2}s, d, 0\right)\,.
\]
By taking into account of the monotonicity with respect to the first argument we deduce
\[
\phi_{\eps_n^1}(s,d,0)\ge \phi_{\eps_n^2}\left(s, d, 0\right)\,.
\]
By passing to the limit as $n\to +\infty$, the latter entails
\[
\liminf_{\eps\to 0}\phi_{\eps}(s,d,0)\ge \limsup_{\eps\to 0}\phi_{\eps}(s,d,0)\,.
\]
This proves the existence of $\lim_{\eps\to 0}\phi_\eps(s,d,0)$.
To conclude the proof we are left to show that 
\[
\lim_{\eps\to 0}\phi_{\eps}(t,d,0)=\lim_{\eps\to 0}\phi_\eps(T,d,0)\,,
\]
for every $t,T>0$. 
Assume without loss of generality that $t<T$. We can use the same rescaling argument as above. More precisely,  for every $w\in\calSF_{\eps}$ which is admissible  for $\phi_\eps(T,d,x)$
we define $v\in\calSF_{\frac{t}{T}\eps}$ by letting
\[
v(i)=w\left(\frac{T}{t}\,i\right)\,.
\]
As above $v$ is admissible  for the minimization problem $\phi_{\frac{t}{T}\eps} (t, d,\frac{t}{T} x)$, and
\[
E_{\frac{t}{T}\eps}\left(v, B_t\left(\frac{t}{T} x\right)\right)=E_{\eps}(w, B_T(x))\,.
\] 
The latter in turn entails
\[
\phi_{\frac{t}{T}\eps} \left(t, d,\frac{t}{T} x\right)\le \phi_{\eps} (T, d, x)\,.
\]
Now, exploiting the independence from the third argument, we get
\[
\phi_{\frac{t}{T}\eps} (t, d,x)\le  \phi_{\eps} (T, d, x)\,.
\]
By taking the limit as $\eps\to 0$ on both sides of the previous inequality we deduce
\[
\lim_{\eps\to 0}\phi_{\eps} (t, d,x)\le  \lim_{\eps\to 0}\phi_{\eps} (T, d, x)\,.
\]
The conclusion follows observing that the opposite inequality is a consequence of the monotonicity of $\phi_\eps(\cdot,d,x)$. 
\end{proof}

We are now in a position to  prove Theorem \ref{th:main}.  For reader convenience we recall that $\mathscr{E}_\eps:\calM(\Omega)\to [0,\infty]$ are defined (see \ref{df.scrEeps}) as
\[
\mathscr{E}_\eps(\mu)=\inf\{E_\eps(u):\ u\in\calSF_\eps,\,  \mu\in\calM_{adm}(\Omega; u) \}\,,
\]
where the convention $\inf\emptyset=\infty$ is adopted. Furthermore, we notice that the existence of the limit in \eqref{df:cell} is a consequence of Lemma \ref{lm:cell} above, which also implies
\[
\psi(d)=\lim_{\eps\to 0}\phi_{\eps}(s,d,x)\,,
\]
for every $s>0$ and $x\in\mathR^2$, being $\phi_{\eps}(s,d,x)$ defined as in \eqref{df:cell2}.

\begin{proof}[Proof of Theorem \ref{th:main}]
\textit{$\Gamma\text{-}\liminf$ {\rm inequality}.} Let $(\mu_\eps)\subset\calM(\Omega)$ be such that $\mu_\eps\wto\mu $ and assume $\mathscr{E}_\eps(\mu_\eps)<\infty$, since otherwise there is nothing to prove. This implies that we can select, for every $\eps>0$, a spin field $u_\eps\in\calSF_\eps(\Omega)$  such that 
\[
\mu_\eps\in\calM_{adm}(\Omega; u_\eps)\,,\quad \mathscr{E}_\eps(\mu_\eps)\leq E_\eps(u_\eps)\leq \mathscr{E}_\eps(\mu_\eps)+\eps\,.
\]
The first condition above also provides the existence of points $\{x_R\}_{R\in\calR_\eps(u)}$ and a selection of admissible interpolation regions for $u_\eps$, $\{S_R(u_\eps)\}_{R\in\calR_\eps(u)}$ such that
\begin{equation}\label{eq:mainmu}
\mu_\eps=\sum_{R\in\calR(u_\eps)}\sigma(\gamma_{R, u_\eps}, S_{R}(u_\eps))\mathcal{H}^2(S_R(u_\eps))\delta_{x_R}\,,
\end{equation}
according to Definition \ref{df:DisTopCha}.
Furthermore, Theorem \ref{th:discomp} gives that for $N>0$, and for some subset $\{x_1,\dots, x_N\}\subset\Omega$ and $d_1,,\dots, d_N\in\mathZ$ it holds:
\[
\mu=4\pi d_1\delta_{x_1}+\dots+4\pi d_{N}\delta_{x_N}\,.
\]

Let $\nu_\eps\in\calM(\Omega)$ be defined as 
\[
\nu_\eps(G)=E_\varepsilon(u_\eps, G) \text{ for every measurable set }G\subset\mathR^2\,.
\] 
Up to subsequences it is not restrictive to suppose that  $\nu_\eps\wto\nu\in\calM(\Omega)$. We fix $R>0$ to be such that:
\begin{equation}\label{eq:mainR1}
R<\min \{|x_n-x_m|: n,m\in\{1\dots, N\}, n\neq m\}\,.
\end{equation}
and
\begin{equation}\label{eq:mainR2}
\nu(\partial B_R(x_n))=0\,,\quad \text{ for }n\in\{1,\dots, N\}\,.
\end{equation}
Notice that thanks to \eqref{eq:mainmu} and \eqref{eq:mainR1} it holds $\mu_{\eps\, \mres B_R(x_n)}\wto 4\pi d_n \delta_{x_n}$ for every $n=1,\dots, N$.
We now  fix $x_n\in \{x_1,\dots, x_N\}$. Our goal is to construct a spin field $w_{\eps,n}\in\calSF_{\eps}$ which is admissible for the minimization problem $\phi_\eps(R,d_n,x_n)$, and whose energy $E_{\eps}(w_{\eps,n}, B_R(x_n))$ is close to $E_{\eps}(u_\eps, B_R(x_n))$.
In order not to burden the notation we suppose $x_n=0$, and we do not write the dependence on $n$, when it is does not generate confusion. For every $0<\delta< \frac 1 2$ we let $r$, $0<r<R$ be such that
\begin{equation}\label{eq:mainr1}
\nu(\partial B_r)=0
\end{equation}
and 
\begin{equation}\label{eq:mainr2}
\nu( \bar A_{R,r})\le \delta\,.
\end{equation}
Combining \eqref{eq:mainR2}, \eqref{eq:mainr1} and \eqref{eq:mainr2} we deduce
\begin{equation}\label{eq:mainRr}
\lim_{\eps\to 0}E_\eps(u_\eps, \bar A_{R,r} )\le \delta\,.
\end{equation}
We define two additional parameters $r_1<r_2$ as 
\begin{equation*}
r_1=r+\frac{R-r}{3}\,,\quad r_2=r+\frac{2(R-r)}{3}\,,
\end{equation*}
and consider a radial cut-off funtion $h\in C^{\infty}_c(B_{r_2})$ such that
\[
0\le h\le 1 \text{ for every }x\in B_{r_2}\,, \quad h=1 \text{ for every  }\in B_{r_1}\,,\quad \|\nabla h\|_{L^{\infty}(B_{r_2})}\le \frac{3}{R-r}\,.
\]
Note that $A_{r_2,r_1}\subset\joinrel\subset A_{R,r}$. Moreover for $\eps<\frac{R-r}{3}$ we have that if $i\in\calL_\eps(\bar A_{r_2, r_1})$ then all the nearest neighborhoods of $i$ lie in $A_{R,r}$. 
This means that for every $i\in\calL_\eps(\bar A_{r_2, r_1})$ it holds
\[
|u_\eps(i)-\boldn|^2\le \sum_{\{T\in\calT_\eps, i\in T\}}E_\eps(u_\eps, T)\le E_\eps(u_\eps, \overline{A}_{R, r})\,.
\]
For $\eps<\eps_0=\eps_0(R, r, \delta)$, by exploiting \eqref{eq:mainRr}, the latter implies
\[
|u_\eps(i)-\boldn|^2\le\delta\,, \text{ for every }i\in\calL_\eps(\bar A_{r_2, r_1})\,,
\]
so that
\begin{equation}\label{eq:maincutoff0}
|h(i)u_\eps(i)+(1-h(i))\boldn|\ge 1-\sqrt{\delta}>\frac{1}{2}\,,\quad  \text{ for every }i\in\calL_\eps(\bar A_{r_2, r_1})\,.
\end{equation}
This justify the following definition: $w_\eps:\calL_\eps\to\mathS^2$
\[
w_\eps (i)=
\begin{cases}
u_\eps(i) & \text{ if } i\in  \calL_\eps( B_{r_1})\,,\\
\frac{h(i)u_\eps(i)+(1-h(i))\boldn}{|h(i)u_\eps(i)+(1-h(i))\boldn|} & \text{ if } i\in  \calL_\eps( \bar A_{r_2,r_1})\,,\\
\boldn & \text{ if } i\in \calL_\eps(\mathR^2\setminus \bar B_{r_2})\,.
\end{cases}
\]
Notice that
\begin{equation}\label{eq:maratona}
|w_\eps(i)-\boldn|^2
\le 4|u_\eps(i)-\boldn|^2\,\quad \text{for every }i\in\calL_\eps\,.
\end{equation} 
This is trivial whenever $i\in\calL_\eps(B_{r_1})$ or $i\in\calL_\eps(\mathR^2\setminus B_{r_2})$, while if $i\in\calL_\eps( \bar A_{r_2,r_1})$ it is a consequence of the Lipschitzianity of the projection to the sphere $\mathrm{Proj}_{\mathS^2}$ far away from the origin and \eqref{eq:maincutoff0} as follows
\[
\begin{split}
|w_\eps(i)-\boldn|^2&=|\mathrm{Proj}_{\mathS^2}(h(i)u_\eps(i)+(1-h(i))\boldn)-\mathrm{Proj}_{\mathS^2}\boldn|^2\\&\le 4|(h(i)u_\eps(i)+(1-h(i))\boldn)-\boldn|^2\le 4|(h(i)(u_\eps(i)-\boldn)|^2\le  4|u_\eps(i)-\boldn|^2\,.
\end{split}
\]
From \eqref{eq:maratona} we also deduce 
\begin{equation}\label{eq:maratona2}
|w_\eps(i)-w_\eps(j)|^2\le 2 |u_\eps(i)-\boldn|^2+2|u_\eps(j)-\boldn|^2\, \text{ for every }i,j\in\calL_\eps\,.
\end{equation}
By combining \eqref{eq:maratona}, \eqref{eq:maratona2} and \eqref{eq:mainRr} we get 
\begin{equation}\label{eq:mainRrw}
\lim_{\varepsilon \to 0} E_\varepsilon(w_\varepsilon,A_{R,r}) <C\delta\,.
\end{equation}
and
\begin{equation}\label{eq:civuolesotto}
\lim_{\eps\to 0}E_{\eps}(w_\eps, B_{R})\le \lim_{\eps\to 0}E_{\eps}(u_\eps, B_r)+C\delta\,.
\end{equation}

Now, the inequalities \eqref{eq:mainRr} and \eqref{eq:mainRrw} in particular show that, for $\eps<\eps_0(R, r, \delta)$ whenever $R\in\calR_\varepsilon(u_\eps)$, $R \subset\joinrel\subset B_R$, is  such that for some $T\in\calT_\eps(R)$ it holds $T\subset\joinrel\subset A_{R,r}$,  we have that 
\[
R\in \calT_\eps\,,\quad R\in\calR_\eps(w_\eps)\,,
\]
and also that $S_R(u_\eps)$ and $S_R(w_\eps)$ are both uniquely defined. 
Thus, thanks to the fact that $u_\eps$ and $w_\eps$ coincides $B_{r_1}$, we deduce that
\begin{equation*}
\{R\in \calR_\eps(u_\eps)\colon R\subset\joinrel\subset B_R(x_i)\}=\{R\in \calR_\eps(w_\eps): R\subset\joinrel\subset B_R(x_i)\}\,.
\end{equation*}

We define the following family of Radon measures $\zeta_\eps$:
\[
\zeta_\eps:=\sum_{T\in\calR_\eps(w_\eps), T\cap B^c_{r_1}\neq\emptyset} \sigma(\gamma_T, S_T(w_\eps))\mathcal{H}^2(S_T(w_\eps))\delta_{x_T}+\sum_{R\in\calR_\eps(u_\eps),\, R\subset\joinrel\subset B_{r_1}, R\subset\joinrel\subset B_R}\mu_{\eps\, \mres R}\,.
\]
By construction it holds  $\zeta_\eps\in\calM_{adm}(\Omega; w_\eps)$ (recall that $w_\eps\equiv \boldn$ outside $B_{r_2}$). We claim that 
\begin{equation}\label{eq:claim}
\zeta_\eps\wto d_n\delta_{0}\,\text{ as }\eps\to 0\,.
\end{equation} 
To prove the claim we evoke again Theorem \ref{th:discomp} to get that, as $\eps\to 0$ and up to subsequences, it holds
\[
\zeta_\eps\wto 4\pi \sum_{m=1}^M a_m\delta_{y_m}\,,
\]
for some $M>0$, $\{y_1,\dots, y_M\}\subset\Omega$ and $\{a_1,\dots, a_M\}\in\mathZ\setminus \{0\}$. Assume by contradiction that, for some $m=1,\dots, M$ it holds $y_m\in \Omega\setminus B_{r_1}$, and choose $\varphi\in C^{\infty}_c(B_s(y_m))$, with $\varphi(y_m)=1$ and $s>0$  such that $B_s(y_m)\subset\joinrel\subset \mathR^2\setminus B_{r_1}$. 
Let $\overarc w_\eps, \overarc u_\eps\in\calU(\Omega)$ defined as in Theorem \ref{theo:interpolation} and set, with a slight abuse of notation \[A_{\eps, r_1}=\bigcup \{T:\, T\in\calR_\eps(w_\eps);\  T\cap B^c_{r_1}\neq\emptyset\}\,.\]
 We have
\[
\begin{split}
&\left\|\zeta_{\eps\mres A_{\eps, r_1}}-\mu_{\eps\mres A_{\eps, r_1}}\right\|_{\mathrm{flat},\, A_{R,r}}\\
&\le \left\|\zeta_{\eps\mres A_{\eps, r_1}}-q(\overarc w_\eps)dx_{\mres A_{\eps, r_1}}\right\|_{\mathrm{flat},\, A_{R,r}}+\left\|\mu_{\eps\mres A_{\eps, r_1}}-q(\overarc u_\eps)dx_{\mres A_{\eps, r_1}}\right\|_{\mathrm{flat},\, A_{R,r}}\\
&\quad +\|q(\overarc w_\eps)dx_{\mres A_{\eps, r_1}}-q(\overarc u_\eps)dx_{\mres A_{\eps, r_1}}\|_{\mathrm{flat},\, A_{R,r}}
\end{split}
\]
 and therefore by using \eqref{eq:intflat}, \eqref{eq:mainRr} and \eqref{eq:mainRrw} we deduce
\[
\begin{split}
\left\|\zeta_{\eps\mres A_{\eps, r_1}}-\mu_{\eps\mres A_{\eps, r_1}}\right\|_{\mathrm{flat},\, A_{R,r}}
\le C\delta\,.
\end{split}
\]
Thus
\[
|a_m|\le \lim_{\eps\to 0}\left|\int_{B_s}\varphi\, \mathrm{d}\zeta_\eps\right|\le \lim_{\eps\to 0}\left|\int_{B_s}\varphi\, \mathrm{d} \mu_\eps\right|+C\delta\le C\delta\,.
\] 
Being $a_m\in\mathZ$, the latter implies $a_m=0$. Therefore every $y_m$ has to lie into $B_{r_1}$. On the other hand $\zeta_{\eps\, \mres B_{r_1}}=\mu_{\eps\, \mres B_{r_1}}$, hence \eqref{eq:claim} follows.
By exploiting \eqref{eq:intcarica} we can also show that $\zeta_\eps(B_R)=d$ for $\eps$ small enough. In conclusion the spin field $w_{\eps,n}\in\calSF_{\eps}$ is admissible for the minimization problem $\phi(R, d_n, x_n)$. 
Finally, localizing the construction for every $B_R(x_n)$, we have
\[
\begin{split}
\liminf_{\eps\to 0}\mathscr{E_\eps}(\mu_\eps)&\ge 
\liminf_{\eps\to 0} E_{\eps}(u_\eps,\Omega)\ge \sum_{n=1}^N\lim_{\eps\to 0} E_{\eps}(u_\eps,B_R(x_n))\\
&\ge \sum_{n=1}^N\liminf_{\eps\to 0} E_{\eps}(w_{\eps, n},B_R(x_n)+\sum_{n=1}^N\liminf_{\eps\to 0}\left( E_{\eps}(u_\eps,B_R(x_n))-E_{\eps}(w_{\eps, n},B_R(x_n))\right)\\
&\ge \sum_{n=1}^N\liminf_{\eps\to 0} E_{\eps}(w_{\eps, n},B_R(x_n))-C\delta\ge \sum_{n=1}^N \lim_{\varepsilon \to 0}\phi_\varepsilon(R,d_n,x_n)-C\delta
\end{split}
\]
where we also used \eqref{eq:civuolesotto}. Passing to the limit for $\delta\to 0$ we obtain
\[
\liminf_{\eps\to 0}\mathscr{E_\eps}(\mu_\eps)\ge \sum_{n=1}^N\lim_{\varepsilon \to 0}\phi_\varepsilon(R, d_n, x_n)\,.
\]
By exploiting Lemma \ref{lm:cell}, it is now easy to infer that 
\[
\liminf_{\eps\to 0}\mathscr{E_\eps}(\mu_\eps)\ge \sum_{n=1}^N\psi(d_n)\,,
\]
which conclude the proof of the $\Gamma\text{-}\liminf$ inequality. \\

\noindent\textit{$\Gamma\text{-}\limsup$ {\rm inequality}.} Let $\mu\in\calM(\Omega)$ be such that $\mathscr{E}(\mu)<\infty$, otherwise there is nothing to prove. This entails
\[
\frac{\mu}{4\pi}=\sum_{n=1}^Nd_n\delta_{x_n}\,,
\]
for some $d_1,\dots,d_N\in\mathZ$ and $\{x_1,\dots, x_N\}\subset\Omega$, and 
\[
\mathscr{E}(\mu)=\sum_{n=1}^N \psi(d_n)\,. 
\]
Let $R>0$ be given by
\[
R=\min \{|x_n-x_m|: n,m\in\{1\dots, N\}, n\neq m\}\,.
\]
For every $n=1,\dots, N$ we select a sequences $(r_\eps)$ and $(w_{\eps,r_\eps,n})\subset\calSF_\eps$ such that $r_\eps\searrow 0$ and $0<r_\eps<R$, and   $w_{\eps,r_\eps,n}$ is an admissible test spin field for $\phi_\eps(r_\eps,d,x_n)$ satisfing
\[
 E_\eps(w_{\eps,r_\eps,n}, B_r(x_n))\leq \phi_\eps(r_\eps,d_n,x_n)+\eps \,.
\]
This provides also the existence of a sequence of Radon measures $(\zeta_{\eps,r,n})\subseteq\calM_{adm}(\Omega, w_{\eps,r_\eps,n})$, such that
\[
\zeta_{\eps,r_\eps,n}(B_{r_\eps}(x_n))=d_n\,.
\]
Exploiting this latter identity, Theorem \ref{th:discomp} and the fact that $w_{\eps,r_\eps,n}$ is constantly equal to $\boldn$ outside $B_{r_\eps}(x_n)$, it is easy to show that $\zeta_{\eps,r_\eps,n}\wto 4\pi d_i\delta_{x_i}$.
We set $W_\eps\in\calSF_\eps(\Omega)$ by letting 
\[
w_\eps(k)=
\begin{cases}
w_{\eps,r_\eps, n}(k) & \text{ if }k\in B_{r_\eps}(x_n)\,,\\
\boldn & \text{otherwise}.
\end{cases}
\]
We also set $\zeta_\eps=\sum_{n=1}^N\zeta_{\eps,n\mres B_R(x_n)}$. 
Clearly $\zeta_\eps\in\calM_{adm}(\Omega, w_\eps)$. Moreover it holds
\[
 E_\eps(w_{\eps},\Omega)=\sum_{n=1}^N E_\eps(w_{\eps,r_\eps n}\,, B_{r_\eps}(x_n))\,, \text{ and }\zeta_\eps\wto\mu\,.
\] 
Thanks to Lemma \ref{lm:cell}, we get
\[
\begin{split}
\mathscr{E}(\mu)=\sum_{n=1}^N \psi(d_n)\ge \lim_{\eps\to 0}\left(\sum_{n=1}^N \phi_\eps(r_\eps,d_n,x_n)\right) =\lim_{\eps\to 0}\left(\sum_{n=1}^N E_\eps(w_{\eps,r_\eps, n}, B_{r_\eps}(x_n))\right)=\lim_{\eps\to 0} E_\eps(w_{\eps},\Omega)\,.
\end{split}
\]
By the very definition of $\mathscr{E}_\varepsilon$ in  \eqref{df.scrEeps} this implies
\[
\mathscr{E}(\mu)\ge\limsup_{\eps\to 0}\mathscr{E}_\varepsilon(\zeta_\eps)\,,
\]
and concludes the proof. 
\end{proof}

\appendix
\section{}\label{sec:app1}

\begin{lemma}\label{lm:cutoff}
Let $u\in H^1(B_1;\mathR^3)$. 
There exists $U\in   H^1(B_2;\mathR^3)$, such that $U=u$ on $B_1$, $U$ is constant  on $\partial B_2$ and
\[
\begin{split}
\int_{B_2} |\nabla U|^2\,\mathrm{d}x\le C\int_{B_1}|\nabla u|^2\,\mathrm{d}x\,, \quad \text{and} \quad 
\|U\|_{L^{\infty}(B_2)}\le C \|u\|_{L^\infty(\partial B_1)}\,,
\end{split}
\]
where $C$ is a dimensional constant. 
\end{lemma}

\begin{proof}
By arguing component-wise we can assume that $u$ is scalar valued. 
Let $H(u)$ be the $H^1$ harmonic extension of $u_{|\partial B_1}$ into $B_1$, namely
\[
\Delta H(u)=0 \hbox{ in }B_1\,,\quad H(u)=u\hbox{ on }\partial B_1\,.  
\]
Minimality and the maximum principle give
\[
\int_{B_1}|\nabla H(u)|^2\,\mathrm{d}x\le\int_{B_1} |\nabla u|^2\,\mathrm{d}x\,,\quad \|H(u)\|_{L^\infty(B_1)}\le\|u\|_{L^\infty(\partial B_1)}\,.
\]
We extend $H(u)$ to $B_2$ via spherical inversion, that is
\[
\begin{split}
H(u)(x)=(H(u)\circ k)(x) \hbox{ for } x\in B_2\setminus B_1\,,
\hbox{ with } k(x)=\frac{x}{|x|^2}\,.
\end{split}
\]
Then, routine computations show
\begin{equation}\label{eq:A1}
\int_{B_2} |\nabla H(u)|^2\,\mathrm{d}x\le C \int_{B_1}|\nabla u|^2\,\mathrm{d}x\,,\quad \|H(u)\|_{L^\infty(B_2)}\le\|u\|_{L^\infty(\partial B_1)}\,,
\end{equation}
for some dimensional constant $C>0$. 
By the coarea formula we select $\rho\in (\frac 4 3,\frac 5 3)$ such that
\begin{equation}\label{eq:truccocutoff}
\int_{\partial B_\rho}|\nabla H(u)|^2\,\mathrm{d}\mathcal{H}^1\le 3 \int_{A_{\frac 4 3,\frac 5 3}}|\nabla H(u)|^2\,\mathrm{d}x\,,
\end{equation} 
and we let
\[
u_0=\frac 1 {2\pi \rho}\int_{\partial B_\rho} H(u)\,\mathrm{d}\mathcal{H}^1\,.
\]
Notice that \eqref{eq:A1} gives also 
\begin{equation}\label{eq:A11}
|u_0|\le\|u\|_{L^\infty(\partial B_1)}\,.
\end{equation}
We define $\phi_\rho:A_{\rho,2}\to [0,1]$ and $\psi_\rho:\mathR^2\setminus\{0\}\to \partial B_\rho$ to be 
\[
\phi_\rho(x)=\left(\frac{2-|x|}{2-\rho}\right)\,,\quad  \psi_\rho(x)=\rho\frac{x}{|x|}\,,
\]
respectively and we set
\[
v(x)=\phi_\rho(x)H(u)(\psi_\rho(x))+\left(1-\phi_\rho(x)\right)u_0\,,\quad \hbox{for }x\in B_2\setminus B_\rho\,.
\] 
Thanks to \eqref{eq:A1} and \eqref{eq:A11} we have 
\begin{equation}\label{eq:A21}
\|v\|_{L^{\infty}( A_{2,\rho})}\le C\|u\|_{L^\infty(\partial B_1)}\,.
\end{equation}
Moreover, it holds
\begin{equation}\label{eq:A22}
|\nabla v(x)|^2\le C\left|\nabla H(u)\right|^2(\psi_\rho(x))+C\left|H(u)\left(\psi_\rho(x)\right)-u_0\right|^2\,,
\end{equation}
and the following estimate
\begin{equation}\label{eq:truccocutoff2}
(\mathrm{osc}_{\partial B_\rho} H(u))^2\le C \int_{\partial B_\rho}|\nabla H(u)|^2\, \mathrm{d}\mathcal{H}^1\,.
\end{equation}
Therefore, integrating over $A_{\rho,2}$ the inequality \eqref{eq:A22}, taking into account \eqref{eq:truccocutoff}, and \eqref{eq:truccocutoff2}, we obtain
\begin{equation}\label{eq:A2}
\int_{A_{\rho,2}}|\nabla v|^2\,\mathrm{d}x\le C \int_{B_1}|\nabla u|^2\,\mathrm{d}x\,.
\end{equation} 
At last, the lemma is proved by setting  $ U:B_2\to \mathR$ to be
\[
U(x)=
\begin{cases}
u(x) & \hbox{ if }x\in B_1\,,\\
H(u(x)) &\hbox{ if }x\in A_{1,\rho}\,,\\
v(x) &\hbox{ if }x\in A_{\rho,2}\,.
\end{cases}
\]
and by combining \eqref{eq:A1}, \eqref{eq:A21},  \eqref{eq:A2}. 
\end{proof}

\begin{lemma}\label{lm:relper}
Let $\Omega\subset\mathR^2$ be a Lipschitz open set. There exists a constant $C=C_\Omega$, such that for every $E\subset\Omega$ connected open set of finite perimeter, it holds
\[
\mathrm{diam}(E)\le C (|E|^{1/2}+P(E;\Omega))\,,
\] 
where $C=C(\Omega)>0$.
\end{lemma}

\begin{proof}
We argue by contradiction, assuming that the statement it false, so that for every $n\in\mathN$, there exists a connected open set with finite perimeter $E_n\subset \Omega$, satisfying
\begin{equation}\label{A:contraddiction}
(|E_n|^{1/2}+P(E_n;\Omega))\le\frac{\diam(E_n)}{n}\le\frac{\diam(\Omega)}{n}\,.
\end{equation}

\textit{STEP 1.} Let $D>0$ be a constant to be fixed. We prove that for $D$  large enough, for every $x\in E_n$ it holds $E_n\subseteq Q_{\frac{D}{n}}(x)$ for every $n\ge \max\{4D, \diam(\Omega)\}$.  Indeed, assume that $\left(\Omega\setminus Q_{\frac{D}{n}}\right)\cap E_n$ is nonempty. Then  being $E_n$ connected we have that
\[
\partial Q_t(x)\cap E_n\neq\emptyset \quad \text{ for every }t\in \left(0, \frac{D}{n}\right)\,.
\]
In fact, being $\partial Q_t(x)$ connected and having $\Omega$ Lipschitz boundary, one easily infers that there exists a constant $C_\Omega>0$ such that  for every $t\in \left(0, \frac{D}{n}\right)$ satisfying $\mathcal{H}^1(\partial Q_t(x) \cap E_n) < C_\Omega t $ there holds $\Omega\cap \partial Q_t(x)\cap \partial E_n\neq \emptyset$. 
We set
\[
\begin{split}
& I_1=\left\{t\in \left(\frac{D}{2n}, \frac{D}{n}\right)\colon \mathcal{H}^1(\partial Q_t(x) \cap E_n) \geq C_\Omega t\right\}\,,\\
& I_2=\left\{t\in \left(\frac{D}{2n},\frac{D}{n}\right)\colon \mathcal{H}^0(\Omega\cap \partial Q_t(x)\cap \partial E_n)\ge 1\right\}\,.
\end{split}
\]
From the above observations we have that $|I_1|+|I_2|\geq\frac{D}{2n}$. 
Let $f:\mathR^2\to\mathR^1$, be given by $f(z)=\max\{|(z-x)\cdot e_1|,|(z-x)\cdot e_2|\}$. It holds $|\nabla f|\le 1$.
The coarea formula gives 
\begin{align*}
|E_n|\ge \int_{Q_{\frac{D}{n}}(x)} |\nabla f| \chi_{E_n}(x)\,\mathrm{d}x &\ge \int_0^{\frac{D}{n}} \int_{\partial Q_{t}(x)}\chi_{E_n}(x)\,\mathrm{d}\mathcal{H}^1(x)\,\mathrm{d}t\\& = \int_0^{\frac{D}{n}}\mathcal{H}^1(\partial Q_t(x) \cap E_n)\,\mathrm{d}t \ge \int_{I_1}C_\Omega t\, \mathrm{d}t \ge C_\Omega |I_1|\frac{D}{2n} \geq C_\Omega |I_1|^2 \,. 
\end{align*} 
Similarly, denoting by $\nu(y) \in \mathbb{S}^2$ the outer unit normal to $E$ at $y\in \partial^* E$, by the slicing-property of sets of finite perimeter, we have
\[
\sqrt{2}\mathcal{H}^1(\partial E_n\cap\Omega)\ge \int_{\partial E_n\cap\Omega\cap Q_{\frac{D}{n}}(x)} |\nu\cdot e_1| +|\nu\cdot e_2| \,\mathrm{d}\mathcal{H}^1\ge \int_0^{\frac{D}{n}}\mathcal{H}^0(\partial E_n\cap\Omega\cap \partial Q_{t}(x)) \,\mathrm{d}t\ge |I_2|\,,
\]
where in the last inequality we used that $n\ge 4D $.
Combining these two inequalities we get 
\[
\begin{split}
& \sqrt{|E_n|}+\mathcal{H}^1(\partial E_n\cap\Omega)\ge \sqrt{C_\Omega} |I_1|+\frac{1}{\sqrt{2}}|I_2| \geq \min\left\{\sqrt{C_\Omega},\frac{1}{\sqrt{2}}\right\}\frac{D}{2n} \,.
\end{split}
\] 
By choosing $D>2 \min\left\{\sqrt{C_\Omega},\frac{1}{\sqrt{2}}\right\}^{-1} \diam(\Omega)$, the latter contradicts \eqref{A:contraddiction}.\\

\noindent\textit{STEP 2.} We show that for every $x\in E_n$ there exists $z\in\partial\Omega$ such that $z\in Q_{\frac{D}{n}}(z)$. This easily follows from the previous point. Indeed, let $x\in E_n$. We have $E_n\subseteq Q_{\frac{D}{n}}(x)$. Assume that $\mathrm{dist}(x, \partial \Omega)>\sqrt{2}\frac{D}{n}$, so that $E_n\subset\joinrel\subset\Omega$. The standard Perimeter-diameter estimate then imply
\[
\tilde C \,\diam(\Omega)\le \sqrt{|E_n|}+P(E_n)=\sqrt{|E_n|}+P(E_n;\Omega)
\]
for some $\tilde C=\tilde C_\Omega$. For $n^{-1}<\tilde C$ this yields a contradiction to \eqref{A:contraddiction}. \\

\noindent\textit{STEP 3.} We now show, that \eqref{A:contraddiction} leads to a contraddiction.
Clearly, it is enough to prove that, for some $n_1>0$ and for every $n>n_1$ it holds: 
\begin{equation}\label{A:provare}
\mathcal{H}^1( E_n\cap \partial\Omega)\le C\, P(E_n;\Omega)
\end{equation}
for some constant $C=C_\Omega$, since this implies $P(E_n)\le P(E_n;\Omega)$ and, again by the standard Perimeter-diameter, estimate contradicts \eqref{A:contraddiction}. We prove \eqref{A:provare} again with a slicing argument. 
We select $r=r_\Omega$ in such a way that, for every $z\in\partial\Omega$, we have that $\partial\Omega\cap Q_r(z)$ is the graph of some Lipschitz function. 
Let $n_1>0$ be such that $\frac{D}{n_1}<\frac{r}{4}$. Let $n>n_1$. Pick  $x\in E_n$, and let $z_x\in Q_{\frac{D}{n}}(x)\cap\partial\Omega$ obtained as in the previous point. Up to translation we can suppose $z_x=0$. 
Let $\varphi$ be a positive Lipschitz function $\varphi: (-r/2,r/2)\to (0,\infty)$, such that $\partial\Omega\cap Q_r(0)$ is the graph of $\varphi$. Furthermore we can suppose that $e_2$ is the normal vector to $\partial\Omega$ in $0$.  Thus we have
\[
\Omega\cap Q_r=\{(t,s): t\in (-l,l),\, s>\varphi(t)\}\,.
\] 
We set 
\[
\omega_n=\{t\in(-l,l): (t,\varphi(t))\in\partial E_n\}\,,
\]
and note that the area formula and the fact that $\varphi$ is Lipschitz-continuous imply
\[\mathcal{H}^1(\partial E_n\cap \partial \Omega)\le C\mathcal{L}^1(\omega_n)\,.\]
For every $x=(t,\varphi(t))\in\omega_n$ we have that $E_n\subset Q_{\frac{D}{n}}(x)$ thanks to the first step of the proof. This entails that there exists a point $(t,s_t)$, such that the vertical line passing trough $(t,0)$ and $(t, \varphi(t))$ meets $\partial E_n\cap \Omega$. Furthermore by construction it holds $(t,s_t)\in Q_r$, and hence by the graphicality of $\partial\Omega$ inside $Q_r$ it cannot be $(t,s_t)\in\partial\Omega$. Therefore, the slicing-properties of sets of finite perimeter imply 
\[
\mathcal{L}^1(\omega_n)\le \int_{\omega_n} \mathcal{H}^0(\partial E_n \cap \Omega \cap \{(t,s) \colon s\in \mathbb{R}\})\,\mathrm{d}t\leq  \int_{\partial E_n \cap \Omega} |\nu_2|\,\mathrm{d}\mathcal{H}^1 \leq  \mathcal{H}^1(\partial E_n\cap \Omega)\,.
\]
We therefore have
\[
\mathcal{H}^1(\partial E_n\cap \partial \Omega)\le C\mathcal{H}^1(\partial E_n\cap \Omega)\,,
\]
proving \eqref{A:provare}. This concludes the proof. 
\end{proof}

\section{}\label{sec:app2}
In this section we prove two elementary facts that we used in the proof of Theorem \ref{theo:interpolation}.
\begin{lemma}\label{lm:costruzioneu_R}
Let $R\in\calR_\eps(u)$ and $S_R(u)\subseteq \mathS^2$ be an admissible interpolation surface for $u$ in $R$. Then, there exists $u_R\in\mathS^2$ such that
\begin{equation*}
u_R\in S_R(u),\quad u_R\cdot z\ge 0,\hbox{ for every }z\in S_R(u)\,.
\end{equation*}
\end{lemma}
\begin{proof}
If $\mathrm{cone}_R(u)$ is a linear subspace of $\mathR^3$, then $S_R(u)$ is given as in \eqref{df:SR2} and it is enough to set $u_R=h$.
If instead $\mathrm{cone}_R(u)$ is not a linear subspace of $\mathR^3$, we can consider the so-called \textit{polar cone}:
\[
\begin{split}
\mathrm{cone}^*_R(u)=\{z\in\mathR^3: z\cdot w\ge 0\text{ for all }w\in\mathrm{cone}_R(u)\}\,.
\end{split}
\]
Then, it is a known fact, see for instance \cite[Theorem 2.1]{Ga}, that there exists
\[
x\in \left(\mathrm{cone}^*_R(u)\cap\mathrm{cone}_R(u)\right)\setminus\{0\}\,.
\]
By letting $u_R=\frac{x}{|x|}$ and recalling \eqref{df:SR1} we conclude.  
\end{proof}

\begin{lemma} \label{lm:uHat} Let $u \in \mathcal{SF}_\varepsilon(\Omega)$. Let $R \in \mathcal{R}_\varepsilon(u)$ and let $\overarc u \colon R \to \mathbb{S}^2$ be defined as in \eqref{eq:intext}. Then, the following statements hold true:
\begin{itemize}
\item[(i)] $\overarc u(R) = S_R(u)$;
\item[(ii)]  for $\mathcal{H}^2$-a.e. $y \in S_R(u)$ there holds
\begin{align*}
\#\{x \in R \colon \overarc u(x) = y\}=1\,.
\end{align*}
\end{itemize}
\end{lemma} 
\begin{proof}
Let $N=\#u\left(\mathcal{L}_\varepsilon(R)\right)$. By Lemma~\ref{lm:mistero3} we have that $N \in \{1,2,3,4\}$.

In the cases where either $N=1$ or $N=2$ or $\cone_R(u)$ is not a plane but it lies on a plane, the conclusion is straightforward. In fact, in all these cases (ii) is trivial being $\mathcal{H}^2(S_R(u))=0$. Thus suppose that we are not dealing with one of the aforementioned cases. In particular it has to be $N=3$ and $N=4$. We describe only the latter case, the other being similar. 
To fix the notation we have
\[
u(\mathcal{L}_\varepsilon(R))=\{p_1,p_2,p_3,p_4\}\] 
with $p_1=-p_2$. By Lemma~\ref{lm:mistero1}, there exist two triangles $T_1=[i_1,i_2,i_3],T_2=[j_1,j_2,j_3] \in \mathcal{T}_\varepsilon(R)$ such that $u(\mathcal{L}_\varepsilon(T_1)) = \{p_1,-p_1,p_3\}$,  $u(\mathcal{L}_\varepsilon(T_1)) = \{p_1,-p_1,p_4\}$ and all the other triangles $T \in \mathcal{T}_\varepsilon(R)$ such that $T \notin \{T_1,T_2\}$ satisfy $u(\mathcal{L}_\varepsilon(T))=\{p_1,-p_1\}$. For all $T \notin \{T_1,T_2\}$ we can easily argue  that $\overarc u(T) =  \gamma_{p_1,u_R}\cup \gamma_{-p_1,u_R}$. For such $T$ we have $\mathcal{H}^2(\overarc u(T))=\mathcal{H}^2(\gamma_{p_1,u_R}\cup \gamma_{-p_1,u_R})=0$. In order to treat the triangles $T_1,T_2$ we distinguish the two possible cases:
\begin{itemize}
\item[]\emph{Case 1:} $p_1,p_3,p_4$ are linear independent;
\item[]\emph{Case 2:} $p_1,p_3,p_4$ are linear dependent, $\mathrm{cone}_R(u)$ is a plane.
\end{itemize}

\noindent\textit{Case 1}. In this case, by  Lemma \ref{rem:geometrico}, we have $S_R(u) = S(p_1,p_3,p_4) \cup S(-p_1,p_3,p_4)$. Notice that $u_R \in S_R(u)$ satisfies $u_R \cdot p_1=0$, $u_R\cdot p_i \geq 0$ for $i=3,4$. We assume without loss of generality that the triples $(p_1,u_R,p_3)$ and  $(p_1,u_R,p_4)$ are linear independent. We show ({\rm i}) and ({\rm ii}) for $u\vert_{T_1\cup T_2}$. Once this is done, this concludes the proof. Using the definition of $\overarc u$ we have that
\begin{align*}
\begin{split}
\overarc u(T_1) = S\left(p_1,\frac{p_1+p_3}{|p_1+p_3|},u_R\right) &\cup S\left(p_3,\frac{p_1+p_3}{|p_1+p_3|},u_R\right) \\&\cup  S\left(-p_1,\frac{p_3-p_1}{|p_3-p_1|},u_R\right) \cup S\left(p_3,\frac{p_3-p_1}{|p_3-p_1|},u_R\right)\,.
\end{split}
\end{align*}
Lemma~\ref{lm:triangle}, applied with $\lambda=\frac{1}{2}$ and the triples $(p_1,p_3,u_R)$ and $(-p_1,p_3,u_R)$ respectively, implies
\begin{align}\label{eq:imageT1}
\overarc u(T_1) = S(p_1,p_3,u_R) \cup  S(-p_1,p_3,u_R)\,.
\end{align}
Similarly, we obtain
\begin{align}\label{eq:imageT2}
\overarc u(T_2) = S(p_1,p_4,u_R) \cup  S(-p_1,p_4,u_R)\,.
\end{align}
Note that for every $T'\in \mathcal{T}_\varepsilon'(T_1\cup T_2)$ the function $\overarc u\vert_{T'}$ is either bijective or $\overarc u(T') \subseteq \gamma_{p_1,u_R}\cup \gamma_{-p_1,u_R} \cup \gamma_{p_3,p_4}$. Furthermore, it is easy to verify that $\mathcal{H}^2(\overarc{u}(T')\cap \overarc{u}(T''))=0$ for all $T', T''\in \mathcal{T}_\varepsilon'(T_1\cup T_2)$. This shows ({\rm ii}) once ({\rm i}) is proven. In order to do so, using \eqref{eq:imageT1} and \eqref{eq:imageT2}, it suffices to show
\begin{align*}
\begin{split}
 S(p_1,p_3,u_R) &\cup  S(-p_1,p_3,u_R)\\&\cup S(p_1,p_4,u_R) \cup  S(-p_1,p_4,u_R) =  S(p_1,p_3,p_4) \cup  S(-p_1,p_3,p_4)\,.\end{split}
\end{align*}
Assume without loss of generality that $u_R \in S(p_1,p_3,p_4)$, i.e., $u_R = \lambda_1 p_1 + \lambda_2 p_3 + \lambda_3 p_4$ for $\lambda_i \geq 0$. Now, let  $y =\mu_1 p_1 + \mu_2 p_3 +\mu_3 p_4$ with $\mu_1 \in \mathbb{R}$ and $\mu_2,\mu_3 \geq 0$. We have
\begin{align*}
y = \begin{cases} 
-\left(\frac{\mu_2}{\lambda_2}\lambda_1 -\mu_1 \right) p_1+\frac{\mu_2}{\lambda_2}u_R + \left(\mu_3 - \frac{\mu_2}{\lambda_2}\lambda_3\right)p_4  &\text{if } \frac{\mu_1}{\lambda_1} \leq\frac{\mu_2}{\lambda_2}\leq \frac{\mu_3}{\lambda_3}   \,,\\
-\left(\frac{\mu_3}{\lambda_3}\lambda_1 -\mu_1 \right) p_1+\frac{\mu_3}{\lambda_3}u_R + \left(\mu_2 - \frac{\mu_3}{\lambda_3}\lambda_2\right)p_3  &\text{if } \frac{\mu_1}{\lambda_1} \leq\frac{\mu_3}{\lambda_3}\leq \frac{\mu_2}{\lambda_2}\,,\\
\left(\mu_1 - \frac{\mu_3}{\lambda_3}\lambda_1\right) p_1 + \left(\mu_2 - \frac{\mu_3}{\lambda_3}\lambda_2\right)p_3 +\frac{\mu_3}{\lambda_3}u_R &\text{if } \frac{\mu_3}{\lambda_3} \leq \min\left\{\frac{\mu_i}{\lambda_i} \colon i=1,2,3\right\}\,,\\
\left(\mu_1 - \frac{\mu_2}{\lambda_2}\lambda_1\right) p_1+\frac{\mu_2}{\lambda_2}u_R + \left(\mu_2 - \frac{\mu_2}{\lambda_2}\lambda_3\right)p_4  &\text{if } \frac{\mu_2}{\lambda_2} \leq \min\left\{\frac{\mu_i}{\lambda_i} \colon i=1,2,3\right\}\,,
\end{cases}
\end{align*}
which is to say
\begin{align*}
y\in \begin{cases} S(-p_1,p_4,u_R) &\text{if }  \frac{\mu_1}{\lambda_1} \leq\frac{\mu_2}{\lambda_2}\leq \frac{\mu_3}{\lambda_3} \,,\\
S(-p_1,p_3,u_R) &\text{if }\frac{\mu_1}{\lambda_1} \leq\frac{\mu_3}{\lambda_3}\leq \frac{\mu_2}{\lambda_2}\,, \\
S(p_1,p_3,u_R)  &\text{if } \frac{\mu_3}{\lambda_3} \leq \min\left\{\frac{\mu_i}{\lambda_i} \colon i=1,2,3\right\}\,,\\
S(p_1,p_4,u_R)&\text{if } \frac{\mu_2}{\lambda_2} \leq \min\left\{\frac{\mu_i}{\lambda_i} \colon i=1,2,3\right\}\,.
\end{cases}
\end{align*}
This shows 
\begin{align*}
  S(p_1,p_3,p_4) \cup  S(-p_1,p_3,p_4)\subseteq S(p_1,p_3,u_R) \cup  S(-p_1,p_3,u_R)\cup S(p_1,p_4,u_R) \cup  S(-p_1,p_4,u_R)  \,.
\end{align*}
In order to show the reverse inclusion, observe first, as $u_R \in S(p_1,p_3,p_4)$, that $S(p_1,p_4,u_R) \cup S(p_1,p_3,u_R) \subset S(p_1,p_3,p_4)$. It remains to show that 
\begin{align*}
   S(-p_1,p_3,u_R)\cup  S(-p_1,p_4,u_R) \subset  S(p_1,p_3,p_4) \cup  S(-p_1,p_3,p_4)  \,.
\end{align*}
Let $y \in  S(-p_1,p_3,u_R)$, i.e., $y=-\mu_1 p_1 +\mu_2 p_3 + \mu_3 u_R$ with $\mu_1 \in \mathbb{R}$ and $\mu_2,\mu_3 \geq 0$. Then $y = (\mu_3\lambda_1+\mu_1)p_1 + (\mu_2+\mu_3\lambda_2)p_3 + \mu_3\lambda_3p_4$ which is to say that 
\begin{align*}
y \in \begin{cases} S(-p_1,p_3,p_4) &\text{if } \mu_3\lambda_1+\mu_1\leq 0\,,\\
 S(p_1,p_3,p_4) &\text{if } \mu_3\lambda_1+\mu_1\geq 0\,.
\end{cases}
\end{align*}
proving the desired inclusion. 

\noindent\textit{Case 2}. In this case, there exists $h \in \mathbb{S}^2$ such that $S_R(u) = \mathbb{S}^2 \cap \{x \in \mathbb{R}^3 \colon x\cdot h \geq 0\}$ and $u_R =h$. Following the same argument as in the previous step we have that
\begin{align*}
 S\left(p_1,p_4, h\right) \cup S\left(-p_1,p_4, h\right) \cup S\left(p_1,p_3, h\right) \cup S\left(-p_1,p_3, h\right) \subseteq S_R(u)\,.
\end{align*}
and that, in order to conclude, it is enough to show the opposite set inclusion. This comes as a conseuqence of Caratheodory's Theorem: indeed given $y \in S_R(u)$, $y = y_\Pi + \lambda u_R$ with $y_\Pi \in \{x\in \mathbb{R}^3\colon h\cdot x =0\}$ and $\lambda \geq 0$, we have that $y_\Pi = \lambda_1 q_1 + \lambda_2 q_2$ with $q_1, q_2 \in \{-p_1, p_1,p_3,p_4\}$ and $\lambda_i \geq 0$.
\end{proof}

\begin{lemma}\label{lm:triangle} Let $p_1,p_2,p_3 \in \mathbb{S}^2$ be linear independent and such that $p_i \neq -p_j$ for all $i\neq j$. For all $\lambda \in [0,1]$ there holds 
\begin{align*}
 S\left(p_1,\frac{\lambda p_1+(1-\lambda)p_2}{|\lambda p_1+(1-\lambda)p_2|}, p_3\right) \cup  S\left(p_2,\frac{\lambda p_1+(1-\lambda)p_2}{|\lambda p_1+(1-\lambda)p_2|}, p_3\right)=  S\left(p_1,p_2, p_3\right)\,.
\end{align*}
\end{lemma}
\begin{proof} We can assume that $\lambda \in (0,1)$ as otherwise there is nothing to prove. For the remainder of this proof we set $p_\lambda = \frac{\lambda p_1+(1-\lambda)p_2}{|\lambda p_1+(1-\lambda)p_2|} $.  As $p_\lambda \in S\left(p_1,p_2, p_3\right)$ we have 
\begin{align*}
 S\left(p_1,p_\lambda, p_3\right) \cup  S\left(p_2,p_\lambda, p_3\right)\subseteq  S\left(p_1,p_2, p_3\right)\,.
\end{align*}
It remains to show the opposite inclusion. To this end let $y \in S(p_1,p_2,p_3)$, i.e. there exist $\lambda_i \geq 0$ such that $w_y =\sum_{i}\lambda_i p_i$ and $\frac{w_y}{|w_y|} = y$.
 We can write
\begin{align*}
w_y = \begin{cases} \mu_1 p_1 + \mu_2p_\lambda +\lambda_3 p_3 &\text{with } \mu_1= \lambda_1-\frac{\lambda}{1-\lambda}\lambda_2\,,\, \mu_2= \frac{|\lambda p_1+(1-\lambda)p_2|}{1-\lambda} \lambda_2\,, \text{ if } \lambda_1 \geq \frac{\lambda}{1-\lambda}\lambda_2\,, \\
\mu_1 p_1 + \mu_2p_\lambda +\lambda_3 p_3 &\text{with } \mu_1=\frac{|\lambda p_1+(1-\lambda)p_2|}{\lambda}\lambda_1\,,\,\mu_2=\lambda_2-\frac{1-\lambda}{\lambda}\lambda_1\,,  \text{ if } \lambda_1 \leq \frac{\lambda}{1-\lambda}\lambda_2\,.
\end{cases}
\end{align*} 
In other words this shows that
\begin{align*}
y = \frac{w_y}{|w_y|} \in \begin{cases}  S\left(p_1,p_\lambda, p_3\right) &\text{if } \lambda_1 \geq \frac{\lambda}{1-\lambda}\lambda_2\,, \\
S\left(p_2,p_\lambda, p_3\right) &\text{if } \lambda_1 \leq \frac{\lambda}{1-\lambda}\lambda_2\,.
\end{cases}
\end{align*}
This concludes the proof.
\end{proof}

\vskip10pt
\textsc{Acknowledgments:}  The research of L.~Kreutz and L.~Briani was supported by the DFG through the Emmy Noether Programme (project number 509436910).

\end{document}

%% file: coppole.pdf_tex
\begingroup%
  \makeatletter%
  \providecommand\color[2][]{%
    \errmessage{(Inkscape) Color is used for the text in Inkscape, but the package 'color.sty' is not loaded}%
    \renewcommand\color[2][]{}%
  }%
  \providecommand\transparent[1]{%
    \errmessage{(Inkscape) Transparency is used (non-zero) for the text in Inkscape, but the package 'transparent.sty' is not loaded}%
    \renewcommand\transparent[1]{}%
  }%
  \providecommand\rotatebox[2]{#2}%
  \newcommand*\fsize{\dimexpr\f@size pt\relax}%
  \newcommand*\lineheight[1]{\fontsize{\fsize}{#1\fsize}\selectfont}%
  \ifx\svgwidth\undefined%
    \setlength{\unitlength}{469.69498281bp}%
    \ifx\svgscale\undefined%
      \relax%
    \else%
      \setlength{\unitlength}{\unitlength * \real{\svgscale}}%
    \fi%
  \else%
    \setlength{\unitlength}{\svgwidth}%
  \fi%
  \global\let\svgwidth\undefined%
  \global\let\svgscale\undefined%
  \makeatother%
  \begin{picture}(1,0.34338349)%
    \lineheight{1}%
    \setlength\tabcolsep{0pt}%
    \put(0,0){\includegraphics[width=\unitlength,page=1]{coppole.pdf}}%
    \put(0.86877425,0.20721035){\color[rgb]{0,0,0}\makebox(0,0)[lt]{\lineheight{1.25}\smash{\begin{tabular}[t]{l}$p_1$\end{tabular}}}}%
    \put(0.81581058,0.14561482){\color[rgb]{0,0,0}\makebox(0,0)[lt]{\lineheight{1.25}\smash{\begin{tabular}[t]{l}$p_2$\end{tabular}}}}%
    \put(0.74605643,0.20191096){\color[rgb]{0,0,0}\makebox(0,0)[lt]{\lineheight{1.25}\smash{\begin{tabular}[t]{l}$p_3$\end{tabular}}}}%
    \put(0.95369842,0.20149436){\color[rgb]{0,0,0}\makebox(0,0)[lt]{\lineheight{1.25}\smash{\begin{tabular}[t]{l}$p_4$\end{tabular}}}}%
    \put(0,0){\includegraphics[width=\unitlength,page=2]{coppole.pdf}}%
    \put(0.62674849,0.20832121){\color[rgb]{0,0,0}\makebox(0,0)[lt]{\lineheight{1.25}\smash{\begin{tabular}[t]{l}$p_1$\end{tabular}}}}%
    \put(0.57847569,0.14780827){\color[rgb]{0,0,0}\makebox(0,0)[lt]{\lineheight{1.25}\smash{\begin{tabular}[t]{l}$p_2$\end{tabular}}}}%
    \put(0.52567991,0.15718319){\color[rgb]{0,0,0}\makebox(0,0)[lt]{\lineheight{1.25}\smash{\begin{tabular}[t]{l}$p_3$\end{tabular}}}}%
    \put(0.70806572,0.17589565){\color[rgb]{0,0,0}\makebox(0,0)[lt]{\lineheight{1.25}\smash{\begin{tabular}[t]{l}$p_4$\end{tabular}}}}%
    \put(0,0){\includegraphics[width=\unitlength,page=3]{coppole.pdf}}%
    \put(0.36325003,0.13707136){\color[rgb]{0,0,0}\makebox(0,0)[lt]{\lineheight{1.25}\smash{\begin{tabular}[t]{l}$p_1$\end{tabular}}}}%
    \put(0.44001295,0.14791297){\color[rgb]{0,0,0}\makebox(0,0)[lt]{\lineheight{1.25}\smash{\begin{tabular}[t]{l}$p_2$\end{tabular}}}}%
    \put(0,0){\includegraphics[width=\unitlength,page=4]{coppole.pdf}}%
    \put(0.11437559,0.13598196){\color[rgb]{0,0,0}\makebox(0,0)[lt]{\lineheight{1.25}\smash{\begin{tabular}[t]{l}$p_1$\end{tabular}}}}%
    \put(0.16393364,0.22688697){\color[rgb]{0,0,0}\makebox(0,0)[lt]{\lineheight{1.25}\smash{\begin{tabular}[t]{l}$p_2$\end{tabular}}}}%
    \put(0,0){\includegraphics[width=\unitlength,page=5]{coppole.pdf}}%
    \put(0.40918632,0.22933799){\color[rgb]{0,0,0}\makebox(0,0)[lt]{\lineheight{1.25}\smash{\begin{tabular}[t]{l}$p_3$\end{tabular}}}}%
    \put(0,0){\includegraphics[width=\unitlength,page=6]{coppole.pdf}}%
  \end{picture}%
\endgroup%

%% file: Skyrme-final.bbl
\vskip10pt

\begin{thebibliography}{100}


\bibitem{ABCS23} {\sc R.~Alicandro, A.~Braides, M.~Cicalese, and M.~Solci}. \textit{Discrete Variational Problems with Interfaces.}
Cambridge University Press, 2023.

\bibitem{AC09} {\sc R.~Alicandro and M.~Cicalese}. Variational analysis of the asymptotics of the XY model. \textit{Arch. Ration. Mech. Anal.} {\bf 192} (2009), no. 3, 501--536.

\bibitem{ADGP14} {\sc R.~Alicandro, L.~De Luca, A.~Garroni, and M.~Ponsiglione}. Metastability and dynamics of discrete topological singularities in two dimensions: a $\Gamma$-convergence approach. \textit{Arch. Ration. Mech. Anal.} {\bf 214} (2014), no. 1, 269--330.

\bibitem{BCKO22} {\sc A.~Bach, M.~Cicalese, L.~Kreutz and G.~Orlando}. The antiferromagnetic XY model on the triangular lattice: topological singularities. \textit{Indiana Univ. Math. J.} {\bf 71} (2022), no. 6, 2411--2475.

\bibitem{BCGO24} \textsc{A.~Bach, M.~Cicalese, A.~Garroni, and G.~Orlando}. Stacking faults in the limit of a discrete model for partial edge dislocations. \textit {arXiv:2407.03975v1}.

\bibitem{BCDP18} {\sc R.~Badal, M.~Cicalese, L.~De Luca, and M.~Ponsiglione}. $\Gamma$-convergence analysis of a generalized $XY$ model: fractional vortices and string defects. \textit{Commun. Math. Phys.} {\bf 358} (2018), no. 2, 705--739. 


\bibitem{berg1981} {\sc B.~Berg and M.~L{\"u}scher}. Definition and statistical distributions of a topological number in the lattice O (3) $\sigma$-model. \textit{Nuclear Physics B} {\bf 190} (1981),  no. 2, 412--424.

\bibitem{Bernand2021} {\sc A.~Bernand-Mantel, C.B.~Muratov, and T.M.~Simon}. A quantitative description of skyrmions in ultrathin ferromagnetic films and stability of degree \(\pm1\) harmonic maps from \(\mathbb{R}^2\) to \(\mathbb{S}^2\). \textit{Arch. Rat. Mech. Anal.} {\bf 239} (2021), no. 1, 219--299.

\bibitem{Bernand2022} {\sc A.~Bernand-Mantel, C.B.~Muratov, and V.V.~Slastikov}. A micromagnetic theory of skyrmion lifetime in ultrathin ferromagnetic films. \textit{Proc. Natl. Acad. Sci. USA} {\bf 119} (2022), no. 29, e2122237119.

\bibitem{Bethuel1994} {\sc F.~Bethuel, H.~Brezis, and F.~Hélein}. \textit{Ginzburg-Landau vortices}. Progress in Nonlinear Differential Equations and their Applications, 13. Birkhäuser Boston Inc., Boston, MA, 1994.

\bibitem{Bra} {\sc A.~Braides}. \textit{$\Gamma$-convergence for Beginners}. Oxford Lecture Series in Mathematics and its Applications, 22. Oxford University Press, Oxford, 2002.
   
\bibitem{BraCicSol} {\sc A.~Braides, M.~Cicalese and F.~Solombrino}. $Q$-tensor continuum energies as limits of head-to-tail symmetric spins systems. \textit{SIAM J. Math. Anal.} {\bf 47} (2015), no. 4, 2832--2867.

\bibitem{BrCo} {\sc H.~Brezis, and J.-M.~Coron}. Large solutions for harmonic maps in two dimensions. \textit{Comm. Math. Phys.} {\bf 92} (1983), no. 2, 203--215.

\bibitem{Bogdanov1994} {\sc A.~Bogdanov and A.~Hubert}. The properties of isolated magnetic vortices. \textit{Physica Status Solidi (B)} {\bf 186} (1994), no. 2, 527--543.

\bibitem{Bogdanov1999} {\sc A.~Bogdanov and A.~Hubert}. The stability of vortex-like structures in uniaxial ferromagnets. \textit{J. Magn. Magn. Mater.} {\bf 195} (1999), no. 1, 182--192.

\bibitem{Bogdanov1989} {\sc A.N.~Bogdanov, M.V.~Kudinov, and D.A.~Yablonskii}. Theory of magnetic vortices in easy-axis ferromagnets. \textit{Sov. Phys. --Solid State} {\bf  31} (1989), no. 10, 1707--1710.

\bibitem{BogdanovJETP1989} {\sc A.N.~Bogdanov and D.A.~Yablonskii}. Thermodynamically stable “vortices” in magnetically ordered crystals. The mixed state of magnets. \textit{Sov. Phys. – JETP} {\bf 68} (1989), no. 1, 101--103.


\bibitem{CanSeg} {\sc G.~Canevari and A.~Segatti}. Defects  in nematic shells: a $\Gamma$-convergence discrete-to-continuum approach. \textit{Arch. Ration. Mech. Anal.} {\bf 229} (2018), no. 1, 125--186.

\bibitem{COR22CPAM} {\sc M.~Cicalese, G.~Orlando, and M.~Ruf}.	Emergence of concentration effects in the variational analysis of the $N$-clock model. \textit{Comm. Pure Appl. Math.} {\bf 75} (2022), no. 10, 2279--2342.
	%
\bibitem{COR22ARMA} {\sc M.~Cicalese, G.~Orlando, and M.~Ruf}.	The N-Clock Model: Variational analysis of fast and slow convergence rates of N. \textit{Arch. Ration. Mech. Anal.} {\bf 245} (2022), no. 2, 1135--1196.
	%
\bibitem{DM} {\sc G.~Dal Maso}. \textit{ An introduction to $\Gamma$-convergence}. Progress in Nonlinear Differential Equations and their Applications, 8. Birkh\"auser Boston, Inc., Boston, MA, 1993. 
	
\bibitem{DiFratta2024} {\sc G.~Di Fratta, C.B.~Muratov, and V.V.~Slastikov}. Reduced energies for thin ferromagnetic films with perpendicular anisotropy. \textit{Math. Models Methods Appl. Sci.} {\bf 34} (2024), no. 10, 1861--1904.

\bibitem{Doring2017} {\sc L.~Döring and C.~Melcher}. Compactness results for static and dynamic chiral skyrmions near the conformal limit. \textit{Calc. Var. Partial Differential Equations} {\bf 56} (2017), no. 3, 30.

\bibitem{Do} {\sc M.P.~do Carmo}. \textit{ Differential geometry of curves and surfaces.} Dover Publications, Inc., Mineola, New York, 2016.

\bibitem{Foster2019} {\sc D.~Foster, C.~Kind, P.J.~Ackerman, J.-S.~Tai, B. M.R.~Dennis, and I.I.~Smalyukh}. Two-dimensional skyrmion bags in liquid crystals and ferromagnets. \textit{Nature Phys.} {\bf 15} (2019), no. 7, 655--659.

\bibitem{Ga} {\sc J.W.~Gaddum}. A Theorem on Convex Cones with Applications to Linear Inequalities. \textit{Proceedings of the American Mathematical Society} {\bf 3} (1952), no. 6, 957--960.

\bibitem{Gustafson2021} {\sc  S.~Gustafson, and L.~Wang}. Co-rotational chiral magnetic skyrmions near harmonic maps. \textit{J. Funct. Anal.} {\bf 280} (2021), no. 4, 108867. 

\bibitem{Heinze2011} {\sc S.~Heinze, K.~von Bergmann, M.~Menzel, J.~Brede, A.~Kubetzka, R.~Wiesendanger, G.~Bihlmayer, and S.~Blügel}. Spontaneous atomic-scale magnetic skyrmion lattice in two dimensions. \textit{Nature Phys.} {\bf 7} (2011), no. 9, 713--718.


\bibitem{Kiselev2011} {\sc N.S.~Kiselev, A.N.~Bogdanov, R.~Sch\"afer, and U.K.~R\"ossler}. Chiral skyrmions in thin magnetic films: new objects for magnetic storage technologies? \textit{J. Phys. D: Appl. Phys.} {\bf 44} (2011), no. 39, 392001.

\bibitem{Komineas2020} {\sc S.~Komineas, C.~Melcher, and S.~Venakides}. The profile of chiral skyrmions of small radius. \textit{Nonlinearity} {\bf 33} (2020), no. 7, 3395--3408.

\bibitem{Komineas2021} {\sc S.~Komineas, C.~Melcher, and S.~Venakides}. Chiral skyrmions of large radius. \textit{Physica D} {\bf 418} (2021), 132842.

\bibitem{Kuchkin2020} {\sc V.M.~Kuchkin, B.~Barton-Singer, F.N.~Rybakov, S.~Bl\"ugel, B.J.~Schroers, and N.S.~Kiselev}. Magnetic skyrmions, chiral kinks, and holomorphic functions. \textit{Phys. Rev. B} {\bf 102} (2020), 144422. 

\bibitem{Kuchkin2023} {\sc V.M.~Kuchkin, N.S.~Kiselev, F.N.~Rybakov, and P.F.~Bessarab}. Tailed skyrmions--an obscure branch of magnetic solitons. \textit{Front. Phys.} {\bf  11} (2023), 1171079.


\bibitem{Lin2004}  {\sc F.~Lin and Y.~Yang}. Existence of two-dimensional skyrmions via the concentration-compactness method. \textit{Comm. Pure Appl. Math.} {\bf 57} (2004), no. 10, 1332--1351. 


\bibitem{Mel} {\sc C.~Melcher}. Chiral skyrmions in the plane. \textit{Proc. R. Soc. Lond. Ser. A} {\bf 470} (2014), 20140394.

\bibitem{Monteil2023} {\sc A.~Monteil, C.B.~Muratov, T.M.~Simon, and V.V.~Slastikov}. Magnetic skyrmions under confinement. \textit{Commun. Math. Phys.} {\bf 404} (2023), no. 3, 1571--1605.

\bibitem{Muratov2024} {\sc C.B.~Muratov, T.M.~Simon, and V.V.~Slastikov}. Existence of higher degree minimizers in the magnetic skyrmion problem. \textit{arXiv:2409.07205v1}

\bibitem{Muhlbauer2009} {\sc S.~Mühlbauer et al.} Skyrmion lattice in a chiral magnet. \textit{Science} {\bf 323} (2009), no. 5916,  915--919.


\bibitem{P07} {\sc M.~Ponsiglione}. Elastic energy stored in a crystal induced by screw dislocations: from discrete to continuous. \textit{SIAM J. Math. Anal.} {\bf 39} (2007), no. 2, 449--469.
	
\bibitem{SS} {\sc E.~Sandier and S.~Serfaty}. \textit{ Vortices in the Magnetic Ginzburg-Landau Model}. Progress in Nonlinear Differential Equations and Their Applications, vol. 70, Birkh\"auser Boston, 2007.
	%
\bibitem{ScUh} {\sc R.~Schoen and K.~Uhlenbeck}. Boundary regularity and the Dirichlet problem for harmonic maps. \textit{J. Differential Geom.} {\bf 18} (1983), no. 2, 253--268. 

\bibitem{Zhang2020} {\sc X.~Zhang et al.} Skyrmion-electronics: writing, deleting, reading and processing magnetic skyrmions toward spintronic applications. \textit{J. Phys. -- Condensed Matter}  {\bf 32} (2020), no. 14, 143001.

\bibitem{Zu2010} {\sc X.Z.~Yu, Y.~Onose, N.~Kanazawa, J.H.~Park, J.H.~Han, Y.~Matsui, N.~Nagaosa, and Y.~Tokura}. Real-space observation of a two-dimensional skyrmion crystal. \textit{Nature} {\bf 465} (2010), no. 7300, 901--904.


\end{thebibliography}
